\documentclass[brochure,english,12pt]{bourbaki}
\usepackage[matrix,arrow]{xy}
\usepackage{amssymb,amsfonts,amsmath,footnote}
\usepackage{graphicx}
\usepackage[francais]{babel}
\numberwithin{equation}{section} 

\newcommand{\field}[1]{\mathbb{#1}}
\newcommand{\Z}{\field{Z}}
\newcommand{\R}{\field{R}}
\newcommand{\C}{\field{C}}
\newcommand{\N}{\field{N}}

\def\mG{\mathcal{G}}
\def\mH{\mathcal{H}}

\def\mL{\mathcal{L}}
\def\mM{\mathcal{M}}

\def\mP{\mathcal{P}}

\def\kb{\mathfrak{b}}
\def\kg{\mathfrak{g}}
\def\kk{\mathfrak{k}}
\def\kn{\mathfrak{n}}
\def\kp{\mathfrak{p}}
\def\kq{\mathfrak{q}}
\def\kt{\mathfrak{t}}
\def\kz{\mathfrak{z}}

\def\Re{{\rm Re}}
\def\Im{{\rm Im}}

\def\ul{\underline}

\DeclareMathOperator{\ad}{ad}
\DeclareMathOperator{\Ad}{Ad}
\DeclareMathOperator{\End}{End}

\DeclareMathOperator{\Ker}{Ker}

\DeclareMathOperator{\tr}{Tr}

\DeclareMathOperator{\SL}{SL}

\DeclareMathOperator{\vol}{Vol}

\newcommand{\norm}[1]{\lVert#1\rVert}
\newcommand{\abs}[1]{\lvert#1\rvert}

\newcommand{\comment}[1]{}

\date{Mars 2017}
\bbkannee{69\`eme ann\'ee, 2016-2017}
\bbknumero{1130}
\title{Geometric hypoelliptic Laplacian\\ and orbital integrals}
\subtitle{after Bismut, Lebeau and  Shen}
\author{Xiaonan MA}
\address{Universit\'e Paris-Diderot Paris 7\\
Institut de Math\'ematiques de Jussieu--Paris Rive Gauche\\
UMR CNRS 7586\\
B\^atiment Sophie Germain\\
Case 7012\\
F--75205 Paris Cedex 13}
\email{xiaonan.ma@imj-prg.fr}

\begin{document}
\maketitle

\tableofcontents

%%%%%%%%%%%%%%%%%%%%%%%%%%%%%%%%%%
\noindent{\bf INTRODUCTION}\label{s0}

In 1956, Selberg expressed 
the trace of an invariant kernel acting   %reasonable operator 
on a locally symmetric space $Z=\Gamma\backslash G/K$ as 
a sum of certain integrals on the orbits of $\Gamma$ in $G$,
%of the function defining its kernel on orbits, 
the so called ``orbital integrals", and 
he gave a geometric expression for such orbital integrals 
for the heat kernel when $G= \SL_{2}(\R)$, and 
the corresponding locally symmetric space is 
%Selberg %\cite[(3.2)]{Selberg56} 
%he evaluated precisely the trace of the heat operator on 
a compact  Riemann surface of constant negative curvature.
%via the sum of the orbital integrals associated with the  
%homotopy classes of closed geodesics,
In this case, the orbital integrals are one to one correspondence
with the closed geodesics in $Z$.
%which is the original Selberg trace formula.
In the general case, Harish-Chandra worked on the evaluation of 
orbital integrals from the 1950s until the 1970s. He could %formulate 
give an algorithm to reduce the computation of 
an orbital integral to lower dimensional Lie groups by %via
the discrete series method. Given a reductive Lie group, 
in a finite number of steps, 
%theoretically we can always compute it.
there is  %an algorithm for 
a formula for such orbital integrals.
See Section \ref{s2.6} for a brief description of
Harish-Chandra's Plancherel theory.

It is important to understand 
the different properties of  orbital integrals even
without knowing their explicit values.
The orbital integrals appear naturally in Langlands program. 
% To understand the Langlands correspondence for a pair of 
% Langlands dual groups, many questions could be reduced to 
% understand the corresponding question for the orbital integrals.
% Of course, in Langlands program, we work on 
% the reductive Lie group, not only over real number,
% but also on positive characteristic number fields.

About 15 years ago, Bismut gave a natural construction of a
Hodge theory whose corresponding Laplacian 
is a hypoelliptic operator acting on the total space 
of the cotangent bundle of a Riemannian manifold. This operator
interpolates formally between the classical elliptic Laplacian
on the base and the generator of the geodesic flow.
We will describe recent developments in the theory
of the hypoelliptic Laplacian, and %in particular
we will explain two consequences of this program, 
the explicit formula
obtained by Bismut for %semisimple 
orbital integrals, and the recent solution
by Shen of Fried's conjecture (dating back to 1986) 
for locally symmetric spaces.
The conjecture predicts the equality of the analytic torsion and 
of the value at $0$ of the Ruelle dynamical  zeta function
associated with the geodesic flow.

We will describe in more detail these two last results.

Let $G$ be a connected reductive Lie group, let $ \kg$ 
be its Lie algebra,
let $\theta\in {\rm Aut}(G)$ be the Cartan involution of $G$. 
Let $K\subset G$ be the maximal compact subgroup of $G$ 
%that is 
given by the fixed-points of  $\theta$,
and let $\kk$ be its Lie algebra. Let $\kg=\kp\oplus \kk$
be the corresponding Cartan decomposition of $\kg$.

Let $B$ be a nondegenerate bilinear symmetric form on $\kg$
which is invariant under the adjoint action of $G$ on $\kg$ 
and also under $\theta$. 
We assume $B$ is positive on $\kp$ and negative on $\kk$.
Then $\left\langle\cdot, \cdot\right\rangle=-B(\cdot,\theta\cdot)$
is a $K$-invariant scalar product on $\kg$ that is 
such that the Cartan decomposition is an orthogonal splitting.

Let $C^\kg\in U(\kg)$ be the Casimir element of $G$.
If $\{e_i\}_{i=1}^m$ is an orthonormal basis of $\kp$ 
and $\{e_i\}_{i=m+1}^{m+n}$ is an orthonormal basis of $\kk$,
set 
\begin{align} \label{3.11}\begin{split}
&B^*(\kg)= -\frac{1}{2}\sum_{1\leq i,j\leq m} \Big|[e_i,e_j]\Big|^2 
-\frac{1}{6}\sum_{m+1\leq i,j\leq m+n} \Big|[e_i,e_j]\Big|^2,\quad
\mL =\frac{1}{2}C^\kg +\frac{1}{8} B^*(\kg).
\end{split}\end{align}

Let $E$ be a finite dimensional Hermitian vector space, 
let $\rho^E:K\to {\rm U}(E)$ be a unitary representation of $K$. 
Let $F=G\times_K E$ be the corresponding vector bundle over 
the symmetric space $X=G/K$.
Then $\mL  $ descends to 
a second order differential operator  $\mL^X$ acting 
on $C^\infty(X,F)$.
For $t>0$, let $e^{-t\mL^X}(x,x')$ be the smooth kernel of
 the heat operator $e^{-t\mL^X}$.

Assume $\gamma\in G$ is semisimple. 
Then up to conjugation, there exist $a\in \kp, k\in K$ 
such that $\gamma= e^a k^{-1}$ and ${\rm Ad}(k)a=a$.
Let $\tr^{[\gamma]}\left[e^{-t\mL^X}\right]$ denote 
the corresponding orbital integral of  $e^{-t\mL^X}$
(cf. (\ref{eq2.2.2}), (\ref{eq2.3.15})).
If $\gamma=1$, then the orbital integral associated with $1\in G$ 
is given by
 \begin{align} \label{3.14}
\tr^{[\gamma=1]}\left[e^{-t\mL^X}\right]
=\tr^F\left[e^{-t\mL^X}(x,x)\right] \, \, 
%\text{Êfor any } x\in X.
\end{align}
which does not depend on $x\in X$.

Let $Z(\gamma)\subset G$ be the centralizer of $\gamma$, 
and let $\kz(\gamma)$ be its Lie algebra. Set 
$ \kp(\gamma)= \kz(\gamma)\cap \kp$, 
$\kk(\gamma)= \kz(\gamma)\cap \kk$.
Then $\kz(\gamma)= \kp(\gamma)\oplus\kk(\gamma)$.

Set $\kz_{0}= \Ker (\ad(a)), \kk_{0}=\kz_{0}\cap \kk$. 
Let $\kz_{0}^{\bot}$ be the orthogonal space to 
$\kz_{0}$ in $\kg$. Let $\kk_{0}^{\bot}(\gamma)$ be the 
orthogonal space to $\kk(\gamma)$ in $\kk_{0}$, and 
$\kz_{0}^{\bot}(\gamma)$ be the orthogonal space to 
$\kz(\gamma)$ in $\kz_{0}$, so that
$\kz_{0}^{\bot}(\gamma)
= \kp_{0}^{\bot}(\gamma)\oplus\kk_{0}^{\bot}(\gamma)$.
For a self-adjoint matrix $\Theta$,
set $\widehat{A}(\Theta)= \det^{1/2}
\Big[\frac{\Theta/2}{\sinh(\Theta/2)}\Big]$.
For $Y\in \kk(\gamma)$, set 
\begin{multline}\label{3.15}
J_{\gamma}(Y)= \left|\det (1-\Ad(\gamma))
|_{\kz_{0}^{\bot}} \right|^{-1/2}
\frac{\widehat{A}(i\ad(Y)|_{\kp(\gamma)})}
{\widehat{A}(i\ad(Y)|_{\kk(\gamma)})}\\
\times \left[    \frac{1}{\det 
(1-{\Ad}(k^{-1}))|_{\kz^{\bot}_{0}(\gamma)}}
\frac{\det \Big(1-e^{-i\ad(Y)}\Ad(k^{-1})
\Big)|_{\kk^{\bot}_{0}(\gamma)}}
    {\det\Big(1-e^{-i\ad(Y)}\Ad(k^{-1})
\Big)|_{\kp^{\bot}_{0}(\gamma)}}
\right]^{1/2}.
\end{multline}
If $\gamma=1$, then the above equation reduces to
$J_{1}(Y)=\frac{\widehat{A}(i\ad(Y)|_{\kp})}
{\widehat{A}(i\ad(Y)|_{\kk})}$ for $Y\in \kk=\kk(1)$.

\begin{theo}\label{t0.1} {\rm(Bismut's orbital integral formula 
\cite[Theorem 6.1.1]{B11b})} Assume $\gamma\in G$ is 
semisimple.Then for any $t>0$, we have
 \begin{multline} \label{3.17}
\tr^{[\gamma]}\left[e^{-t\mL^X}\right]
= (2\pi t)^{-\dim \kp(\gamma)/2} e^{-\frac{|a|^{2}}{2t}}\\
\int_{\kk(\gamma)} J_{\gamma}(Y)
\tr^{E}\left[\rho^{E}(k^{-1}) e^{-i\rho^{E}(Y)}\right]
e^{-\frac{|Y|^{2}}{2t}} 
\frac{dY}{(2\pi t)^{\dim \kk(\gamma)/2}}.
\end{multline}
\end{theo}

There are some striking similarities of Equation (\ref{3.17})
%Theorem \ref{t0.1} 
with the Atiyah-Singer index formula, where 
the $\widehat{A}$-genus of the tangent bundle appears.
%as the contribution of the manifold. Here, 
%even we do not compute 
%an integer, such as the index of an elliptic operator, 
Here the $\widehat{A}$-function of both $\kp$ and $\kk$ parts 
(with different roles) appear naturally in the integral (\ref{3.17}). 

A more refined version %reformulation 
of Theorem \ref{t0.1} for the orbital
integral associated with the wave operator is given in 
\cite[Theorem 6.3.2]{B11b} (cf. Theorem \ref{t2.10}).

Let $\Gamma\subset G$ be a discrete cocompact torsion free 
subgroup.
The above objets constructed  on $X$ descend to 
the locally symmetric space
$Z=\Gamma\backslash X$ and $\pi_{1}(Z)=\Gamma$.
We denote by $\mL^Z$ the
corresponding differential  operator on $Z$.
% Note that for $\gamma\in \Gamma$,
% $X(\gamma) = Z(\gamma)/K\cap Z(\gamma)$
% is a symmetric space.
 Let $[\Gamma]$ be the set of conjugacy classes in $\Gamma$. 
The Selberg trace formula (cf. (\ref{eq2.22}), (\ref{eq2.4.10})) for 
the heat kernel of the Casimir operator on $Z$ says that
\begin{equation}\label{0.5}
\tr[e^{-t\mL^Z}]=\sum_{[\gamma]\in [\Gamma]} 
\vol \Big(\Gamma\cap Z(\gamma)\backslash 
Z(\gamma) \Big) 
\tr^{[\gamma]}[e^{-t \mL^X}].
\end{equation}
Each term $\tr^{[\gamma]}[\cdot]$ in (\ref{0.5}) is evaluated in 
\eqref{3.17}.

Assume $m= \dim X$ is odd now.
Let $\rho:\Gamma\to U({\bf q})$ be a unitary representation. 
Then $F=X\times_{\Gamma} \C^{\bf q}$ is a flat Hermitian vector bundle 
on $Z=\Gamma\backslash X$.  Let $T(F)$ be the analytic torsion 
associated with $F$ on $Z$ (cf. Definition \ref{d3.1}), 
which is a regularized determinant of the Hodge Laplacian 
for the de Rham complex associated with $F$. 

In 1986, Fried discovered a surprising relation of 
the analytic torsion to dynamical systems. In particular, 
for a compact orientable hyperbolic manifold, he identified 
the value at zero of the Ruelle dynamical zeta function
associated with the closed geodesics in $Z$ and with $\rho$, 
to the corresponding analytic torsion, and he conjectured that 
a similar result should hold for general compact 
locally homogenous manifolds. 
In 1991, Moscovici-Stanton \cite{MS91} made an important 
progress in the proof of 
Fried's conjecture for locally symmetric spaces. 
The following recent result of Shen establishes Fried's conjecture
for arbitrary  locally symmetric spaces, and Theorem \ref{t0.1}
is one important ingredient in Shen's proof.

 Given $[\gamma]\in [\Gamma]\backslash \{1\}$, let
 $B_{[\gamma]}$  be the space of closed geodesics in $Z$ 
 which lie in the homotopy class $[\gamma]$,
 and let $l_{[\gamma]}$ be the length of the geodesic
associated with $\gamma$ in $Z$.  
The group $\mathbb{S}^1$ acts on 
$B_{[\gamma]}$ by rotations. This action is locally free. 
Denote by $\chi_{\mathrm{orb}}(\mathbb{S}^1\backslash 
B_{[\gamma]})\in \mathbb{Q}$ the orbifold Euler characteristic
number for the quotient orbifold 
$\mathbb{S}^1\backslash B_{[\gamma]}$. Let
\begin{align}\label{4.3.5}
n_{[\gamma]}=\left|\Ker\big(\mathbb{S}^1\to 
\mathrm{Diff}(B_{[\gamma]})\big)\right|
\end{align}
be the generic multiplicity of $B_{[\gamma]}$.
\begin{theo}
\label{t0.2}\cite{Shen16}
%If $F$ is an acyclic flat vector bundle on $Z$ with holonomy 
For any unitary representation $\rho:\Gamma\to U({\bf q})$, 
\begin{align}
 \label{4.3.6}
R_\rho(\sigma)=  \exp\left(
\sum_{[\gamma]\in [\Gamma]\backslash \{1\}}
 \tr[\rho(\gamma)]\frac{\chi_{\mathrm{orb}}
 (\mathbb{S}^1\backslash B_{[\gamma]})}{n_{[\gamma]}}
 e^{-\sigma l_{[\gamma]}}\right)
 \end{align}
is a well-defined meromorphic function on $\C$.
If  $H^\bullet(Z,F)=0$, then $R_\rho(\sigma)$
is holomorphic at $\sigma=0$ and 
\begin{align}
\label{4.3.9}
R_\rho(0)=T(F)^2.
\end{align}
\end{theo}

This article is organized as follows.
In Section \ref{s00}, we describe %introduce  
Bismut's program on the geometric hypoelliptic Laplacian
in de Rham theory, and we give its applications. 
In Section \ref{s1}, we introduce the heat kernel on 
smooth manifolds and the basic ideas in the heat equation proof 
of the Lefschetz fixed-point formulas,
which will serve as a model for the proof of Theorem \ref{t0.1}.
In Section \ref{s2}, we review %the definition of the 
orbital integrals, their relation
to  Selberg trace formula,
and  we state Theorem \ref{t0.1}.
In Section \ref{s3}, we give the basic ideas in how 
to adapt the  construction of the hypoelliptic Laplacian of 
Section \ref{s00}
in the context of locally symmetric spaces in order to 
establish Theorem \ref{t0.1}. 
In Section \ref{s4}, 
we concentrate on Shen's solution of Fried's conjecture. 
%on the equality of the analytic torsion and 
%of the values at $0$ of the dynamical zeta function.

{\bf Notation}: If $A$ is a $\Z_{2}$-graded algebra, if $a,b\in A$,
the supercommutator $[a,b]$
is given by
\begin{align}\label{0.7}
[a,b]= ab -(-1)^{\deg a\cdot \deg b} ba.
\end{align}
If $B$ is another $\Z_{2}$-graded algebra, we denote by 
$A\widehat{\otimes} B$ the $\Z_{2}$-graded tensor product,
such that the $\Z_{2}$-degree 
of $a\widehat{\otimes} b$ is given by 
$\deg a +\deg b$, and
where the product is given by
\begin{align}\label{0.8}
(a\widehat{\otimes} b)\cdot (c\widehat{\otimes} d)
= (-1)^{\deg b\cdot \deg c} ac\, \widehat{\otimes} \, bd.
\end{align}
If $E=E^+\oplus E^-$ is a $\Z_{2}$-graded vector space,
and $\tau=\pm 1$ on $E^\pm$, 
for $u\in \End(E)$, the supertrace $\tr_{s}[u]$
is given by
\begin{align}\label{0.9}
\tr_{s}[u]= \tr[\tau u].
\end{align}

In what follows,
we will often add a superscript to indicate
where the trace or supertrace is taken.
%of this vector space on the trace symbol . 

\comment{  For $\beta$ a number or a matrix, we denote
\begin{equation}\label{0.10}
\sinh(\beta)=\frac{1}{2}(e^\beta-e^{-\beta}),\;\cosh(\beta)
=\frac{1}{2}(e^\beta+e^{-\beta}),\; \tanh(\beta)
=\frac{\sinh(\beta)}{\cosh(\beta)}.
\end{equation}
}

\textbf{\emph{Acknowledgments.}}
I thank Professor Jean-Michel Bismut very heartily for his help 
and advice during the preparation of this manuscript. 
It is a pleasure to thank Laurent Clozel, Bingxiao Liu, 
George Marinescu and Shu Shen for their help and remarks.

%%%%%%%%%%%%%%%%%%%%%%%%%%%%%
\section{From hypoelliptic 
Laplacians to the trace formula}\label{s00}
In this section, we describe some basic ideas %of
taken from Bismut's program 
on the geometric hypoelliptic Laplacian and its applications 
to geometry and dynamical systems.

A differential operator $P$ is hypoelliptic if 
for every distribution $u$ defined on an open set $U$
 such that $Pu$ is smooth, then $u$ is smooth on $U$.
Elliptic operators are hypoelliptic, but there are 
hypoelliptic differential operators which are not elliptic.
Classical examples are Kolmogorov operator 
$\frac{\partial ^2}{\partial y^2}- y \frac{\partial}{\partial x}$
on $\R^2$ \cite{K34} and H\"ormander's generalization
$\sum_{j=1}^k X_{j}^2 + X_{0}$  
on Euclidean spaces \cite{Hormander67}.
Along this line, see for example 
Helffer-Nier's \cite{HelfferNier05} recent book
and Lebeau's work \cite{Lebeau05} on the hypoelliptic estimates
and Fokker-Planck operators.

In 1978, Malliavin \cite{Malliavin78}  introduced the so-called 
`Malliavin calculus' to reprove
H\"ormander's regularity result  \cite{Hormander67}
from a probabilistic  point of view.
Malliavin calculus was further %completed
developed by Bismut \cite{B81d}
and Stroock \cite{Stroock81}. 

%Hypoelliptic operator is well know for analysts, but not for 
%geometers. 
About 15 years ago, Bismut initiated a program whose purpose 
is to study the applications of hypoelliptic second order
differential operators to differential geometry.

In \cite{B05a},
Bismut constructed a (geometric)  hypoelliptic Laplacian on 
the total space of the cotangent bundle 
$T^*M$ of a compact Riemannian manifold $M$,
that depends on a parameter $b>0$.
This hypoelliptic Laplacian %on $T^*M$ 
is a deformation of the usual Laplacian on $M$.
More precisely, when $b\to 0$, it %this hypoelliptic Laplacian
converges to the Laplacian on $M$ in a suitable sense,
and when $b\to +\infty$, it converges to the generator of 
the geodesic flow. %on $T^*M$. 
In this way, %dynamical system 
properties of the geodesic flow on $M$ 
are potentially related to
the spectral properties of the Laplacian on $M$.

We now explain briefly %what is
Bismut's %(geometric) 
hypoelliptic Laplacian in de Rham theory.
Let $(M, g^ {TM})$ be a compact Riemannian manifold 
of dimension $m$.
Let $(\Omega^\bullet (M), d)$ be the de Rham complex of $M$,
let $d^*$ be the formal $L_{2}$ adjoint of $d$, and 
let $\square^M= (d+d^{*})^2$ be the Hodge Laplacian 
acting on $\Omega^\bullet (M)$.

 Let $\pi: \mM \to M$ 
 be the total space of the cotangent bundle $T^{*}M$.
%, which depends on a parameter $b>0$. 
 %Let be the natural projection.
Let $\Delta^{V}$ be the  Laplacian along the fibers 
$T^{*}M$, and let $\mH$ be the function on $\mM$ defined by
\begin{align}\label{3.1}
\mH(x,p)= \frac{1}{2}\, |p|^{2} 
\quad \text{ for } p\in T_{x}^{*}M, x\in M.
\end{align}
Let $Y^{\mH}$ be the Hamiltonian vector field on $\mM$ 
associated with $\mH$ and  with the canonical symplectic form
on $\mM$.
Then $Y^{\mH}$ is the generator of the geodesic flow. %on $\mM$.
Let $L_{Y^{\mH}}$ denote the corresponding 
Lie derivative operator acting on $\Omega^\bullet (\mM)$.
For $b>0$, the Bismut hypoelliptic Laplacian on $\mM$
is given by
%on the total space of $T^*M$ has form:
\begin{align}\label{3.2}
\mL_b =\frac{1}{b^2} \alpha + \frac{1}{b} \beta + \vartheta,
\end{align}
with 
\begin{align}\label{3.3}\begin{split}
& \alpha = \frac{1}{2} (-\Delta^V + |p|^2 -m +\cdots),\quad
  \beta = -L_{Y^{\mH}} + \cdots, 
\end{split}\end{align}
where the dots and  $\vartheta$ are geometric terms 
which we will not be made explicit. 
%precise here. Thus basically, Essentially
The operator $\mL_b$ is essentially the weighted sum of 
the harmonic oscillator along the fiber,
 minus the generator of the geodesic flow 
 $- L_{Y^{\mH}}$ along the horizontal 
 direction.\footnote{On Euclidean spaces, all geometric terms
 vanish and the operator $\mL_b$ acting on functions  reduces to
 the Fokker-Planck operator.}

The vector space $\Ker (\alpha)$ is spanned by the function 
$\exp(-|p|^2/2)$.
We identify $\Omega^\bullet (M)$ to $\Ker (\alpha)$ by the map 
$s \to \pi^* s \exp(-|p|^2/2)/ \pi^{m/4}$. 
Let $P$ be the standard $L_{2}$-projector from 
$\Omega^\bullet (\mM)$ on $\Ker (\alpha)$.
Then  by \cite[Theorem 3.14]{B05a}, 
\begin{align}\label{3.4}
P(\vartheta - \beta\alpha^{-1}\beta)P
= \frac{1}{2}\square^M.
\end{align}
In \cite{B05a}, equation \eqref{3.4}  is used to prove that 
%gives strong indication that $\mL_b$ is a deformation of 
%$\frac{1}{2}\square^M$ in the sense that
as $b\to 0$,
we have the formal convergence of resolvents 
\begin{align}\label{3.4a}
(\lambda - \mL_b)^{-1}
\to P \Big(\lambda - \frac{1}{2}\square^M\Big)^{-1}P.
\end{align}

Bismut-Lebeau \cite{BL08} set up the proper analysis foundation
for the study of the hypoelliptic Laplacian $\mL_b$.
%developed, among other important results, 
%the analytic details  of the above interpolation when  $b\to 0$.
They not only proved a corresponding version of 
the Hodge theorem,
but they also studied the precise properties of its resolvent 
and of the corresponding heat kernel. 
Since $\mM$ is noncompact, they needed to refine %for example
the hypoelliptic estimates of 
H\"{o}rmander in order to control hypoellipticity at infinity.
%in the cotangent bundle. 
They developed the adequate theory 
of semiclassical pseudodifferential operators with parameter
$\hbar= b$ and obtained the proper version of the convergence
 of resolvents in (\ref{3.4a}).
%the convergence as $b\to 0$ of $\mL_b$ to 
%$\frac{1}{2}\square^M$ by studying the convergence
% of the resolvent (\ref{3.4a}) in a proper sense.
 They developed also a hypoelliptic local index theory
 which is itself a deformation of classical elliptic local index theory.
 
In  \cite{BL08}, Bismut-Lebeau defined a hypoelliptic version of 
the analytic torsion of  Ray-Singer \cite{RS71} associated with 
the elliptic Hodge Laplacian  in (\ref{3.4}). 
 %and of the analytic torsion form of Bismut-Lott \cite{BLo95}. 
%  They compared explicitly  their 
%  hypoelliptic torsion form with the Bismut-Lott torsion fom,
%  and extended it to the group action case.
 %In particular, they proved 
 The main result in  \cite{BL08} is the proof of 
 the equality of the hypoelliptic
 torsion with the Ray-Singer analytic torsion.
% Note that one of basic difficulties is that the hypoelliptic Laplacian 
% is not self-adjoint. Moreover, they developed also the 
%  appropriate local index theory  for the associated heat kernel, 
%  and give an explicit formula relating the analytic torsion
%  objects associated to the hypoelliptic Laplacian
% to the classical Ray-Singer analytic torsion.

In his thesis \cite{Shen16a}, Shen studied the Witten deformation 
of the hypoelliptic Laplacian for a Morse function 
on the base manifold, and  identified the hypoelliptic 
 torsion to the combinatory torsion. 
%In particular, combined with Bismut-Zhang's result \cite{BZ92}, 
Shen's work gives a new proof 
of Bismut-Lebeau's result on the equality of the hypoelliptic 
 torsion and the Ray-Singer analytic  torsion. 

%Beside the results on the orbital integrals I will review here, 
 This article concentrates on applications of the hypoelliptic 
 Laplacian to orbital integrals. We will briefly summarize other 
 applications.
%  the readers will find a few more 
% applications of the geometric hypoelliptic Laplacians.
 
A  version of Theorem \ref{t0.1} for compact Lie groups 
can be found in \cite{B08a}.
In \cite[Theorem 4.3]{B08a}, as a test of his ideas, 
Bismut gave a new proof of 
the classical explicit formula for the scalar heat 
kernel in terms of the coroots lattice
\cite{Frenkel84} for a simple simply connected compact Lie group, 
by using the hypoelliptic Laplacian on the total 
space of the cotangent bundle of the group.
% Note that in \cite{Frenkel84}, using Poisson summation,
% Frenkel showed that the Kac character formula for finite energy 
% representations of affine Lie algebras can be reexpressed in terms
% of the scalar heat kernel on the group.
In \cite{B08c}, Bismut also constructed a hypoelliptic 
Dirac operator which is a hypoelliptic deformation of the usual 
Dirac operator.

In \cite[Theorem 0.1]{B13}, Bismut established a 
Grothendieck-Riemann-Roch theorem for a proper holomorphic 
submersion $\pi: M\to B$ of complex manifolds in %the 
Bott-Chern cohomology. 
For compact K\"{a}hler manifolds, Bott-Chern cohomology
coincides with de Rham cohomology.
In the general situation considered in \cite{B13}, 
the elliptic methods of \cite{B86}, \cite{BGS88b} are known to fail,
and hypoelliptic methods seem to be 
the only way to obtain this result.
% Note that if the base manifold $B$ is K\"{a}hler,  it is 
% a consequence of the Atiyah-Singer 
% family index theorem which holds in the de Rham cohomology,
% but if $B$ is not K\"{a}hler, it seems that the usual local family 
% index technique \cite{B86} fails in this situation, thus Bismut 
%had to use his hypoelliptic Laplacian formalism.

%We would like to notice that 
As in the case of the Dirac operator, there does not 
exist a universal 
hypoelliptic Laplacian which works for all situations, %problems.
there are several hypoelliptic Laplacians. To attack a specific 
(geometric) problem, we need to construct %or find 
the corresponding hypoelliptic Laplacian.
Still all the hypoelliptic Laplacians have naturally the same 
structure, but the geometric terms depend on the situation. 
% Even if the harmonic 
% oscillator along the fiber and the generator of the geodesic flow 
% iappears always in the final form, the other geometric 
% terms depend really on the problems at hand %we are working. 
% For example, to establish Theorem \ref{t0.1}, already the base 
% manifold is a noncompactsymmetric space,
% so we have to enlarge  the total space of the cotangent bundle,
% and in \cite{B11b}, \cite{B13}, there are also positive polynomials
% of degree $4$ appearing along the fiber.
%Another important point is that  
Probability theory plays an important role, both 
formally and technically in its construction and in its use.
% at the formal aspect and at the technical level. Sometimes,
%it is `easier' and natural to get the estimates from Malliavin calculus, 
%but hard to prove them by  the classical analysis (including  
%functional analysis) techniques for the hypoelliptic Laplacian.

In this article,  we will not touch the analytic and 
probabilistic aspects of the proofs.
We will explain how to give a natural construction of the 
%essentially how to construct the natural 
hypoelliptic Laplacian which is needed in order to establish 
Theorem \ref{t0.1}. 
The method consists in giving a cohomological interpretation
to orbital integrals, so as to reduce 
their evaluation to methods related to the proof of 
Lefschetz fixed-point formulas.
Theorem \ref{t0.1} gives a direct link of
 the trace formula to  index theory.
% We will emphasize its cohomological aspect: how to view it as a 
% Lefschetz fixed-point formula, both as formula and as proof.
% Theorem \ref{t0.1} gives a direct link between 
% the trace formula and the index theory,  so that we can now 
% roughly say that the trace formula is an index type theorem.

%We invite the readers to read 
We hope this article can be used as an invitation to the original  
papers \cite{B05a, B08a, B08c, B11b, B13, B16b}
and to several surveys on this topic
\cite{B08b, B08d, B11a, B12, B15, B16} and \cite{Lebeau08}.

%\newpage 

%%%%%%%%%%%%%%%%%%%%%%%%%%%%%
\section{Heat kernel and Lefschetz fixed-point formula}\label{s1}

This section is organized as follows. In Section \ref{s1.1}, 
we explain some basic facts about heat kernels.
%of Laplace operators.
In Section \ref{s1.2}, we review the heat  equation 
proof of the Lefschetz fixed-point formula. 
This proof will be %served 
used as a model for the proof of the main theorem of this article. 
% We try to emphasize the similarity in the spirit of the two proofs, 
% and how it helps to find the right approach.

%%%%%%%%%%%%%%%%%%%%%%%%%%%%%%%%%%
\subsection{A brief introduction to the heat kernel}\label{s1.1}
Let $M$ be a compact manifold of dimension $m$. 
Let $TM$ be the tangent bundle, $T^*M$ be the cotangent 
bundle, and let $g^{TM}$ be a Riemannian metric on 
%the tangent bundle $TM$ of 
$M$. Let $F$ be a complex vector bundle over $M$, and 
let $h^F$ be a Hermitian metric on $F$. 
Let $C^\infty(M,F)$ be the space of smooth sections 
of $F$ on $M$. Let $\left\langle \cdot ,\cdot\right\rangle$
be the $L_{2}$-Hermitian product on $C^\infty(M,F)$ 
defined by the integral of the pointwise product 
with respect to the Riemannian volume form $dx$. We denote by 
$L_{2}(M,F)$ the vector space of $L_{2}$-integrable sections 
of $F$ on $M$. %and naturally $C^\infty(M,F)\subset L_{2}(M,F)$. 

Let $\nabla^F: C^\infty(M,F)\to C^\infty(M,T^*M\otimes F)$ be 
a Hermitian connection on $(F, h^F)$ and let $\nabla^{F,*}$ 
be its formal adjoint. 
Then the (negative) Bochner Laplacian $\Delta^F$ 
acting on $C^\infty(M,F)$, is defined by
\begin{equation}
\label{eq1.1}
-\Delta^F=\nabla^{F,*}\nabla^F.
\end{equation}

The operator $-\Delta^F$ is 
an essentially self-adjoint second order elliptic operator. 
Let $\nabla^{TM}$ be the Levi-Civita connection on 
$(TM,g^{TM})$. We can rewrite it as
\begin{equation}
\label{eq1.2}
-\Delta^F=-\sum^m_{i=1} \Big((\nabla_{e_i}^{F})^2
-\nabla^F_{\nabla^{TM}_{e_i}e_i}\Big),
\end{equation} 
where $\{e_i\}^m_{i=1}$ is a local smooth orthonormal frame
of $(TM, g^{TM})$.

%Hermitian 
For a self-adjoint section $\Phi\in C^\infty(M,\text{End}(F))$
(for any $x\in M$ %which means
that  $\Phi_x\in \text{End}(F_x)$ 
is self-adjoint), set
%We define a self-adjoint differential operator $-\Delta^F_\Phi$ by
\begin{equation}
\label{eq1.3}
-\Delta^F_\Phi=-\Delta^F-\Phi.
\end{equation} 
Then the heat operator $e^{t \Delta^F_\Phi}: L_{2}(M,F)
\rightarrow L_{2}(M,F)$ for $t>0$ of $-\Delta^F_\Phi$ 
is the unique solution of 
\begin{equation}
\label{eq1.4}
\left\{
\begin{matrix}
&\big(\frac{\partial}{\partial t}-\Delta^F_\Phi \big )
e^{t\Delta^F_\Phi}=0&\\ 
&\lim_{t\rightarrow 0} 
e^{t\Delta^F_\Phi}s=s \in L_{2}(M,F)& 
\text{ for any } s\in L_{2}(M,F).
\end{matrix}
\right.
\end{equation}
For $x,x'\in M$, let $e^{t\Delta^F_\Phi}(x,x')\in F_x\otimes F^*_{x'}$ 
be the Schwartz kernel of the operator $e^{t\Delta^F_\Phi}$
with respect to the Riemannian volume element $dx'$. 
Classically, $e^{t\Delta^F_\Phi}$ is smooth in $x,x'\in M, t>0$. 

Since $M$ is compact, the operator $-\Delta^F_\Phi$ has 
discrete spectrum, consisting of eigenvalues
$\lambda_1\leq \lambda_2\leq \cdots
\leq \lambda_k\leq \cdots$ counted with multiplicities, 
with $\lambda_k\rightarrow +\infty$
as $k\rightarrow +\infty$. Let $\{\varphi_j\}^{+\infty}_{j=1}$ 
be a system %the set 
of orthonormal eigenfunctions such that 
$-\Delta^F_\Phi \varphi_j=\lambda_j \varphi_j$. Then
$\{\varphi_j\}_{j=1}^{+\infty}$ is an %complete 
orthonormal basis 
of $L_{2}(M,F)$. The heat kernel can also be written as 
(cf. \cite[Proposition 2.36]{BeGeVe04}, \cite[Appendix D]{MM07})
\begin{equation}
\label{eq1.5}
e^{t\Delta^F_\Phi}(x,x')=\sum^{+\infty}_{j=1} e^{-t\lambda_j}
\varphi_j(x)\otimes \varphi_j(x')^*
\end{equation} 
%for $t>0$, 
where $\varphi_j(x')^*\in F^*_{x'}$ is 
the metric dual of $\varphi_j(x')\in F_{x'}$.

The trace of the heat operator is given by
\begin{equation}
\label{eq1.7}
\tr[e^{t\Delta^F_\Phi}]=\sum^{+\infty}_{j=1} e^{-t\lambda_j}.
\end{equation}
The (heat) trace $\tr[e^{t\Delta^F_\Phi}]$ involves 
the full spectrum information of  operator $\Delta^F_\Phi$ 
and has many applications.
%for example, it appears naturally in the heat equation %kernel 
%approach to the Atiyah-Singer index theorem. 
%so called, the local index theorem, 
%It also relates to the secondary spectral invariants such 
%as the analytic torsion and the eta-invariant.

In general, it is difficult to evaluate explicitly 
$\tr[e^{t\Delta^F_\Phi}]$ for $t>0$. 
However, we will explain the explicit formula obtained by 
Bismut for locally symmetric spaces and
its connection with Selberg trace formula. 

\begin{rema}\label{r1.1}
Let $\pi: \widetilde{M}\rightarrow M$ be the universal cover of $M$
with fiber $\pi_1(M)$, the fundamental group of $M$. 
Then geometric data on $M$ lift to %on 
$\widetilde{M}$, and we will %simply 
add a $\;\widetilde{}\;$ 
to denote the corresponding objets on $\widetilde{M}$. 
%We put $\widetilde{M}$ the canonical projection , 
%By a classical result, 
It's well-known  (see for instance \cite[(3.18)]{MM15}) that
if $\widetilde{x}, \widetilde{x}'\in \widetilde{M}$ are such that
$\pi (\widetilde{x})=x, \pi (\widetilde{x}')=x'$, we have
\begin{equation}
\label{eq1.6}
e^{t\Delta^F_\Phi}(x,x')=\sum_{\gamma\in\pi_1(M)}\gamma 
e^{t\widetilde{\Delta}^F_\Phi}
(\gamma^{-1}\widetilde{x},\widetilde{x}'),
\end{equation}
%which shows in particular that 
where the right-hand side is uniformly convergent.
\end{rema}

%%%%%%%%%%%%%%%%%%%%%%%%%%%%%%%%%%
\subsection{The Lefschetz fixed-point formulas}\label{s1.2}

Let $\Omega^\bullet(M)=\oplus_j \Omega^j(M)
=\oplus_j C^\infty(M, \Lambda^j(T^*M)) $ be the vector space 
of smooth differential forms on $M$ (with values in $\R$), which is 
$\Z$-graded by degree. %their degrees. 
Let $d:\Omega^j(M)\rightarrow \Omega^{j+1}(M)$ 
be the exterior differential operator. Then $d^2=0$
so that $(\Omega^\bullet(M),d)$ forms %a complex, the so-called 
the de Rham complex. The %($j$th) 
de Rham cohomology groups of $M$ are defined by %as
\begin{equation}
\label{eq3.1.1}
H^j(M,\R)=\frac{\Ker(d|_{\Omega^j(M)})}
{\Im (d|_{\Omega^{j-1}(M)})},\quad
H^\bullet (M,\R)= \bigoplus_{j=0}^m H^j(M,\R).
\end{equation}
They are canonically isomorphic to the singular cohomology %group
of $M$. %by the de Rham theorem.

Let $d^*: \Omega^\bullet(M)\rightarrow \Omega^{\bullet-1}(M)$
be the %corresponding 
formal adjoint of $d$ with respect to the scalar product 
$\left\langle \cdot , \cdot \right\rangle$ on $\Omega^\bullet(M)$,
i.e., for all $s,s'\in\Omega^\bullet(M)$,
\begin{equation}
\label{eq3.1.2}
\langle d^*s,s'\rangle :=\langle s, ds'\rangle.
\end{equation}

Set 
\begin{equation}
\label{eq3.1.3}
D=d+d^*.
\end{equation}
Then $D$ is a first order elliptic differential operator, 
%and as $(d^*)^2=0$, 
and we have
\begin{equation}
\label{eq3.1.4}
D^2=dd^*+d^*d.
\end{equation}
The operator $D^2$ is called the Hodge Laplacian, 
it is an operator of the type \eqref{eq1.3}
%which is a Laplace type operator as in \eqref{eq1.3}
for $F=\Lambda^\bullet(T^*M)$, which preserves the $\Z$-grading 
on $\Omega^\bullet(M)$. By Hodge theory, we have 
the isomorphism,
\begin{equation}
\label{eq3.1.5}
\Ker (D|_{\Omega^j(M)})=\Ker (D^2|_{\Omega^j(M)})
\simeq H^j(M,\R),\;\mathrm{for}\; j=0,1,\cdots,m.
\end{equation}

We give here a baby example to explain %what is 
the heat equation
proof of the Atiyah-Singer index theorem (cf. \cite{BeGeVe04}).

Let $H$ be a compact Lie group acting on $M$ on the left. 
Since the exterior differential commutes with the %induced group
action of $H$ on $\Omega^\bullet(M)$, %we know that 
$H$ acts naturally on $H^j(M,\R)$ for any $j$.
The Lefschetz number for $h\in H$
is given by
\begin{equation}
\label{eq1.10}
\chi_h(M) = \sum^m_{j=0} (-1)^j \tr[h|_{H^j(M,\R)}]
= \tr_{s}[h|_{H^\bullet(M,\R)}].
\end{equation}
The Lefschetz fixed-point formula computes %simply 
$\chi_h(M)$ in term of %via the 
geometric data on the fixed-point set of $h$.

Instead of working on $H^j(M,\R)$, we will work on 
the much larger space $\Omega^\bullet(M)$ to establish 
the Lefschetz fixed-point formulas.

Since $H$ is compact, by an averaging argument on $H$, 
we can assume that the metric $g^{TM}$ is $H$-invariant. 
Then the operator $D$ defined above is also $H$-invariant.
We have the following result 
(cf. \cite[Theorem 3.50, Proposition 6.3]{BeGeVe04}),
\begin{theo}[McKean-Singer formula]
\label{t1.1}
For any $t>0$,
	\begin{equation}
	\label{eq1.11}
	 \chi_h(M)=\tr_s [he^{-tD^2}].
	\end{equation}
%Here $\tr_s=\tr|_{\Omega^\mathrm{even}}
%-\tr|_{\Omega^\mathrm{odd}}$.
\end{theo}
\begin{proof}
For any $t>0$, we have
\begin{equation}\label{eq1.12}
	\begin{split}
\frac{\partial}{\partial t}\tr_s[he^{-tD^2}]&=
-\tr_s[hD^2 e^{-t D^2}]\\
&=-\frac{1}{2}\tr_s[[D,hDe^{-tD^2}]]=0.
	\end{split}
\end{equation}
Here $[\cdot,\cdot]$ is a supercommutator
defined as in (\ref{0.7}),
and as in the case of matrices, the supertrace of 
a supercommutator vanishes by a simple algebraic argument. 
	
By \eqref{eq1.7} and (\ref{eq3.1.5}), we have
\begin{equation}
	\label{eq1.13}
\lim_{t\rightarrow +\infty} \tr_s[he^{-tD^2}]=\chi_h(M).
\end{equation}
Combining \eqref{eq1.12} and \eqref{eq1.13}, we get \eqref{eq1.11}.
\end{proof}

A simple %careful 
analysis shows that only the  fixed-points of $h$
contribute to the limit of $\tr_s[he^{-tD^2}]$ as 
$t\rightarrow 0$. %by a localization phenomenon. This yields
Further simple work then leads to  the Lefschetz fixed-point formulas.

% Now by studying carefully the limit of $\tr_s[he^{-tD^2}]$ as 
% $t\rightarrow 0$, which will be localized on the fixed-point sets
% of $h$, we get the Lefschetz fixed-point formula.

Even though we will work on a more refined object
the trace of a heat operator, %we will explain that 
the above philosophy still applies.

%%%%%%%%%%%%%%%%%%%%%%%%%%%%%%%%
\section{Bismut's explicit formula for the orbital integrals}\label{s2}

In this section, we give an introduction to orbital integrals
and to Selberg trace formula, and we present the main result of this 
article: Bismut's explicit evaluation of the orbital integrals.
% For the set-up of the problem and Bismut's proof,
% the only fact on Lie groups we needed is the Cartan 
% decomposition $G= e^{\kp} K$ induced by the Cartan involution 
% $\theta$ on the group $G$. 
Also, we compare Harish-Chandra's Plancherel theory with 
Bismut's explicit formula for the orbital integrals.

This section is organized as follows. In Section \ref{s2.1},
we recall some basic facts on symmetric spaces, and
we explain how the Casimir operator for a reductive Lie group
induces a Bochner Laplacian on the associated symmetric space. 
In Section \ref{s2.2}, we give an introduction to
%an intuitive motivation for the introduction of %to introduce
 orbital integrals and to Selberg trace formula,
and in Section \ref{s2.3}, 
we describe the geometric definition of  
 orbital integrals given by Bismut.
In Section \ref{s2.4}, 
we present the main result of this article, 
Bismut's explicit evaluation of the orbital integrals,
%Theorem \ref{t2.7}, 
and give some examples.
Finally in Section \ref{s2.6}, we present briefly Harish-Chandra's 
Plancherel theory for comparison with Bismut's result.

%%%%%%%%%%%%%%%%%%%%%%%%%%%%%%%%
\subsection{Casimir operator and Bochner Laplacian}\label{s2.1}

Let $G$ be a connected real reductive Lie group with Lie algebra 
$\kg$ and Lie bracket $[\cdot,\cdot]$. 
Let $\theta \in \text{Aut}(G)$ be its Cartan involution. 
Let $K$ be the subgroup of $G$ fixed by $\theta$, 
with Lie algebra $\kk$.  
Then $K$ is a maximal compact subgroup of $G$, and
$K$ is connected.
%as a consequence of the Cartan decomposition.

The Cartan involution $\theta$ acts naturally as a Lie algebra 
automorphism of $\kg$.  %we also denote it by $\theta$. 
Then the Cartan decomposition of $\kg$ is given by
\begin{equation}
\label{CartanDecom}
\kg=\kp\oplus \kk, \;\mathrm{with}\;\kp=\{a\in\kg \;:\;
\theta a=-a\},\;\kk=\{a\in\kg \;:\; \theta a=a\}. 
\end{equation} 
From \eqref{CartanDecom}, we get
\begin{equation}
\label{eq2.2}
[\kp,\kp]\subset\kk,\qquad[\kk,\kk]
\subset \kk,\qquad [\kp,\kk]
\subset \kp.
\end{equation}
Put $n=\dim \kk, m=\dim \kp$. 
Then $\dim \kg=m+n$.

If $ g,h \in G, \;  u \in \kg$, let
$\Ad(g)h=ghg^{-1}$ be the adjoint action of $g$ on $h$,
and let $\mathrm{Ad}(g)u\in\kg$
denote the action of $g$ on $u$ via the adjoint representation.
If $u,v\in \kg$, set 
\begin{equation}\label{eq2.1}
\mathrm{ad}(u)v=[u,v],
\end{equation}
then $\mathrm{ad}$ is the derivative of the map 
$g\in G\rightarrow \mathrm{Ad}(g)\in \mathrm{Aut}(\kg)$.

Let $B$ be a real-valued nondegenerate symmetric bilinear 
form on $\kg$ which is invariant under the adjoint action 
of $G$ on $\kg$, and also under the action of $\theta$. 
Then \eqref{CartanDecom} is an orthogonal splitting of 
$\kg$ with respect to $B$. 
We assume that $B$ is positive on $\kp$ and negative
on $\kk$. Put 
$\langle\cdot,\cdot\rangle=-B(\cdot,\theta \cdot)$ 
the associated scalar product on $\kg$, 
which is invariant under the adjoint action of $K$.
Let $|\cdot|$ be the corresponding norm on $\kg$. 
The splitting \eqref{CartanDecom} is also orthogonal 
with respect to $\langle\cdot,\cdot\rangle$.  

\begin{rema}\label{r2.1}
 For $G=\mathrm{GL}^+({\bf q},\R)
 =\{A\in \mathrm{GL}({\bf q},\R), \det A>0\}$,
 the Cartan involution is given by $\theta(g)={}^tg^{-1}$, 
 where ${}^t\cdot$ denotes the transpose %operation 
 of a matrix. Then $K=\mathrm{SO}({\bf q})$, 
 the special orthogonal group, and $\kk$ is the vector space of 
 anti-symmetric matrices and 
 $\kp$ is the vector space of  symmetric matrices. 
 We can take $B(u,v)=2\tr^{\R^{\bf q}}[uv]$ for
 $u,v\in \kg=\mathfrak{gl}({\bf q},\R)=\mathrm{End}(\R^{\bf q})$.
\end{rema}

Let $U(\kg)$ be the enveloping algebra of $\kg$ which will 
be identified with the algebra of left-invariant differential operators
on $G$. Let $C^\kg\in U(\kg)$ be the Casimir element. 
 If $\{e_i\}^m_{i=1}$ is an orthonormal basis of 
$(\kp,\langle\cdot,\cdot\rangle)$ and 
if $\{e_i\}^{m+n}_{i=m+1}$ is an orthonormal basis of 
$(\kk,\langle\cdot,\cdot\rangle)$, then 
\begin{equation}
\label{Casimir}
C^\kg= C^\kp+C^\kk, \;\mathrm{with}\; 
C^\kp=-\sum^m_{i=1}e_i^2,\; C^\kk=\sum^{m+n}_{i=m+1} e_i ^2.
\end{equation}
Then $C^\kk$ is the Casimir element of $\kk$ with respect to
the bilinear form induced by $B$ on~$\kk$. 
Note that $C^\kg$ lies in the center of $U(\kg)$.

Let $\rho^V: K\rightarrow \mathrm{Aut}(V)$ be an orthogonal 
or unitary representation of $K$ on a finite dimensional Euclidean 
or Hermitian vector space $V$.
%which also induces a representation of Lie algebra $\kk$. 
We denote by $C^{\kk,V}\in\text{End}(V)$ 
the corresponding Casimir operator acting on $V$, given by
\begin{equation}
\label{eq2.4}
C^{\kk,V}=\sum^{m+n}_{i=m+1} \rho^{V,2}(e_i).
\end{equation}

Let 
\begin{align}\label{eq:2.3a}
p: G\to X=G/K
\end{align}
be the quotient space.
%If $X=G/K$ is the quotient space, 
Then $X$ is contractible. 
More precisely, $X$ is a symmetric space and the exponential map 
$\exp: \kp\rightarrow G/K,\; a\mapsto  p e^{a}$ is 
a diffeomorphism. We have a natural identification of
vector bundles on $X$:
\begin{equation}
\label{eq2.3}
TX=G\times_K \kp,
\end{equation}
where $K$ acts on  $\kp$ via the adjoint representation.
%where the vector bundle $G\times_K \kp$ is induced by
%the adjoint representation of $K$ on $\kp$.
The scalar product of $\kp$ descends to a Riemannian metric 
$g^{TX}$ on $TX$. 
Let $\omega^\kg$ be the canonical left-invariant $1$-form on $G$ 
with values in $\kg$, and let $\omega^\kk$ be the 
$\kk$-component of $\omega^\kg$. Then $\omega^\kk$ 
defines a connection on the $K$-principal bundle $G\rightarrow G/K$.
The connection $\nabla^{TX}$ on $TX$
induced by $\omega^\kk$ and by \eqref{eq2.3} is precisely 
the Levi-Civita connection on $(TX, g^{TX})$.

Note since the adjoint representation of $K$ preserves $\kp$ 
and $\kk$, %in \eqref{CartanDecom}, and this defines 
we obtain $C^{\kk,\kp}\in\mathrm{End}(\kp),\; C^{\kk,\kk}\in
\mathrm{End}(\kk)$.
In fact, $\tr^\kp[C^{\kk,\kp}]$ 
is the scalar curvature of $X$, and 
$-\dfrac{1}{4}\tr^\kk[C^{\kk,\kk}]$ is the scalar curvature of $K$ 
for the Riemannian structure induced by $B$ (cf. \cite[(2.6.8)
and (2.6.9)]{B11b}).

Let $\rho^E: K\rightarrow \mathrm{Aut}(E)$ be
a unitary representation of $K$. Then the vector space $E$ 
descends to 
a Hermitian vector bundle $F=G\times_K E$ on $X$, and 
$\omega^\kk$ induces a Hermitian connection $\nabla^F$ on $F$. 
Then $C^\infty(X,F)$ can be identified to $C^\infty(G,E)^K$,
the $K$-invariant part of $C^\infty(G,E)$. 
The Casimir operator $C^\kg$, %a priori 
acting on $C^\infty(G,E)$, descends to an operator 
acting on $C^\infty(X,F)$, which will still be denoted by 
$C^\kg$.

Let $A$ be a self-adjoint endomorphism of $E$ which is 
$K$-invariant. Then $A$ descends to a parallel self-adjoint section 
of $\text{End}(F)$ over $X$.

\begin{defi}\label{d2.1}
	Let $\mL^X, \mL^X_A$ act on 
	$C^\infty(X,F)$ by the formulas,
\begin{align}\label{eq2.5}
    \begin{split}
&\mL^X=\frac{1}{2}C^\kg
+\frac{1}{16}\tr^\kp
[C^{\kk,\kp}]
+\frac{1}{48}\tr^{\kk}[C^{\kk,
\kk}];\\
&\mL^X_A=\mL^X+A.
\end{split}\end{align}
\end{defi}

From \eqref{Casimir}, $-C^\kp$ descends to 
the Bochner Laplacian $\Delta^F$ on $C^\infty(X,F)$, the operator
$C^\kk$ %becomes $C^{\kk,E}$ and then 
descends to a parallel section $C^{\kk,F}$ of $\End(F)$ on $X$. 
If the representation 
$\rho^E$ above is irreducible, then $C^{\kk,F}$ acts as 
$c\,\mathrm{Id}_F$, where $c$ is a constant function on $X$. 
%In particular, if $E=\C$ is a trivial representation of $K$, then 
%$C^{\kk,F}=0$.
Thus from \eqref{eq1.3} and \eqref{eq2.5}, we have
\begin{equation}
\label{eq2.6}
\mL^X=-\frac{1}{2}\Delta^F_{\phi} \;\;\mathrm{with}\;\; 
\phi=-C^{\kk,F}
- \frac{1}{8} \tr^\kp
[C^{\kk,\kp}]
-\frac{1}{24}\tr^\kk
[C^{\kk,\kk}].
\end{equation}

The group $G$ acts on $X$ on the left. This action lifts to 
$F$. More precisely, for any $h\in G$ and $[g,v]\in F$, 
the left action of $h$ is given by
\begin{equation}
\label{eq2.8}
h.[g,v]=[hg,v]\in G\times_K E=F.
\end{equation}
Then the operators $\mL^X,\;\mL^X_A$ 
commute with $G$. %this $G-$action.

Let $\Gamma\subset G$ be a discrete subgroup of $G$ 
such that the quotient space $\Gamma\backslash G$ is compact. 
Set
\begin{equation}
\label{eq2.9}
	Z=\Gamma\backslash X=\Gamma\backslash G/ K. 
\end{equation}
Then $Z$ is a compact locally symmetric space. %Usually, 
In general $Z$ is an orbifold. If $\Gamma$ is torsion-free  (i.e., 
if $\gamma\in\Gamma$, $k\in\N^*$,  
then $\gamma^k=1$ implies $\gamma=1$), 
then $Z$ is a smooth manifold.

From now on, we assume that $\Gamma$ is torsion free, so that 
$\Gamma=\pi_1(Z)$ and $X$ is just the universal cover of $Z$.

A vector bundle like $F$ on $X$ descends to a vector bundle 
on $Z$, which we still denote by $F$. 
Then the operators $\mL^X,\;\mL^X_A$ 
descend to operators $\mL^Z,\;\mL^Z_A$ 
acting on $C^\infty(Z,F)$.

For $t>0$, 
let $e^{-t \mL^X_A}(x,x')\;(x,x'\in X), 
e^{-t \mL^Z_A}(z,z')\;(z,z'\in Z)$ 
be the smooth kernels of the heat operators $e^{-t \mL^X_A},\; 
e^{-t \mL^Z_A}$ with respect to the Riemannian volume 
forms $dx', dz'$ respectively. By \eqref{eq1.6}, we get
\begin{eqnarray}
	\tr[e^{-t \mL^Z_A}] &
	=& \int_Z \tr[e^{-t \mL^Z_A}(z,z)]dz\label{eq2.10}\\
&=& \int_{\Gamma\backslash X} \sum_{\gamma\in\Gamma} 
\tr[\gamma e^{-t \mL^X_A}(\gamma^{-1}\widetilde{z},
\widetilde{z})]dz.
\nonumber
\end{eqnarray}
% We will now  reorganize the integral formula for 
% $\tr[e^{-t \mL^Z_A}]$ from \eqref{eq2.10} as integrals on 
% $\Gamma\backslash G$.

%%%%%%%%%%%%%%%%%%%%%%%%%%%%%
\subsection{Orbital integrals and Selberg trace formula}\label{s2.2}

Let $C^b(X,F)$ be the vector space of continuous bounded sections 
of $F$ over $X$. Let $Q$ be an operator acting on $C^b(X,F)$ 
with a continuous kernel $q(x,x')$ with respect to
the volume form $dx'$. It is convenient to view $q$ as
a continuous function $q(g,g')$ defined on $G\times G$ 
with values in $\End(E)$ which satisfies  for 
any $k,k'\in K$,
\begin{equation}
\label{eq2.11}
q(gk,g'k')=\rho^E(k^{-1})q(g,g')\rho^E(k').
\end{equation}

Now we assume that the operator $Q$ commutes with 
the left action of $G$ on $C^b(X,F)$ defined in \eqref{eq2.8}. 
This is equivalent to 
\begin{equation}
\label{eq2.12}
q(gx,gx')=gq(x,x')g^{-1} \quad \text{ for any }\, x,x'\in X,\; g\in G,
\end{equation} 
where %$x,x'\in X,\; g\in G$ and 
the action of $g^{-1}$ maps 
$F_{gx'}$ to $F_{x'}$, the action of $g$ maps $F_{x}$ to $F_{gx}$.

If we consider instead the kernel $q(g,g')$, then this implies 
that for all $g''\in G$,
\begin{equation}
\label{eq2.13}
q(g''g,g''g')=q(g,g')\in \End(E).
\end{equation}
Thus the kernel $q$ is determined %characterized 
by $q(1,g)$. Set
\begin{equation}
\label{eq2.14}
q(g)=q(1,g).
\end{equation}
Then we obtain from \eqref{eq2.11} and \eqref{eq2.13} 
that for $g\in G, k\in K$,
\begin{equation}
\label{eq2.15}
q(k^{-1}gk)=\rho^E(k^{-1})q(g)\rho^E(k).
\end{equation}
This implies that $\tr^E[q(g)]$ is invariant when replacing
$g$ by $k^{-1}gk$. 

In the sequel, we will use the same notation $q$ for the various 
versions of the corresponding kernel $Q$.

\begin{defi}\label{d2.4}
The element $\gamma\in G$ is said to be elliptic if it
is conjugate in 
$G$ to an element of $K$. We say that $\gamma$ is hyperbolic 
if it is conjugate in $G$ to $e^a,\; a\in \kp$. 

For $\gamma\in G$, $\gamma$ is semisimple if there exist 
$g\in G$, $a\in \kp$, $k\in K$ such that 
\begin{align}
\label{eq2.2.2a}
\Ad(k)a=a, \quad \gamma=\Ad(g)\left(e^ak^{-1}\right).
\end{align}
\end{defi}

By \cite[Theorem 2.19.23]{Eberlein96}, if $\gamma\in G$ is
a semisimple element, $\Ad(g)e^a$ and $\Ad(g)k^{-1}$ are 
uniquely determined by $\gamma$ (i.e., they do not depend on 
$g\in G$ such that \eqref{eq2.2.2a} holds), and 
\begin{align} \label{eq2.2.3a}
Z(\gamma)=Z\left(\Ad(g)e^a\right)\cap Z\left(\Ad(g)k^{-1}\right),
\end{align}
where $Z(\gamma)\subset G$ is the centralizer of $\gamma$ in $G$.

Let $dk$ be the Haar measure on $K$ that gives 
volume $1$ to $K$. Let $dg$ be measure on $G$ 
(as a $K$-principal bundle on $X=G/K$) given by
\begin{equation}
\label{volumeform}
dg=dx\;dk.
\end{equation}
Then $dg$ is a left-invariant Haar measure on $G$. 
Since $G$ is unimodular, it is also a right-invariant Haar measure.

For $\gamma\in G$ semisimple, $Z(\gamma)$ is reductive and 
$K(\gamma)$, the fixed-points set of $\Ad(g) \theta \Ad(g)^{-1}$ 
 in $Z(\gamma)$ (cf. (\ref{eq2.2.3a})),
is a maximal compact subgroup.
Let $dy$ be the volume element on
the symmetric space $X(\gamma)=Z(\gamma)/K(\gamma)$ 
induced by $B$. Let $dk'$ be the Haar measure on 
$K(\gamma)$    %=K\cap Z(\gamma)$
that gives volume $1$ to 
$K(\gamma)$. Then $dz=dydk'$ is a left and right Haar measure 
on $Z(\gamma)$. Let $dv$ be the canonical measure on 
$Z(\gamma)\backslash G$ that is canonically associated with $dg$ 
and $dz$ so that 
\begin{align}
\label{eq2.2.1}
dg=dzdv. 
\end{align}

\begin{defi}
 [Orbital integral] For $\gamma\in G$ semisimple, 
 we define the orbital integral associated with $Q$ and $\gamma$
 by 
 \begin{align}
 \label{eq2.2.2}
 \tr^{[\gamma]}[Q]=\int_{Z(\gamma)\backslash G} 
 \tr^E[q(v^{-1} \gamma  v)]dv,
 \end{align}
 once the integral converges. 
\end{defi}

Note that the map
\begin{align}\label{eq2.2.3}
Z(\gamma) \backslash G\to \mathcal{O}_{\gamma} 
=\Ad_G\gamma \, \text{ given by }\, v\to v^{-1}\gamma v
\end{align}
identifies $Z(\gamma)\backslash G$ as the orbit 
$\mathcal{O}_{\gamma}$ of $\gamma$ with the adjoint action 
of $G$ on $G$. This justifies the name ``orbital integral" 
for \eqref{eq2.2.2}.

Let $\Gamma\subset G$ be a discrete torsion free 
cocompact subgroup as in Section \ref{s2.1}.
Since the operator $Q$ commutes with the left action of $G$,
$Q$ descends to an operator $Q^Z$ acting on $C^\infty(Z,F)$. 
We assume that the sum 
$\sum_{\gamma\in\Gamma} q(g^{-1}\gamma g') $
is uniformly and absolutely convergent on $G\times G$.
%and that $Q^Z$ is a trace class operator.

Let $[\Gamma]$ be the set of conjugacy classes in $\Gamma$. 
If $[\gamma]\in [\Gamma]$, set 
\begin{equation}
\label{kernelgamma}
	q^{X,[\gamma]}(g,g')=\sum_{\gamma'\in [\gamma]} 
	q(g^{-1}\gamma' g').
\end{equation}
Then from (\ref{eq2.13})--(\ref{kernelgamma}), we get 
\begin{equation}
\label{eq2.17}
q^Z(z,z')=\sum_{[\gamma]\in [\Gamma]} q^{X,[\gamma]}(g,g'),
\end{equation}
with $g,g'\in G$ fixed lift %lifting 
of $z,z'\in Z$. Thus as in \eqref{eq2.10},
\begin{equation}
\label{eq2.19}
\tr[Q^Z]=\sum_{[\gamma]\in [\Gamma]}
\tr[Q^{Z,[\gamma]}] \;\;\mathrm{with}\;\;
\tr[Q^{Z,[\gamma]}]
=\int_Z \tr[q^{X,[\gamma]}(z,z)]dz.
\end{equation}

From \eqref{volumeform}, \eqref{kernelgamma}, \eqref{eq2.19},
and the fact that 
$[\gamma]\simeq \Gamma\cap Z(\gamma)\backslash \Gamma$, 
we have
\begin{align}
\label{eq2.21}\begin{split}
\tr[Q^{Z,[\gamma]}]&=\int_{\Gamma\cap Z(\gamma)
\backslash G} \tr^E[q(g^{-1}\gamma g)]dg\\
&=\vol \Big(\Gamma\cap Z(\gamma)\backslash Z(\gamma)\Big)
\tr^{[\gamma]}[Q]\\
&=\vol\Big(\Gamma\cap Z(\gamma)\backslash X(\gamma)\Big)
\tr^{[\gamma]}[Q].
\end{split} \end{align}
 
From \eqref{eq2.19} and  \eqref{eq2.21},  we get 
\begin{theo}[Selberg trace formula]
 \label{t2.5a} \begin{align}
 \label{eq2.22}
 \tr\left[Q^Z\right]
 =\sum_{[\gamma]\in [\Gamma]} 
 \vol \Big(\Gamma\cap Z(\gamma)\backslash X(\gamma)\Big) 
 \tr^{[\gamma]}[Q]. 
\end{align}
\end{theo}
 Selberg \cite[(3.2)]{Selberg56} 
was the first to give a closed formula for %evaluated first precisely 
the trace of the heat operator on
a compact hyperbolic Riemann surface via \eqref{eq2.22}, 
which is the original Selberg trace formula. 
Harish-Chandra's Plancherel theory,
developed from the 1950s until the 1970s, is 
an algorithm to reduce the computation of 
an orbital integral to a lower dimensional group by
the discrete series method, cf. Section \ref{s2.6}.

To understand
better the structure of each integral in \eqref{eq2.2.2},
we first reformulate it in more geometric terms.  %way.

%%%%%%%%%%%%%%%%%%%%%%%%%%%%%%%%%%%
\subsection{Geometric orbital integrals}\label{s2.3}

Let $d(\cdot,\cdot)$ be the Riemannian distance on $X$.
If $\gamma\in G$, the displacement function $d_\gamma$ 
is given by for $x\in X$,
\begin{equation}
\label{eq2.3.1}
d_\gamma(x)=d(x,\gamma x).
\end{equation}
By \cite[\S 6.1]{BaGSc85}, the function $d_\gamma$ is 
convex on $X$, i.e., for any geodesic $t\in \R\to x_t\in X$ with 
constant speed, the function $d_\gamma(x_{t})$ is convex
on $t\in \R$. 

Recall that $p: G\rightarrow X=G/K$
is the natural projection in (\ref{eq:2.3a}). 
%Let $Z(g)\subset G$ denote the centralizer of $g$ in $G$.
We have the following geometric %algebraic 
description on the semisimple elements in $G$.
\begin{theo}\label{t2.5}
\cite[Theorem 3.1.2]{B11b}.
The element  $\gamma\in G$ is semisimple if and only if 
the function $d_\gamma$ %has a minimum value in $X$. 
attains its minimum in $X$. If $\gamma\in G$ is semisimple, and 
\begin{align}
\label{eq2.3.2a}
X(\gamma)=\{x\in X: d_\gamma(x)=m_\gamma:
=\inf_{y\in X}d_\gamma(y)\},
\end{align}
for $g\in G$, $x=pg\in X$, then $x\in X(\gamma)$ if and only if
there exist $a\in\kp,\; k\in K$ 
such that 
\begin{equation}\label{eq2.3.2}
\gamma=\Ad(g)(e^ak^{-1})\quad  \text{ and }\, \Ad(k)a=a.
\end{equation}
If $g_t=ge^{ta}$, 
then $t\in [0,1]\rightarrow x_t=pg_t$ is the unique geodesic 
connecting $x\in X(\gamma)$ and $\gamma x$ in $X$. 
Moreover, we have
	\begin{equation}
	\label{eq2.3.3}
	m_\gamma=|a|.
	\end{equation}
\end{theo}

Since the integral \eqref{eq2.21} depends only on the conjugacy
class of $\gamma$, from Theorem \ref{t2.5} or (\ref{eq2.2.2a}),
we may and we will assume that 
\begin{equation}
\label{eq2.3.4}
\gamma=e^ak^{-1},\; \quad \Ad(k)a=a,\; \quad
	a\in \kp,\;k\in K.
\end{equation} 
Furthermore,  by (\ref{eq2.2.3a}), we have
\begin{equation}
\label{eq2.3.5}
	Z(\gamma)=Z(e^a)\cap Z(k),\quad \kz(\gamma)
	=\kz(e^a)\cap \kz(k),
\end{equation}
where we use the symbol $\mathfrak{z}$ to denote the 
corresponding Lie algebras of the centralizers.

Put 
\begin{equation}
\label{eq2.3.6}
	\kp(\gamma)=\kz(\gamma)\cap \kp,\quad
	\kk(\gamma)=\kz(\gamma)\cap \kk.
\end{equation}

From \eqref{eq2.2} and \eqref{eq2.3.5}, we get
\begin{equation}
\label{eq2.3.7}
	\kz(\gamma)=\kp(\gamma)\oplus \kk(\gamma).
\end{equation}
Thus the restriction of $B$ to $\kz(\gamma)$ is non-degenerate.
Let $\kz^\perp(\gamma)$  be the orthogonal space to 
$\kz(\gamma)$ in $\kg$ with respect to $B$. 
Then $\kz^\perp(\gamma)$ splits as %follow,
\begin{equation}
\label{eq2.3.8}
\kz^\perp(\gamma)=\kp^\perp(\gamma)\oplus 
\kk^\perp(\gamma),
\end{equation} 
where $\kp^\perp(\gamma)\subset \kp,\; 
\kk^\perp(\gamma)\subset \kk$ are the orthogonal spaces to 
$\kp(\gamma),\;\kk(\gamma)$ in $\kp,\;\kk$ with respect to 
the scalar product induced by $B$.

Set 
\begin{equation}
\label{eq2.3.9}
K(\gamma)=K\cap Z(\gamma),
\end{equation}
then from \eqref{eq2.3.5} and \eqref{eq2.3.6}, $\kk(\gamma)$ is 
just the Lie algebra of $K(\gamma)$.

\begin{theo}
    \cite[Theorems 3.3.1, 3.4.1, 3.4.3]{B11b}
\label{t2.6}
The set $X(\gamma)$ is a submanifold of~$X$.
In the geodesic coordinate system centered at $p1$, 
we have the identification
\begin{equation}	\label{eq2.3.10}
	X(\gamma)=\kp(\gamma).
\end{equation}
The action of $Z(\gamma)$ on $X(\gamma)$ is transitive and 
we have the identification of $Z(\gamma)$-manifolds,
\begin{equation}\label{eq2.3.11}
	X(\gamma)\simeq Z(\gamma)/K(\gamma).
\end{equation}
	The map
	\begin{equation}
	\label{eq2.3.12}
\rho_\gamma: (g,f,k')\in Z(\gamma)\times_{K(\gamma)}
(\kp^\perp(\gamma)\times K)\rightarrow ge^fk'\in G
	\end{equation}	
is a diffeomorphism of left $Z(\gamma)$-spaces, 
	and of right $K$-spaces.	
The map $(g,f,k')\mapsto (g,f)$ corresponds to the projection 
$p: G\rightarrow X=G/K$. In particular,  the map
\begin{equation}
\label{2.3.14}
 \rho_{\gamma}: (g,f)\in Z(\gamma)\times_{K(\gamma)} 
 \kp^\perp(\gamma)\rightarrow p(ge^f)\in X
\end{equation}
is a diffeomorphism.

Moreover, under the diffeomorphism (\ref{eq2.3.12}), 
we have the identity of right $K$-spaces,
\begin{equation}
	\label{eq2.3.13}
\kp^\perp(\gamma)_{K(\gamma)}\times K
= Z(\gamma)\backslash G.
\end{equation}
Finally, there exists $C_\gamma>0$ such that if 
$f\in \kp^\bot(\gamma)$, $|f|>1$,
\begin{align}
\label{eq2.3.14a}
d_{\gamma}(\rho_\gamma(1,f))\ge |a|+C_\gamma |f|. 
\end{align}
\end{theo}

% From Theorem \ref{t2.6}, the map
% \begin{equation}
% \label{2.3.14}
%  \rho_{\gamma}: (g,f)\in Z(\gamma)\times_{K(\gamma)} 
%  \kp^\perp(\gamma)\rightarrow p(ge^f)\in X
% \end{equation}
% is a diffeomorphism. This
The map $\rho_\gamma$ in (\ref{2.3.14}) is the normal coordinate
system on $X$ based at $X(\gamma)$.
\begin{figure}[htbp]
\begin{center}
    \includegraphics[scale=1]{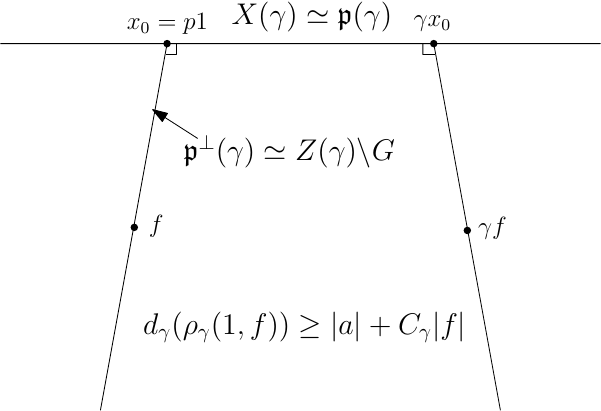}
\caption{Normal coordinate}
\label{default}
\end{center}
\end{figure}

Recall that $dy$ is the volume element on $X(\gamma)$ 
(cf. Section \ref{s2.2}).
%induced by the bilinear form $B$ via \eqref{eq2.3.11}, 
Let $df$ be the volume element on $\kp^\perp(\gamma)$. 
Then $dydf$ is a volume form on 
$Z(\gamma)\times_{K(\gamma)}\kp^\perp(\gamma)$ 
that is $Z(\gamma)$-invariant.  Let $r(f)$ be the smooth function 
on $\kp^\perp(\gamma)$ that is $K(\gamma)$-invariant 
such that we have the identity of volume element on $X$
via \eqref{2.3.14},
\begin{equation}
\label{eq2.3.14}
	dx= r(f)dydf, \;\mathrm{with}\; r(0)=1.
\end{equation}	

In view of \eqref{eq2.3.13}, \eqref{eq2.3.14}, Bismut could
reformulate geometrically the orbital integral \eqref{eq2.2.2} 
as an integral along the normal direction of $X(\gamma)$ in $X$.
\begin{prop}[Geometric orbital integral]
\label{d2.7}
The orbital integral for the operator $Q$ in 
Section \ref{s2.2} and a semisimple element $\gamma\in G$ 
is given by 
\begin{equation}
	\label{eq2.3.15}
\tr^{[\gamma]}[Q]=\int_{\kp^\perp(\gamma)} 
\tr^E[q(e^{-f}\gamma e^f)]r(f)df.
\end{equation}
\end{prop}

Equation  (\ref{eq2.3.15})  gives a geometric interpretation 
for orbital integrals.
It is remarkable that even before its explicit computation,
the variational problem connected with the minimization of the 
displacement function $d_{\gamma}$ is used in %formula 
(\ref{eq2.3.15}).
% to observe from Bismut's geometric formulation 
% of orbital integrals (\ref{eq2.3.15}) that the orbital integral
% associated with $\gamma\in G$ is an integral along 
% the normal direction of the minimum set $X(\gamma)$
% of the displacement function $d_{\gamma}$ in (\ref{eq2.3.1})
% and $d_{\gamma}$ grows at least linearly along the normal direction.
% This explains in a geometric way %why 
% that the orbital integral 
% is well-defined for the heat kernel.
% The set that our final formula should be localized is given 
% before we really start to work, it is $X(\gamma)$ 
% as the minimum set of $d_{\gamma}$.

We need the following criterion for %property on 
the semisimplicity of an element.
\begin{prop} {\rm (Selberg \cite[Lemmas 1, 2]{Selberg60})}
	\label{p2.9}
If $\Gamma\subset G$ is a discrete cocompact subgroup, 
then for any $\gamma\in \Gamma$, $\gamma$ is semisimple,
and $\Gamma \cap Z(\gamma)$ is cocompact in $Z(\gamma)$.
\end{prop}
\begin{proof}
Let $U$ be a compact subset of $G$ such that $G=\Gamma\cdot U$. 
Let $\gamma\in \Gamma$. Let $\{x_k\}_{k\in \N}$ be a %series 
family of points in $X$ such that 
$d(x_k,\gamma x_k)\rightarrow m_\gamma
=\inf_{x\in X}d(x,\gamma x)$ as $k\rightarrow +\infty$.
	
Then there exists $\gamma_k\in \Gamma,\; x'_k\in U$ such that 
$\gamma_k x'_k=x_k$. Since $U$ is compact, 
there is a subsequence 
$\{x'_{k_j}\}_{j\in\N}$ of $\{x'_k\}_{k\in\N}$ such that 
as $j\rightarrow +\infty$, $x'_{k_j}\rightarrow y\in U$. Then
\begin{equation}\label{eq2.3.18}
	\begin{split}
d(y,\gamma^{-1}_{k_j}\gamma\gamma_{k_j}y)&\leq d(x'_{k_j},y)
+d(x'_{k_j}, \gamma^{-1}_{k_j}\gamma\gamma_{k_j}x'_{k_j})
+d(\gamma^{-1}_{k_j}\gamma
\gamma_{k_j}x'_{k_j},\gamma^{-1}_{k_j}\gamma\gamma_{k_j}y)\\
	&=2d(x'_{k_j},y)+d(x_{k_j}, \gamma x_{k_j}),
	\end{split}
\end{equation} 
where the right side tends to $m_\gamma$ 
as $j\rightarrow +\infty$.
	
Since $\Gamma$ is discrete and each 
$\gamma^{-1}_{k_j}\gamma\gamma_{k_j}\in \Gamma$, 
 the set of such $\gamma^{-1}_{k_j}\gamma\gamma_{k_j}$ 
is bounded, so that there exist infinitely many $j$ such that 
$\gamma^{-1}_{k_j}\gamma\gamma_{k_j}=\gamma'\in \Gamma$. 
Then 
\begin{equation}\label{eq2.3.19}
m_\gamma=d(y,\gamma'y)
=d(\gamma_{k_j}y,\gamma \gamma_{k_j}y). 
\end{equation}
This means that $d_\gamma$ reaches its minimum in $X$. 
Therefore $\gamma$ is semisimple.

Since $\Gamma$ is discrete, $[\gamma]$ is closed in $G$,
thus $\Gamma\cdot Z(\gamma)$ as the inverse image
of $[\gamma]$ of the continuous map
$g\in G\to g \gamma g^{-1}\in G$, is closed in $G$.
This implies $\Gamma \cap Z(\gamma)\backslash Z(\gamma)
= \Gamma \backslash \Gamma\cdot Z(\gamma)$ 
is a closed subset of the compact quotient $\Gamma \backslash G$.
Thus  $\Gamma \cap Z(\gamma)$ is cocompact in $Z(\gamma)$.
\end{proof}

Let $\Gamma\subset G$ be a discrete torsion
free cocompact subgroup as in Section \ref{s2.2}. 
\comment{
From \eqref{eq2.21}, \eqref{eq2.3.13} and \eqref{eq2.3.15}, 
we get for $\gamma\in \Gamma$,
\begin{equation}
\label{eq2.3.16}
\tr[Q^{Z,[\gamma]}]=\vol \Big(\Gamma\cap Z(\gamma)
\backslash X(\gamma)\Big)     \tr^{[\gamma]}[Q].
\end{equation}
Combing \eqref{eq2.19} and \eqref{eq2.3.16}, we can reformulate 
Theorem \ref{t2.5a} as  %the Selberg trace formula,
\begin{theo}[Selberg trace formula]\label{t2.8}
	\begin{equation}
	\label{eq2.3.17}
\tr[Q^Z]=\sum_{[\gamma]\in [\Gamma]} 
\vol \Big(\Gamma\cap Z(\gamma)\backslash
X(\gamma)\Big)  \tr^{[\gamma]}[Q].
	\end{equation}
\end{theo}
}
%Note that 
Set $Z= \Gamma\backslash X$, then $\Gamma=\pi_1(Z)$.
For $x\in X(\gamma)$, 
the unique geodesic from $x$ to $\gamma x$ descends to 
the closed geodesic in $Z$ 
in the homotopy class $\gamma\in \Gamma$ 
which has the shortest length %distance 
$m_\gamma$. 
Thus the Selberg trace formula \eqref{eq2.22} %\eqref{eq2.3.17}
relates the trace of an operator 
$Q$ to the dynamical  properties of the geodesic flow on $Z$ via 
orbital integrals.

%%%%%%%%%%%%%%%%%%%%%%%%%%%%%%%%%%
\subsection{Bismut's explicit formula for orbital integrals}\label{s2.4}
By the standard heat kernel estimate, 
for the heat operator $e^{-t\mL_A^X}$ on $X$, 
%there exists $c>0$ such that for every $t_{0}>0$,
there exist $c>0$,  $\lambda, C>0$, $M>0$ such that for any %$t>t_{0}$,
$t>0$, $x,x'\in X$, we have 
(cf. for instance \cite[(3.1)]{MM15})
\begin{align}
\label{eq2.4.1a}
\left|e^{-t\mL_A^X}(x,x')\right|\le C t^{-M}
e^{\lambda  t 
-c \, d^2(x,x') /t}. 
\end{align}
Note also that by Rauch's comparison theorem, there exist
$C_0, C_1>0$ such that for all $f\in \kp^\bot(\gamma)$,
\begin{align}
\label{eq3.6.9}
|r(f)|\le C_0 e^{C_1 |f|}. 
\end{align}
From \eqref{eq2.3.14a}, \eqref{eq2.4.1a} and \eqref{eq3.6.9}, 
the orbital integral $\tr^{[\gamma]}[e^{-t\mL_A^X}]$
is well-defined for any semisimple element $\gamma\in G$.

Let $\gamma\in G$ be the semisimple element as 
in \eqref{eq2.3.4}. Set
\begin{equation}
\label{eq2.4.1}
\kp_0=\kz(a)\cap\kp,\;\quad  \kk_0=\kz(a)\cap \kk,\;\quad
	\kz_0=\kz(a)=\kp_0\oplus\kk_0.
\end{equation}
Let $\kz^\perp_0$ be the orthogonal space to $\kz_0$ in $\kg$
with respect to $B$. 

Let $\kp^\perp_0(\gamma)$ be the orthogonal to $\kp(\gamma)$
in $\kp_0$, and let $\kk^\perp_0(\gamma)$ be 
the orthogonal space to $\kk(\gamma)$ in $\kk_0$. 
Then the orthogonal space to $\kz(\gamma)$ in $\kz_0$ is
\begin{equation}
\label{eq2.4.2}
\kz^\perp_0(\gamma)=\kp^\perp_0(\gamma)\oplus
\kk^\perp_0(\gamma).
\end{equation}

For $Y_0^\mathfrak{k}\in \kk(\gamma)$, we claim that 
\begin{equation}
\label{eq2.4.3}
\det\Big(1-\exp(-i\theta \ad(Y_0^\kk))
\Ad(k^{-1})\Big)|_{\kz^\perp_0(\gamma)} 
\det\Big(1-\Ad(k^{-1})\Big)|_{\kz^\perp_0(\gamma)}
\end{equation}
has a natural square root, which depends analytically on $Y_0^\kk$. 
Indeed, $\ad(Y_0^\kk)$ commutes with 
$\Ad(k^{-1})$, and no eigenvalue of $\Ad(k)$ 
acting on $\kz^\perp_0(\gamma)$ is equal to $1$. 
If $\kz^\perp_0(\gamma)$ is $1$-dimensional, then 
$\Ad(k)|_{\kz_0^\perp(\gamma)}=-1$ and  
$\ad(Y_0^\kk)|_{\kz_0^\perp(\gamma)}=0$, 
the square root is just $2$. If $\kz^\perp_0(\gamma)$ is 
$2$-dimensional, if $\Ad(k)|_{\kz_0^\perp(\gamma)}$ 
is a rotation of angle $\phi$ and  
$\theta\ad(Y_0^\kk)|_{\kz_0^\perp(\gamma)}$ 
acts by an infinitesimal rotation of angle $\phi'$, such a square root
is given by (cf. \cite[(5.4.10)]{B11b})
\begin{equation}
\label{eq2.4.4}
4 \sin\Big(\frac{\phi}{2}\Big)\sin\Big(\frac{\phi+i\phi'}{2}\Big).
\end{equation}

If $V$ is a finite dimensional Hermitian vector space and if 
$\Theta\in \mathrm{End}(V)$ is self-adjoint, then
$\dfrac{\Theta/2}{\sinh(\Theta/2)}$ is a self-adjoint positive 
endomorphism.  Set
\begin{equation}
\label{eq2.4.5}
\widehat{A}(\Theta)=\mathrm{det}^{1/2} 
\Big [\frac{\Theta/2}{\sinh(\Theta/2)}\Big ].
\end{equation}
In \eqref{eq2.4.5}, the square root is taken to be the positive
square root.

For $Y^\kk_0\in \kk(\gamma)$, set
\begin{align}
\label{eq2.4.6}
\begin{split}
J_\gamma(Y^\kk_0) = & \frac{1}
{\big|\det (1-\Ad(\gamma))|_{\kz^\perp_0}\big|^{1/2}}
\cdot \frac{\widehat{A}(i\ad(Y^\kk_0)|_{\kp(\gamma)})}
{\widehat{A}\big(i\ad(Y^\kk_0)|_{\kk(\gamma)}\big)} \cdot\\
& \left[\frac{1}{\det (1-\Ad(k^{-1}))
|_{\kz_0^\perp(\gamma)}} 
\frac{\det\Big(1-\exp(-i \ad(Y_0^\kk))
\Ad(k^{-1})\Big)|_{\kk^\perp_0(\gamma)}}
{\det\Big(1-\exp(-i\ad(Y_0^\kk))
\Ad(k^{-1})\Big)|_{\kp^\perp_0(\gamma)}} \right]^{1/2}.	
\end{split}
\end{align}

From \eqref{eq2.4.3}, we know that \eqref{eq2.4.6} is well-defined. 
Moreover, there exist $c_\gamma,\; C_\gamma>0$ such that
for any $Y^\kk_0\in \kk(\gamma)$
\begin{equation}
\label{eq2.4.7}
|J_\gamma(Y^\kk_0)|\leq c_\gamma \, e^{C_\gamma|Y^\kk_0|}.
\end{equation}

We note that $p=\dim \kp(\gamma),\; q=\dim \kk(\gamma)$ 
and $r=\dim \kz(\gamma)=p+q$.
Now we can restate Theorem \ref{t0.1} as follows.

\begin{theo}\cite[Theorem 6.1.1]{B11b}\label{t2.7}
	For any $t>0$, we have
\begin{align}\label{eq2.4.8}
    \begin{split}
&\tr^{[\gamma]}[e^{-t \mL^X_A}] 
=\frac{e^{-|a|^2/2t}}{(2\pi t)^{p/2}}
%&\hspace{1cm}\cdot 
\int_{\kk(\gamma)} J_\gamma(Y^\kk_0) 
\tr^E\Big[\rho^E(k^{-1})e^{-i \rho^E(Y^\kk_0)-tA}\Big]
%\\&\hspace{40mm}
e^{-|Y^\kk_0|^2/2t}
\frac{dY^\kk_0}{(2\pi t)^{q/2}}.
\end{split}\end{align}
\end{theo}
\begin{rema}\label{r2.8}
For $\gamma=1$, we have $\kk(1)=\kk$, $\kp(1)=\kp$, and for 
$Y^\kk_0\in\kk$, by (\ref{eq2.4.6}),
\begin{equation}\label{2.4.8}
J_1(Y^\kk_0)=\frac{\widehat{A}(i\ad(Y^\kk_0)|_\kp)}
{\widehat{A}(i\ad(Y^\kk_0)|_\kk)}.
\end{equation}
\end{rema}

Let $\mathcal{S}(\R)$ be the Schwartz space of $\R$. 
Let $\tr^{[\gamma]}[\cos(s\sqrt{\mL_A^X})]$ be the even 
distribution on $\R$ determined by the condition
that for any even function $\mu\in \mathcal{S}(\R)$ 
with compactly supported Fourier transformation $\widehat{\mu}$, 
we have
\begin{align}
\tr^{[\gamma]}\left[\mu\left(\sqrt{\mL^X_A}\right)\right]
=\int_\R \widehat{\mu}(s)\tr^{[\gamma]}
\left[\cos\left(2\pi s\sqrt{\mL_A^X}\right)\right] ds.
\end{align}
The wave operator $\cos( \sqrt{2}\pi s\sqrt{\mL_A^X})$ defines 
a distribution on $\R\times X\times X$. 

Let $\Delta^{\kz(\gamma)}$ be the standard Laplacian on 
$\kz(\gamma)$ with respect to the scalar product 
$\langle\cdot,\cdot\rangle=-B(\cdot,\theta\cdot)$. 
%From Theorem \ref{t2.7}, we can get 
Now we can state the following 
microlocal version of Theorem \ref{t2.7} for the wave operator.
\begin{theo}
 \label{t2.10}\cite[Theorem 6.3.2]{B11b}. 
 We have the following identity of even distributions
 on $\R$ supported on
 $\{|s|\ge \sqrt{2}|a|\}$ and with singular support %include 
 in $\pm\sqrt{2}|a|$, 
 \begin{align}
\tr^{[\gamma]}\left[\cos\left(s\sqrt{\mL_A^X}\right)\right]
=\int_{H^\gamma} 
\tr^E\left[\cos\left(s\sqrt{-\frac{1}{2}\Delta^{\kz(\gamma)}+A}
\right)
J_\gamma(Y^\kk_0)\rho^E(k^{-1})e^{-i\rho^E(Y_0^\kk)}\right],
 \end{align}
 where $H^\gamma=\{0\}\times (a,\kk(\gamma))\subset 
 \kz(\gamma)\times \kz(\gamma)$. 
\end{theo}

%Let $e^{t\Delta^{\kz(\gamma)}/2}$ be the corresponding 
%heat operator for $t>0$ on $\kz(\gamma)$ with coordinate 
%$(y,Y^\kk_0)\in \kp(\gamma)\oplus\kk(\gamma)=\kz(\gamma)$. 
%Let $\delta_{y=a}$ be the distribution of Dirac measure at point 
%$y=a$ on $\kp(\gamma)$. 
%Then we can reformulate Theorem \ref{t2.7} as
%\begin{equation}
%\label{eq2.4.9}
%\begin{split}
%\tr^{[\gamma]}[e^{-t\mL^X_A}]
%=&\tr^E\bigg[e^{t\Delta^{\kz(\gamma)}/2-t A}
%\Big[J_\gamma(Y^\kk_0)\rho^E(k^{-1}) e^{-i\rho^E(Y^\kk_0)}
%\delta_{y=a}\Big]\bigg](0).
%\end{split}
%\end{equation}

\begin{rema}\label{t2.9}
We assume that the semisimple element $\gamma$ is nonelliptic, 
i.e., $a\neq 0$. We also assume that 
\begin{equation}
\label{2.4.9}
[\kk(\gamma), \kp_0]=0.
\end{equation}
Then for $Y^\kk_0\in\kk(\gamma)$, 
$\mathrm{ad}(Y^\kk_0)|_{\kp(\gamma)}=0,\;
\mathrm{ad}(Y^\kk_0)|_{\kp_0^\perp(\gamma)}=0$.

Now from \eqref{eq2.4.8}, we have 
\cite[Theorem 8.2.1]{B11b}: for $t>0$,
\begin{equation}\label{2.4.10}
\begin{split}
&\tr^{[\gamma]}\Big[e^{-t\mL^X_A}\Big] 
= \frac{e^{-|a|^2/2t}}{\big|\det(1-\mathrm{Ad}(\gamma))
|_{\kz_0^\perp}\big|^{1/2}}
\frac{1}{\det(1-\mathrm{Ad}(k^{-1}))|_{\kp_0^\perp(\gamma)}}\\
&\quad \frac{1}{(2\pi t)^{p/2}}\cdot \tr^E\Big[ \rho^E(k^{-1})
\exp\Big(-t \Big(A+\frac{1}{48}\tr^{\kk_0}[C^{\kk_0,\kk_0}]
+\frac{1}{2}C^{\kk_0,E}\Big)\Big)\Big].
\end{split}
\end{equation}

Note that if $G$ is of real rank $1$, then $\kp_0$ is the 
vector subspace generated by $a$, so that \eqref{2.4.9} holds. 
Thus \eqref{2.4.10} recovers the result of Sally-Warner 
 \cite{SallyWarner}
where they assume that the real rank of $G$ is $1$.
\end{rema}

From \eqref{eq2.22} and \eqref{eq2.4.8}, we obtain %the precise 
a refined version of the Selberg trace formula 
for the Casimir operator :
\begin{equation}\label{eq2.4.10}
\tr[e^{-t\mL^Z_A}]=\sum_{[\gamma]\in [\Gamma]} 
\vol \Big(\Gamma\cap Z(\gamma)\backslash 
X(\gamma)\Big)  \tr^{[\gamma]}[e^{-t \mL^X_A}],
\end{equation}
and each term $\tr^{[\gamma]}[\cdot]$ is %evaluated precisely in
given by the closed formula \eqref{eq2.4.8}.

We give two examples here to explain the %precise
explicit version of 
the Selberg trace formula \eqref{eq2.4.10}.
 
\begin{exem}[Poisson summation formula]\label{e2.9}
Take $G=\R$ and $A=0$. Then $K=\{0\}$. 
We have $X=\R$
and $\mL^X_A= - \dfrac{1}{2}\Delta^{\R}
=-\dfrac{1}{2}\dfrac{\partial^2}{\partial x^2}$,
where $x$ is the coordinate on $\R$.
 Let $p_t(x,x')$ be the heat kernel associated with 
$e^{t\Delta^{\R}/2}$.

For $a\in \R$, we have $Z(a)=\R, \kk(a)=\{0\}$.
By (\ref{eq2.2.2}) or (\ref{eq2.3.15}), we have
\begin{align}\label{eq:2.4.11a}
\tr^{[a]}\Big[e^{-t\mL^X_A}\Big]= p_{t}(0,a).
\end{align}
From \eqref{eq2.4.8}, we get
\begin{equation}\label{eq2.4.11}
\tr^{[a]}\Big[e^{-t\mL^X_A}\Big]=\frac{1}{\sqrt{2\pi t}}
e^{-\frac{a ^2}{2t}}.
	\end{equation}
Thus %in this case, 
\eqref{eq2.4.8} gives simply an evaluation 
of the heat kernel on $\R$ which is well-known that
\begin{align}
\label{eq3.9.2}
p_t(x,x')=\frac{1}{\sqrt{2\pi t}}e^{-\frac{(x-x')^2}{2t}}.
\end{align}

%$\Z\backslash\R=\mathbb{S}^1$.	
Take $\Gamma=\Z\subset \R$, then  
$Z=\Z\backslash\R=\mathbb{S}^1$.
For any $\gamma\in \Gamma$, 
$X(\gamma)=Z(\gamma)/K(\gamma)=Z(\gamma)=\R$. 
Thus $\Gamma\cap Z(\gamma)\backslash Z(\gamma)
=\Z\backslash \R=\mathbb{S}^1$ and 
$\vol (\mathbb{S}^1)=1$. 	
% In the same time, we know that $e^{2\pi i n\theta}\;(n\in\Z)$ is 
% eigenfunction of $-\frac{1}{2}\dfrac{\partial^2}{\partial\theta^2}$ 
% with eigenvalue $2\pi^2 n^2$. Thus
% 	\begin{equation}
% 	\label{eq2.4.12}
% \tr[e^{-t\mL^X_A}]=\sum_{n\in\Z}e^{-2\pi^2 n^2 t}.
% 	\end{equation}	
The Selberg trace formula  \eqref{eq2.4.10} reduces 
to the Poisson summation formula: %for any $t>0$,
	\begin{equation}
	\label{eq2.4.15}
\sum_{k\in\Z} e^{-2\pi^2 k^2 t}=\sum_{k\in\Z}
\frac{1}{\sqrt{2\pi t}}e^{-\frac{k^2}{2t}}  \quad
\text{ for any } t>0.
	\end{equation}
\end{exem}
\begin{exem}\label{e2.10}
Let $G=\mathrm{SL}_2(\R)$ be the $2\times 2$ real 
special linear group with Lie algebra $\kg=\mathfrak{sl}_2(\R)$. 
The Cartan involution is given by 
$\theta:G\rightarrow G,\; g\mapsto {}^tg^{-1}$. 
Then $K=\mathrm{SO}(2)=\Big\{\begin{bmatrix}
	\cos{\beta}& \sin{\beta}\\
	-\sin{\beta}&\cos{\beta}
\end{bmatrix}\; :\; \beta\in\R\Big\}\simeq\mathbb{S}^1$
is the corresponding maximal compact subgroup and $X=G/K$
is the Poincar\'{e} upper half-plane 
$\mathbb{H}=\{z=x+iy\in \C\,:\, y>0, x\in \R\}$.     %\Im(z)>0\}$.
Precisely, an element $g=\begin{bmatrix}
 a&b\\
 c&d
\end{bmatrix}\in \mathrm{SL}_2(\R)$ acts on $\mathbb{H}$ by
\begin{equation}
 \label{eq2.4.16}
 g z=\frac{az+b}{cz+d}\in\mathbb{H}\quad \text{ for }
 z\in \mathbb{H}.
\end{equation}
	 
The Cartan decomposition of $\mathfrak{sl}_2(\R)$ is
\begin{equation}
\label{eq2.4.17}	
	 \kg=\kp\oplus\kk,
\end{equation}
where $\kk$ is the set of real antisymmetric matrices, 
and $\kp$ is the set of traceless symmetric matrices.	 
Let $B$ be the bilinear form on $\kg$ defined for $u,v\in \kg$ by
\begin{equation}\label{eq2.4.18}	
	 B(u,v)=2\tr^{\R^2}[uv].
\end{equation}
	 
Set 
\begin{equation} \label{eq2.4.19}
 e_1=\begin{bmatrix}
	 \frac{1}{2}&0\\
	 0&-\frac{1}{2}
	 \end{bmatrix},\;
\quad e_2=\begin{bmatrix}
	 0&\frac{1}{2}\\
	 \frac{1}{2}&0
	 \end{bmatrix},\;
\quad e_3=\begin{bmatrix}
	 0&\frac{1}{2}\\
	 -\frac{1}{2}&0
	 \end{bmatrix}.
\end{equation}
Then $\{e_1,e_2\}$ is a basis of $\kp$, and $e_3$ is a basis 
of $\kk$. They together form an orthonormal basis of 
the Euclidean space 
$(\kg,\langle\cdot,\cdot\rangle=-B(\cdot,\theta\cdot))$.
Moreover,we have the relations,
\begin{equation} \label{eq2.4.20}
 [e_1,e_2]=e_3,\;\quad [e_2,e_3]=-e_1,\;\quad [e_3,e_1]=-e_2.
 \end{equation}
	 
%Then we compute 
The metric on $X$ is given by $\frac{1}{y^{2}}(dx^2+ dy^{2})$.
The scalar curvature of $X$ is
\begin{equation}\label{eq2.4.21}
	 \tr^\kp[C^{\kk,\kp}]=-2|[e_1,e_2]|^2=-2.
 \end{equation}
Let $\Delta^X$ be the Bochner Laplacian acting on $C^\infty(X,\C)$. 
Then $\Delta^X= y^{2}(\frac{\partial^{2}}{\partial x^{2}}
+ \frac{\partial^{2}}{\partial y^{2}})$.
Since $\tr^\kk[C^{\kk,\kk}]=0$ here, we have on 
$C^\infty(X,\C)$,
\begin{equation} \label{eq2.4.22}
 \mL^X=\frac{1}{2}C^\kg
+\frac{1}{16}\tr^\kp [C^{\kk,\kp}]
+\frac{1}{48}\tr^{\kk}[C^{\kk, \kk}]
=-\frac{1}{2}\Delta^X-\frac{1}{8}.
\end{equation}
	 
From \eqref{eq2.4.20}, we see that a semisimple 
nonelliptic element $\gamma\in G$ is hyperbolic. 
Thus such $\gamma$ is conjugate to $e^{a e_1}$ 
with some $a \in\R\backslash \{0\}$. 
Note that the orbital integral depends only on 
the conjugacy class of $\gamma$ in $G$

If $\gamma=e^{a  e_1}$ with $a \in\R\backslash \{0\}$, 
then by \eqref{eq2.4.20},  
$\kk(\gamma)=0,\;\kz_0=\kz(\gamma)=\R e_1$, and
we have
\begin{equation}\label{eq2.4.2222}
\det(1-\mathrm{Ad}(\gamma))|_{\kz^\perp_0}
=-(e^{a /2}-e^{-a /2})^2.
\end{equation}
From Theorem \ref{t2.7}, \eqref{eq2.4.22} and \eqref{eq2.4.2222}, 
we get
\begin{equation}\label{eq2.4.23}
\tr^{[\gamma]}\Big[e^{t\Delta^X/2}\Big]
=\frac{1}{\sqrt{2\pi t}}\frac{\exp(-\frac{a ^2}{2t}-\frac{t}{8})}
{2\sinh(\frac{|a |}{2})}.
\end{equation}

For $Y^\kk_0=y_0e_3\in \kk$, the relations \eqref{eq2.4.20}
imply that
\begin{equation}\label{eq2.4.24}
\widehat{A}(i\mathrm{ad}(Y^\kk_0)|_\kp)
=\frac{y_0/2}{\sinh(y_0/2)}.
\end{equation}
From Theorem \ref{t2.7}, \eqref{2.4.8} and \eqref{eq2.4.24}, 
we get
\begin{equation}\label{eq2.4.25}
\tr^{[1]}[e^{t\Delta^X/2}] 
= \frac{e^{-t/8}}{2\pi t} \int_{\R} e^{-y_0^2/2t}
\frac{y_0/2}{\sinh(y_0/2)}\frac{dy_0}{\sqrt{2\pi t}}.
\end{equation}

% Recall the identity
% \begin{equation}
% \frac{1}{\sqrt{2\pi t}}e^{-y^2_0/2t}
% =\frac{1}{2\pi}\int_{\R} e^{-t\rho^2/2-i\rho y_0}d\rho.
% \end{equation} 
By taking the derivative with respect to $y_0$ in both sides of 
$\frac{1}{\sqrt{2\pi t}}e^{-y^2_0/2t}
=\frac{1}{2\pi}\int_{\R} e^{-t\rho^2/2-i\rho y_0}d\rho$, we get
\begin{equation}\label{eq2.4.26}
\frac{1}{\sqrt{2\pi t}} e^{-y^2_0/2t}\frac{y_0}{t}
=\frac{1}{2\pi}\int_{\R} e^{-t\rho^2/2}\rho \sin(\rho y_0)d\rho.
\end{equation}
Thus
\begin{equation}\label{eq2.4.27}
\begin{split}
\frac{1}{ t} \int_{\R} e^{-y_0^2/2t}
\frac{y_0/2}{\sinh(y_0/2)}\frac{dy_0}{\sqrt{2\pi t}}&= 
\frac{1}{4\pi}\int_{\R} e^{-t\rho^2/2}\rho \left(
\int_{-\infty}^{+\infty} \frac{\sin(\rho y_0)}
{\sinh(y_0/2)}dy_0\right)d\rho\\
&=\frac{1}{2}\int_{\R} e^{-t\rho^2/2} 
\rho \tanh(\pi \rho)d\rho,
\end{split}
\end{equation}
where we use the identity 
$\int_{-\infty}^{+\infty} \frac{\sin(\rho y_0)}{\sinh(y_0/2)}dy_0
=2\pi \tanh(\pi \rho)$.

Let $\Gamma\subset \mathrm{SL}_2(\R)$ be 
a discrete torsion-free 
cocompact subgroup. 
Then $Z=\Gamma\backslash X$ is a compact Riemann surface. 
We say that $\gamma\in\Gamma$ is primitive if there does not exist
$\beta\in\Gamma$ and $k\in\N,\;k\geq 2$ such that
$\gamma=\beta^k$.

If $\gamma=e^{a  e_1}\in \Gamma$ is primitive, 
then $|a |$ is the length of the corresponding closed geodesic
in $Z$ and for any $k\in \Z,\; k\neq 0$, 
$Z(\gamma^k)=Z(\gamma)= e^{\R e_1}$, and moreover,
\begin{equation}\label{eq2.4.29}
\vol (Z(\gamma^k)\cap \Gamma\backslash Z(\gamma^k))
=|a |.
\end{equation}

Thus by \eqref{eq2.4.10}, \eqref{eq2.4.23}, \eqref{eq2.4.25}, 
\eqref{eq2.4.27} and \eqref{eq2.4.29}, we get
\begin{equation}\label{eq2.4.30}
\begin{split}
\tr[e^{t\Delta^Z/2}] &=\sum_{\substack{\gamma\in \Gamma \;
\mathrm{primitive},\\  [\gamma]=[e^{a  e_1}],\;a \neq 0}} 
|a |\sum_{k \in \N,\;k \neq 0} \tr^{[e^{k a  e_1}]}
[e^{t\Delta^X/2}]+\vol (Z)\tr^{[1]}[e^{t\Delta^X/2}]\\
&=\sum_{\substack{\gamma\in \Gamma \;\mathrm{primitive},\\  
[\gamma]=[e^{a  e_1}],\;a \neq 0}} |a |
\sum_{k \in \N,\;k \neq 0} \frac{1}{\sqrt{2\pi t}} 
\frac{1}{2\sinh(\frac{k |a |}{2})}
e^{-\frac{k ^2a ^2}{2t}-\frac{t}{8}}\\
&\hspace{10mm}+\frac{\vol (Z)}{4\pi}e^{-t/8}
\int_{\R} e^{-t\rho^2/2}\rho\tanh(\pi \rho)d\rho.
\end{split}
\end{equation}
Formula \eqref{eq2.4.30} is exactly the original Selberg trace 
formula in \cite[(3.2)]{Selberg56} 
(cf. also \cite[p. 233]{McKean72}).
\end{exem}

%%%%%%%%%%%%%%%%%%%%%%%%%%%%%
\subsection{Harish-Chandra's Plancherel Theory}\label{s2.6}
In this subsection, we briefly describe  
Harish-Chandra's approach to orbital integrals. 
This approach can be used to evaluate the orbital integrals
of arbitrary test function, %however only 
for sufficiently
regular semisimple elements.  This formula  contains %certain 
complicated expressions involving infinite sums which do not
converge absolutely, and have no obvious closed form except 
for some special groups.  
An useful reference on 
Harish-Chandra's work on orbital integrals
is Varadarajan's book \cite{Vara77}.

Recall that $G$ is a connected reductive group. 
Denote by $G'\subset G$ the space of regular elements.
Let $C^\infty_c(G)$ be the vector space of smooth functions
with compact support on $G$.
For $f\in C^\infty_c(G)$, attached to each $\theta$-invariant 
Cartan subgroup 
$H$ of $G$, Harish-Chandra introduce a smooth function
${}^\prime\! F_f^H$ (cf. \cite[\S 17]{HCRI}), 
as  an orbital integral %Fourier transformation 
of $f$ in a certain sense, defined on 
$H\cap G'$, which has reasonable limiting behavior on 
the singular set in $H$. 

Let $\gamma$ be a semisimple element such that \eqref{eq2.3.4} 
holds. If $\gamma$ is regular, then up to conjugation there exists 
a unique $\theta$-invariant Cartan subgroup $H$ which contains 
$\gamma$. 
In this case, ${}^\prime\! F_f^H(\gamma)$ is equal to a product of  
$\tr^{[\gamma]}[f]$ and  an explicit Lefschetz like 
denominator of $\gamma$. Now if $\gamma$ is a singular semisimple 
element,  let $H$ be the unique (up to conjugation) 
$\theta$-invariant  Cartan subgroup with maximal compact 
dimension, which contains $\gamma$. 
Following Harish-Chandra \cite{HCD}, there is an explicit  
differential operator 
$D$ defined on  $H$ such that 
%$\tr^{[\gamma]}[f]$ can be obtained from the limit of the function  
%$D\,{}^\prime\! F_f^H(\gamma')$ when $\gamma'\to \gamma$. 
\begin{align}\label{eq:2.6.1}
\tr^{[\gamma]}[f]= \lim_{\gamma'\in H\cap G'\to \gamma}D\,
{}^\prime\! F_f^H(\gamma').
\end{align}
Thus, to determine the orbital integral  $\tr^{[\gamma]}[f]$,
it is enough to  calculate ${}^\prime\! F_f^H$ on the regular set 
$H\cap G'$.  

Take $\gamma\in H\cap G'$ a regular element in $H$. 
Harish-Chandra developed certain techniques  to calculate 
${}^\prime\! F_f^H$, obtaining formulas which are known as 
{\it Fourier inverse formula}. 
Indeed, $f\in C_c^\infty(G)\to {}^\prime\! F_f^H(\gamma)$ defines
an invariant distribution on $G$. The idea is to write 
${}^{\prime}\!F_f^H(\gamma)$ as a combination of  
invariant eigendistributions  
(i.e.,  a distribution on $G$
which is invariant under the adjoint action of  $G$,
and which is an eigenvector of the center of $U(\kg)$), 
like the global character of the discrete series representations 
and the unitary principal series representations of $G$, 
as well as certain singular invariant eigendistributions. 
More precisely, let $H=H_IH_R$ be Cartan decomposition of 
$H$ (cf. \cite[\S 8]{HCRI}), where $H_I$ is a compact 
Abelian group and $H_R$ is a vector space.  
Denote by $\widehat{H}, \widehat{H}_I,\widehat{H}_R$  
the set of irreducible unitary representations  of $H,H_I,H_R$. 
Then $\widehat{H}= \widehat{H}_I \times \widehat{H}_R$. 
Following \cite{HCL,HerbSally1979}, for 
$a^*=(a_I^*,a_R^*)\in \widehat{H}$, we can associate 
an  invariant eigendistribution $\Theta^{H}_{a^*}$ on $G$. 
Note that if $H$ is compact and if $a^*_I$ is regular, 
then $\Theta^H_{a^*}$ is the global  character of 
the discrete series representations of $G$, and that
if $H$ is noncompact and if $a^*_I$ is regular, then 
$\Theta^H_{a^*}$ is the global  character of 
the unitary principal series representations of $G$. 
When $a_I^*$ is singular, $\Theta^H_{a^*}$ is much more 
complicated. It is an alternating sum of some unitary characters, 
which in general are  reducible.

In \cite{HCL}, Harish-Chandra announced the following Theorem. 
%%%%%%%%%%%%%%%%%%%%%%%%%
\begin{theo}\label{tHCr}\cite[Theorem 15]{HCL}.
Let $\{H_1,\cdots, H_l\}$ be the complete set of non conjugated 
$\theta$-invariant Cartan subgroups of $G$. 
Then there exist computable continuous functions $\Phi_{ij}$ 
on $H_i\times \widehat{H}_j$ such that for any regular element 
$\gamma\in H_i\cap G'$, 
 \begin{align}
 \label{eqHCr}
 {}^\prime \!F_f^{H_i}(\gamma)
 =\sum_{j=1}^l\sum_{a^*_I\in \widehat{H}_{jI}}
 \int_{a^*_R\in \widehat{H}_{jR}}\Phi_{ij}(\gamma,a_I^*,a_R^*)
 \Theta^{H_j}_{a^*}(f)da_R^*.
 \end{align}
\end{theo}
In \cite{HCL}, Harish-Chandra only explained the idea of 
%the proof which is obtained 
a proof by induction on $\dim G$. 
A more explicit version is obtained by 
Sally-Warner \cite{SallyWarner} when $G$ is of real rank one 
(cf. Remark \ref{t2.9}), and by Herb \cite{Herb79} 
(cf. also Bouaziz \cite{Bou95}) for general $G$. 
However,   Herb's formula only holds for $\gamma$ in an open 
dense subset of $H_i\cap G'$ and involves  certain infinite sum 
of integrals which converges, but cannot be directly differentiated, 
term by term. In particular, the orbital integral of singular 
semisimple elements could not be obtained from Herb's formula 
by applying term by term the differential operator $D$
in (\ref{eq:2.6.1}). 
%In general, we do not know 
%a closed formula for  $\tr^{[\gamma]}[f]$ for singular $\gamma$. 
When $\gamma=1$, much more is known:
\begin{theo}[Harish-Chandra \cite{HC}]\label{thm:HC}
There exists computable real analytic elementary functions 
$p^{H_j}(a^*)$ defined on $\widehat{H}_j$ such that
for $f\in C^\infty_c(G)$, 
we have 
\begin{align}\label{eq:HC2}
\tr^{[1]}[f]=f(1)=\sum_{j=1}^l\sum_{a^*_I\in \widehat{H}_{jI}, 
\mathrm{regular}}
\int_{a^*_R\in \widehat{H}_{jR}}
\Theta^{H_j}_{a^*}(f)p^{H_j}(a_I^*,a_R^*)da_R^*.
\end{align}
\end{theo}

Theorem \ref{thm:HC} can be applied to more general functions 
such as  Harish-Chandra Schwartz functions, e.g., the trace of  
the heat kernel $q_t\in C^\infty(G,\End(E))$ of 
$e^{-t\mathcal{L}_A^X}$.  Thus,
\begin{align}\label{eq:HC3}
\tr^{[1]}\Big[e^{-t\mathcal{L}_A^X}\Big]
=\sum_{j=1}^l\sum_{a^*_I \in \widehat{H}_{jI},\mathrm{regular}}
\int_{a^*_{R}\in\widehat{H}_{jR}}
\Theta^{H_j}_{a^*}\left(\tr^E[q_t]\right)
p^{H_j}(a_I^*,a_R^*)da_R^*.
\end{align}
For $H_j$, we can associate a cuspidal parabolic subgroup $P_j$ 
with Langlands decomposition $P_j=M_jH_{jR}N_j$ such that 
$H_{jI}\subset M_j$ is a compact Cartan subgroup of $M_j$.  
For $a^*=(a^*_I, a^*_R )\in \widehat{H}_j$ with $a^*_I$ regular, 
denote by $(\varsigma_{a^*_I},V_{a^*_I})$ the discrete series 
representations of $M_j$ associated to $a^*_I$, and  
denote by $(\pi_{a^*},V_{a^*})$ the associated principal series 
representations of $G$ associated to $\varsigma_{a^*_I}$ and 
$a^*_R$. We have 
\begin{align}
\Theta^{H_j}_{a^*}\left(\tr^E[q_t]\right)
=\tr^{V_{{a^*}}\otimes E}\left[\pi_{a^*}(q_t)\right]
\, \text{ with } \pi_{a^*}(q_t)= \int_{G} q_{t}(g) \pi_{a^*}(g) dg.
%\in \End(V_{{a^*}}\otimes E).
\end{align}

It is not difficult to see that  the image of the operator 
$\pi_{a^*}(q_t)$ is 
$(V_{a^*}\otimes E)^K
\simeq (V_{{a^*_I}}\otimes E)^{K\cap M_j}$, and 
$\pi_{a^*}(q_t)$ acts as 
$e^{-t(\frac{1}{2}C^{\kg,\pi_{a^*}}
+\frac{1}{16}\tr^\kp[C^{\kk,\kp}]+\frac{1}{48}\tr^{\kk}[C^{\kk,
\kk}]+A)}$ on its image. We get 
\begin{align}\label{eq:HC4}
\tr^{V_{a^*}\otimes E}\left[\pi_{a^*}(q_t)\right]
= e^{-t(\frac{1}{2}C^{\kg,\pi_{a^*}}
+\frac{1}{16}\tr^\kp[C^{\kk,\kp}]+\frac{1}{48}\tr^{\kk}[C^{\kk,
\kk}])}\tr^{(V_{{a^*_I}}\otimes E)^{K\cap M_j}}[e^{-tA}].
\end{align}
Thus,
\begin{multline}\label{eq:HC}
\tr^{[1]}\Big[e^{-t\mathcal{L}_A^X}\Big]
=\sum_{j=1}^l\sum_{a^*_I \in \widehat{H}_{jI},\mathrm{regular}}
\int_{a^*_{R}\in\widehat{H}_{jR}}e^{-t(\frac{1}{2}C^{\kg,\pi_{a^*}}
+\frac{1}{16}\tr^\kp[C^{\kk,\kp}]+\frac{1}{48}\tr^{\kk}[C^{\kk,
\kk}])}\\
\tr^{(V_{a^*_I}\otimes E)^{K\cap M_j}}
[e^{-tA}]p^{H_j}(a^*_I,a_R^*)da_R^*.
\end{multline}

\begin{rema}\label{r2.6.3}
Equation \eqref{eq:HC} is not as explicit  as \eqref{eq2.4.8}, 
because in general it is not easy to determine
all parabolic subgroups, 
all the discrete series of $M$, and  the Plancherel densities 
$p^{H_j}(a^*)$. 

%From these descriptions, we hope the readers have
We hope that from these descriptions, the readers got an idea on 
Harish-Chandra's Plancherel theory  %formula 
as an algorithm to compute orbital integrals. %Contrast to 
These results  use the full force %stress
of the unitary representation theory (harmonic analysis)
of reductive Lie groups, both at the technical %technique level 
and the representation level. %of the results.

Bismut's explicit formula of the orbital integrals
associated with the Casimir operator gives a closed 
formula in full generality for any semisimple element and any
reductive Lie group. 
%A part from the normal coordinate system on $X$
%based at $X(\gamma)$ (cf. Theorems \ref{t2.5}, \ref{t2.6}),
Bismut avoided completely the use of the harmonic analysis
on reductive Lie groups. The hypoelliptic deformation
allows him to localize the orbital integral for $\gamma$
to any neighborhood of the family of shortest geodesic 
associated with $\gamma$, i.e., $X(\gamma)$.

There is a mysterious connection between
Harish-Chandra's Plancherel theory  and Theorem \ref{t0.1}: 
in Harish-Chandra's Plancherel theory,
the integral are taken  on the $\kp$ part,
but in Theorem \ref{t0.1}, the integral is on the $\kk$ part.
In particular, in Example \ref{e2.10} for $G= \SL_{2}(\R)$,
we obtain the contribution $\int_{\R}e^{-t\rho^2/2}\rho 
\tanh (\pi \rho) d\rho$ from the Plancherel theory
for $\gamma=1$. This coincide with (\ref{eq2.4.25})
by using a Fourier transformation argument as explained in 
(\ref{eq2.4.27}).
\end{rema}

\begin{rema}\label{r2.6.4}
Assume $G=\mathrm{SO}^0(m,1)$ with $m$ odd. 
There exists only one Cartan subgroup $H$,
and $p^H(a^*_I,\cdot)$ is an explicit polynomial. 
In this case, \eqref{eq:HC} becomes completely explicit. 
%which is the essential ingredient in \cite{MP13}.
\end{rema}

%%%%%%%%%%%%%%%%%%%%%%%%%%%%%%%%
\section{Geometric hypoelliptic operator and 
dynamical  systems}\label{s3}

% In this section, we describe some basic ideas of Bismut's program 
% on his geometric hypoelliptic operators and its applications 
% in geometry and dynamical system.

In this section, we explain how to construct geometrically the 
hypoelliptic Laplacians for a symmetric space, %in hoping  
with the goal to prove 
Theorem \ref{t0.1} in the spirit of the heat kernel proof of the 
Lefschetz fixed-point formula (cf. Section \ref{s1.2}).
We introduce a hypoelliptic version of the orbital integral
that depends on $b$. %and the corresponding 
The analogue of the methods of local index theory
are needed to evaluate the limit.
%a hypoelliptic version of the McKean-Singer formula,
Theorem \ref{t3.6.4} identifies the orbital integral associated 
with the Casimir operator to the hypoelliptic orbital integral
for the parameter $b>0$.
As $b\to +\infty$, the hypoelliptic orbital integral localizes
near $X(\gamma)$.

This section is organized as follows. In Section \ref{s3.3}, 
we explain how to compute %realize
the cohomology of a vector space
by using algebraic de Rham complex and 
%then we introduce its analytic counterpart, the 
its Bargmann transformation, whose Hodge Laplacian
is a harmonic oscillator. 
%In Section \ref{s3.4}, we describe the ideas needed in the 
%construction of the geometric hypoelliptic Laplacian in order to 
%establish Theorem \ref{t2.7}.
In Section \ref{s3.5}, 
we recall the construction of the Dirac operator of Kostant,
and in Sections \ref{s3.4a},
we construct the geometric hypoelliptic Laplacian %realize the ideas 
by combining the constructions in Sections \ref{s3.3} and \ref{s3.5}.
In Section \ref{s3.6}, we introduce the hypoelliptic orbital integrals 
and a hypoelliptic version of the McKean-Singer formula
for these  %the elliptic and hypoelliptic
orbital integrals. In Section \ref{s3.7}, we describe the limit 
of the hypoelliptic orbital integrals as $b\to +\infty$.
Finally, in Section \ref{s3.8}, we explain %briefly
some relations of the hypoelliptic heat equation 
to the wave equation on the base manifold, which 
plays an important role in the proof of uniform
Gaussian-like estimates for the hypoelliptic heat kernel.

%%%%%%%%%%%%%%%%%%%%%%%%%%%%%%%%%%%%
\subsection{Cohomology of a vector space and harmonic oscillator}
\label{s3.3}
Let $V$ be a real vector space of dimension $n$, and
let $V^*$ be its dual. Let $Y$ be the tautological section of $V$ 
over $V$. Then $Y$ can be identified with the corresponding radial 
vector field. Let $d^V$ denote the de Rham operator.
%Let $(\Omega^\bullet(V),d^V)$ be the de Rham complex of $V$.

Let $L_Y$ be the Lie derivative associated with $Y$,
and let $i_Y$ be the contraction of $Y$.
By Cartan's formula, we have the identity 
%of operators acting on $\Omega^\bullet(V)$,
\begin{equation}
\label{eq3.3.1}
L_Y=[d^V, i_Y].
\end{equation}

% We work first in algebraic way to compute the cohomology of $V$,
% then we will interpret it in analytic way. 
% Thus instead of working on smooth forms on $V$, 
% we work on the polynomial forms on $V$.

Let $S^\bullet(V^*)=\oplus_{j=0}^\infty S^j(V^*)$ be the symmetric
algebra of $V^*$, which can be canonically identified %naturally as 
with the polynomial algebra of $V$. 
Then $\Lambda^\bullet(V^*)\otimes S^\bullet(V^*)$
is the vector space %subspace of $\Omega^\bullet(V)$ 
of polynomial forms on $V$. Let $N^{S^\bullet(V^*)},\; 
N^{\Lambda^\bullet(V^*)}$ be the number operators on 
$S^\bullet(V^*),\; \Lambda^\bullet(V^*)$, 
which act by multiplication by 
$k$ on $S^k(V^*),\; \Lambda^k(V^*)$. Then
\begin{equation}\label{eq3.3.2}
L_Y|_{\Lambda^\bullet(V^*)\otimes S^\bullet(V^*)}
	=N^{S^\bullet(V^*)}+N^{\Lambda^\bullet(V^*)}.
\end{equation}
By (\ref{eq3.3.1}) and  (\ref{eq3.3.2}), 
the cohomology of %$V$ via 
the polynomial forms 
$(\Lambda^\bullet(V^*)\otimes S^\bullet(V^*), d^V)$
on $V$ is equal to  $\R 1$.

Assume that $V$ is equipped with a scalar product. Then %each 
$\Lambda^\bullet(V^*)$, $S^\bullet(V^*)$ inherit
 associated scalar products. %We normalized it so that 
 For instance, if $V=\R$, then
$\norm{1^{\otimes j}}^2=j!$. 
%We equip $S^\bullet(V^*)$ 
%with the direct sum of the scalar products on the $S^j(V^*)$.
With respect to this scalar product on 
$\Lambda^\bullet(V^*)\otimes S^\bullet(V^*)$, $i_Y$
is the adjoint of $d^V$. Therefore $L_Y$ is the associated Hodge
Laplacian on $\Lambda^\bullet(V^*)\otimes S^\bullet(V^*)$.
Remarkably enough, it does not depend on $g^V$. 
By \eqref{eq3.3.2}, we get
\begin{equation}
\label{eq3.3.3}
	\Ker (L_Y)=\R 1.
\end{equation}
We have given a Hodge theoretic interpretation to the proof that
the cohomology of the complex of polynomial forms %$V$
is concentrated in degree $0$.

Let $\Delta^V$ denote the (negative) Laplacian on $V$. 
Let $L_2(V)$ be the corresponding Hilbert space of 
square integrable real-valued functions on $V$.
\begin{defi}\label{d3.3}
Let $T: S^\bullet(V^*)\rightarrow L_2(V)$ be the map such that  
given $P\in S^\bullet(V^*)$, then
\begin{equation}	\label{eq3.3.4}
(TP)(Y)=\pi^{-n/4} e^{-\frac{\abs{Y}^2}{2}}(e^{-\Delta^V/2}P)
(\sqrt{2}Y).
\end{equation}
\end{defi}
Since $P$ is a polynomial, $e^{-\Delta^V/2}P$ 
is defined by taking the obvious formal expansion of 
$e^{-\Delta^V/2}$. Its inverse, the Bargmann kernel, is given by 
\begin{equation}\label{eq3.3.5}
(Bf)(Y)=\pi^{n/4} e^{\Delta^V/2} \Big( e^{\abs{Y}^2/4}
f(\frac{Y}{\sqrt{2}})\Big).
\end{equation}
Here the operator $e^{\Delta^V/2}$ is defined via 
the standard heat kernel of $V$.

Set
\begin{equation}
\label{eq3.3.6}
 \overline{d}=T d^V B,\quad \overline{d}^*=T i_Y B: 
 \Lambda^\bullet(V^*)\otimes L_2(V)\rightarrow 
 \Lambda^\bullet(V^*)\otimes L_2(V).
\end{equation}
Then by \eqref{eq3.3.4} and \eqref{eq3.3.5}, we get
\begin{equation}\label{eq3.3.7}
\overline{d}=\frac{1}{\sqrt{2}}(d^V+Y^*\wedge)\;,\quad
\overline{d}^*=\frac{1}{\sqrt{2}}(d^{V*}+i_Y).
\end{equation}
Here $Y^*$ is the metric dual of $Y$ in $V^*$, 
and $d^{V*}$ is the usual formal $L_{2}$ adjoint of $d^V$.

Let $\{e_j\}$ be an orthonormal basis of $V$ and let 
$\{e^j\}$ be its dual basis. For $U\in V$, let $\nabla_U$
be the usual differential along the vector $U$. 
Put $Y=\sum^n_{j=1} Y_j e_j$, then
\begin{equation}\label{eq3.3.8}
\begin{split}
&d^V=\sum^n_{j=1} e^j\wedge \nabla_{e_j},\; \quad
d^{V*}=-\sum^n_{j=1} i_{e_j}\nabla_{e_j};\\
&Y^*\wedge=\sum^n_{j=1}Y_j e^j\wedge,\; \quad
i_Y=\sum^n_{j=1} Y_j i_{e_j}.
\end{split}
\end{equation}
From \eqref{eq3.3.7} and \eqref{eq3.3.8}, we get
\begin{equation}\label{eq3.3.9}
\overline{d}^2=(\overline{d}^{*})^2=0,\quad
T L_{Y} T^{-1} = [\overline{d},\overline{d}^*]
=\frac{1}{2}\Big(-\Delta^V+\abs{Y}^2-n\Big)
+N^{\Lambda^\bullet(V^*)}.
\end{equation}
Note that $\frac{1}{2}(-\Delta^V+\abs{Y}^2-n)$ 
is the harmonic oscillator on $V$ already appeared in (\ref{3.3}).
%By \eqref{eq3.3.3}, the kernel of the operator 
%$[\overline{d},\overline{d}^*]$ in  
%$\Lambda^\bullet(V^*)\otimes S^\bullet(V^*)$
In \eqref{eq3.3.3}, we saw that the kernel of $[d^V,i_Y]$ in
$\Lambda^\bullet(V^*)\otimes S^\bullet(V^*)$
is generated by 
$1$ and so it is $1$-dimensional and is concentrated 
in total degree $0$. Equivalently the kernel of 
the unbounded operator $[\overline{d},\overline{d}^*]$ 
acting on $\Lambda^\bullet(V^*)\otimes L_2(V)$ 
is $1$-dimensional and is generated by the function 
$e^{-\abs{Y}^2/2}/{\pi^{n/4}}$.

 %%%%%%%%%%%%%%%%%%%%%%%%%%%%%%%%%%
\subsection{The Dirac operator of Kostant}
\label{s3.5}
Let $V$ be a finite dimensional real vector space of dimension $n$
 and let $B$ be a real valued symmetric bilinear form on $V$.

Let $c(V)$ be the Clifford algebra associated to $(V,B)$. Namely,
 $c(V)$ is the algebra generated over $\R$ by $1,u\in V$ and the
  commutation relations for $u,v\in V$,
\begin{equation}
\label{eq3.5.1}
 uv+vu=-2B(u,v).
\end{equation}
We denote by $\widehat{c}(V)$ the Clifford algebra associated to 
$-B$. Then $c(V),\widehat{c}(V)$ are filtered by length, 
and their corresponding $\mathrm{Gr}^\cdot$ is just 
$\Lambda^\bullet (V)$. Also
 they are $\Z_2$-graded by length.

In the sequel, we assume that $B$ is nondegenerate. 
Let $\varphi: V\rightarrow V^*$ be the isomorphism such that 
if $u,v\in V$,
\begin{equation}
\label{eq3.5.2}
	(\varphi u,v)=B(u,v).
\end{equation}

If $u\in V$, let $c(u),\widehat{c}(u)$ act on 
$\Lambda^\bullet(V^*)$ by 
\begin{equation}\label{eq3.5.3}
c(u)=\varphi u\wedge - i_u,\qquad
\widehat{c}(u)=\varphi u\wedge+ i_u.
\end{equation}
Here $i_u$ is the contraction operator by $u$.

%From the supercommutator relation 
Using supercommutators as in (\ref{0.7}), 
from \eqref{eq3.5.3}, we find  %check easily 
that for $u,v\in V$,
\begin{equation}
\label{eq3.5.5}
[c(u),c(v)]=-2B(u,v),\quad [\widehat{c}(u),\widehat{c}(v)]=2B(u,v),
\quad [c(u),\widehat{c}(v)]=0.
\end{equation}
Equation (\ref{eq3.5.5}) shows that $c(\cdot), \widehat{c}(\cdot)$
are representations of the Clifford algebras 
$c(V)$, $\widehat{c}(V)$ on $\Lambda^\bullet(V^*)$.

We will apply now the above constructions to the vector space 
$(\kg, B)$ of Section \ref{s2.1}.

If $\{e_i\}^{m+n}_{i=1}$ is a basis of $\kg$, we denote by
$\{e^*_i\}^{m+n}_{i=1}$ its dual basis of $\kg$ with respect 
to $B$ (i.e., $B(e_{i},e_{j}^*)=\delta_{ij}$), 
and by $\{e^i\}^{m+n}_{i=1}$ the dual basis of 
$\kg^*$.

Let $\kappa^\kg\in \Lambda^3(\kg^*)$ be such that if 
$a,b,c\in\kg$,
\begin{equation}
\label{eq3.5.6}
	\kappa^\kg(a,b,c)=B([a,b],c).
\end{equation}
Let $\widehat{c}(\kappa^{\kg})\in \widehat{c}(V)$ correspond to 
$\kappa^\kg\in \Lambda^3(\kg^*)$ defined by
\begin{equation}
\label{eq3.5.7}
\widehat{c}(\kappa^\kg)=\frac{1}{6} 
\kappa^\kg(e_i^*,e_j^*,e_k^*) \,
\widehat{c}(e_i)\,\widehat{c}(e_j)\,\widehat{c}(e_k).
\end{equation}

\begin{defi}\label{d3.6}
Let $\widehat{D}^\kg\in \widehat{c}(\kg)\otimes U(\kg)$ be the 
Dirac operator
\begin{equation}\label{eq3.5.8}
\widehat{D}^\kg=\sum^{m+n}_{i=1} \widehat{c}(e_i^*) e_i
-\frac{1}{2}\widehat{c}(\kappa^\kg).
\end{equation}
\end{defi}
Note that $\widehat{c}(\kappa^\kg), \widehat{D}^\kg$ are 
$G$-invariant.
%The operator $\widehat{D}^\kg$ is called 
%the Dirac operator, %of Kostant,
The operator $\widehat{D}^\kg$
acts naturally on $C^\infty(G,\Lambda^\bullet(\kg^*))$.

\begin{theo}[Kostant formula, 
   {\cite{Ko99}, \cite[Theorem 2.7.2, (2.6.11)]{B11b}}]
\label{t3.7}
\begin{equation}	\label{eq3.5.9}
\widehat{D}^{\kg,2}=-C^\kg-\frac{1}{8}\tr^\kp[C^{\kk,\kp}]-
\frac{1}{24}\tr^\kk[C^{\kk,\kk}].
\end{equation}
\end{theo}

%%%%%%%%%%%%%%%%%%%%%%%%%%%%%%%%%%%
\subsection{Construction of geometric hypoelliptic operators}
\label{s3.4a}
The operator $\widehat{D}^\kg$ acts naturally on 
$C^\infty(G,\Lambda^\bullet(\kg^*))$ and also on
$C^\infty(G,  \Lambda^\bullet(\kg^*)\otimes S^\bullet(\kg^*))$. 
As we saw in Section \ref{s3.3}, 
from a cohomological point of view, 
$\Lambda^\bullet(\kg^*)\otimes S^\bullet(\kg^*)\simeq \R$.
This is how ultimately $C^\infty(G,\R)$ (and 
$C^\infty(X,\R)$)  will reappear.

% As the Dirac operator of Kostant acts naturally only on 
% $C^\infty(G,\Lambda^\bullet(\kg^*))$,
% it can act neither on $C^\infty(G,\R)$ (which could descend to 
% $C^\infty(X,\R)$) nor on $C^\infty(G,\Lambda^\bullet(\kp^*))$ 
% (which descends to $\Omega^\bullet(X)$). 
% In view of the consideration of the algebraic de Rham theory 
% in Section \ref{s3.3}, we need to work on the fiber 
% $\kg^*$ or $\kg$, roughly on the space
% $C^\infty(G,  \Lambda^\bullet(\kg^*)\otimes S^\bullet(\kg^*))$,
% %$C^\infty(G,  L_{2}(\kg,\Lambda^\bullet(\kg^*)))$, 
% and we use the operators $\overline{d},\; \overline{d}^*$ 
% in Section \ref{s3.3}.

We denote by $\Delta^{\kp\oplus \kk}$ the standard Euclidean 
Laplacian on the Euclidean vector space $\kg=\kp\oplus \kk$. 
If $Y\in \kg$, we split $Y$ in the form 
\begin{align}\label{eq:3.5O}
Y= Y^\kp +Y^\kk \quad \text{ with }
Y^\kp \in \kp, Y^\kk \in \kk.
\end{align}
If $U\in \kg$, we use the notation
\begin{equation}
\label{eq3.5.16}
\begin{split}
\widehat{c}(\mathrm{ad}(U))&=-\frac{1}{4}B([U,e^*_i],e^*_j) 
\,\widehat{c}(e_i)\,\widehat{c}(e_j),\\
 c(\mathrm{ad}(U))&=\frac{1}{4}B([U,e^*_i],e^*_j) 
\,c(e_i)\,c(e_j).
\end{split}
\end{equation}

%Here is the final operator which
Here is the  operator $\mathfrak{D}_b$
appeared in \cite[Definition 2.9.1]{B11b}
which acts on 
\begin{align}\label{eq:3.51}
C^\infty(G,\Lambda^\bullet(\kg^*)\otimes S^\bullet(\kg^*) )
``\simeq'' C^\infty(G\times \kg,\Lambda^\bullet(\kg^*)).
\end{align}
\begin{defi}\label{d3.9}
	Set
\begin{equation}\label{eq3.5.15}
\mathfrak{D}_b=\widehat{D}^\kg+ic\Big([Y^\kk,Y^\kp]\Big)
+\frac{\sqrt{2}}{b}\Big(\overline{d}^\kp
-i\overline{d}^\kk+\overline{d}^{\kp*}+i\overline{d}^{\kk*}\Big).
\end{equation}
\end{defi}

% From Section \ref{s3.3}, the kernel $H\subset 
% L^2(\kg,\Lambda^\bullet(\kg^*))$ of the operator 
% $\overline{d}^\kp-i\overline{d}^\kk+\overline{d}^{\kp*}
% +i\overline{d}^{\kk*}$ is $1$-dimensional and spanned by 
% $e^{-\abs{Y}^2/2}$.
The introduction of $i$ in the third term in the right-hand side 
of (\ref{eq3.5.15}) is made so that its principal symbol 
anticommutes with the principal symbol of $\widehat{D}^\kg$.

Let $\{e_{j}\}_{j=1}^m$ be an orthonormal basis of $\kp$, and let 
$\{e_{j}\}_{j=m+1}^{m+n}$ be an orthonormal basis of $\kk$.
If $U\in\kk$, $\mathrm{ad}(U)|_\kp$ acts as an antisymmetric 
endomorphism of $\kp$ and by \eqref{eq3.5.16}, we have
\begin{equation}\label{eq3.5.17}
c(\mathrm{ad}(U)|_\kp)=\frac{1}{4}\sum_{1\leq i,j\leq m}
\langle [U,e_i],e_j\rangle c(e_i)c(e_j).
\end{equation}
Finally, if $v\in \kp$, $\mathrm{ad}(v)$ exchanges $\kk$ and 
$\kp$ and is antisymmetric with respect to $B$, i.e., it is symmetric 
with respect to the scalar product on $\kg$. 
Moreover, by \eqref{eq3.5.16}
\begin{equation}\label{eq3.5.18}
c(\mathrm{ad}(v))= - \frac{1}{2}
\sum_{\substack{m+1\leq i\leq m+n \\ 1\leq j\leq m}} 
\langle [v,e_i], e_j\rangle c(e_i)c(e_j).
\end{equation}

If $v\in \kg$, we denote by $\nabla^V_{v}$ the corresponding 
differential operator along $\kg$. In particular,
$\nabla^V_{[Y^\kk,Y^\kp]}$ denotes the differentiation
operator in the direction $[Y^\kk,Y^\kp]\in \kp$.
If $Y\in\kg$, we denote by 
$\underline{Y}^\kp+ i\underline{Y}^\kk$ 
the section of $U(\kg)\otimes_\R \C$ associated with
$Y^\kp+i Y^\kk\in\kg\otimes_\R \C$. 
\begin{theo}\label{t3.10}\cite[Theorem 2.11.1]{B11b}
The following identity holds:
\begin{equation}	\label{eq3.5.19}
\begin{split}
&\frac{\mathfrak{D}^2_b}{2}=\frac{\widehat{D}^{\kg,2}}{2}
+\frac{1}{2}\Big|[Y^\kk,Y^\kp]\Big|^2
+\frac{1}{2b^2} \Big(-\Delta^{\kp\oplus \kk}+\abs{Y}^2-m-n\Big)
+\frac{N^{\Lambda^\bullet(\kg^*)}}{b^2}\\
&+\frac{1}{b}\bigg( \underline{Y}^\kp
+i\underline{Y}^\kk-i\nabla^V_{[Y^\kk,Y^\kp]} 
+\widehat{c}(\mathrm{ad}(Y^\kp+iY^\kk))
 +2ic(\mathrm{ad}(Y^\kk)|_\kp)
-c(\mathrm{ad}(Y^\kp))\bigg).
\end{split}
\end{equation}
\end{theo}		
	
By \eqref{eq2.3},
\begin{equation}
\label{eq3.5.25}
G\times_K\kg=TX\oplus N,\;\mathrm{with}\; N=G\times_K\kk.
\end{equation}	
		
Let $\widehat{\mathcal{X}}$ be the total space of $TX\oplus N$ 
over $X$, and let 
$\widehat{\pi}:\widehat{\mathcal{X}}\rightarrow X$ be 
the natural projection. 
Let $Y=Y^{TX}+Y^N$, $Y^{TX}\in TX$, $Y^N\in N$ be the canonical 
sections of $\widehat{\pi}^*(TX\oplus N)$, $\widehat{\pi}^*(TX)$, 
$\widehat{\pi}^*(N)$ over $\widehat{\mathcal{X}}$.

Note that the natural action of $K$ on 
$C^\infty(\kg,\Lambda^\bullet(\kg^*)\otimes E)$ is given by 
\begin{equation}\label{eq3.5.12}
(k\cdot \phi)(Y)=\rho^{\Lambda^\bullet(\kg^*)\otimes E}(k)\,
\phi(\mathrm{Ad}(k^{-1})Y),\; \mathrm{for}\; 
\phi\in C^\infty(\kg,\Lambda^\bullet(\kg^*)\otimes E).
\end{equation}
Therefore %Thus we have 
\begin{align}\label{eq:3.53}
S^\bullet(T^*X\oplus N^*) \otimes 
\Lambda^\bullet(T^*X\oplus N^*)\otimes F
= G\times_{K} (S^\bullet(\kg^*)\otimes 
\Lambda^\bullet(\kg^*)\otimes E),
\end{align}
and the bundle 
$G\times_K C^\infty(\kg,\Lambda^\bullet(\kg^*)\otimes E)$
over $X$ is just 
$C^\infty(TX\oplus N,\widehat{\pi}^*
(\Lambda^\bullet(T^*X\oplus N^*)\otimes F))$. 

 By (\ref{eq3.5.12}), the $K$ action on  
$C^\infty(G\times \kg, 
\Lambda^\bullet(\kg^*)\otimes E)$ is given  by 
\begin{align}
\label{eq3.6.2}
(k\cdot s) (g,Y)= 
\rho^{\Lambda^\bullet(\kp^*\oplus \kk^*)\otimes E}(k)
s\left(gk,\Ad(k^{-1})Y\right).
\end{align}
If a vector space $W$ is a $K$-representation, we denote
by $W^K$ its K-invariant subspace. Then 
\begin{align}\label{eq:3.54}\begin{split}
C^\infty(G, &S^\bullet(\kg^*)\otimes 
\Lambda^\bullet(\kg^*)\otimes E)^K\\
&= C^\infty(X, S^\bullet(T^*X\oplus N^*) \otimes 
\Lambda^\bullet(T^*X\oplus N^*)\otimes F)\\
&``\simeq'' 
C^\infty(X, C^\infty(TX\oplus N,\widehat{\pi}^*
(\Lambda^\bullet(T^*X\oplus N^*)\otimes F)))\\
&= C^\infty(\widehat{\mathcal{X}}, \widehat{\pi}^*
(\Lambda^\bullet(T^*X\oplus N^*)\otimes F)).
\end{split}\end{align}

% Since  $K$ acts unitarily on 
% $\Lambda^\bullet(\kg^*)\otimes L_2(\kg)\otimes E
% =L_2(\kg,\Lambda^\bullet(\kg^*)\otimes E)$, we can regard 
% $G\times_K L_2(\kg,\Lambda^\bullet(\kg^*)\otimes E)$ 
% as an infinite dimensional vector bundle over $X$. 
% The smooth version of this bundle is given by 
% $G\times_K C^\infty(\kg,\Lambda^\bullet(\kg^*)\otimes E)$, 
% which is just 
% $C^\infty(TX\oplus N,\widehat{\pi}^*
% (\Lambda^\bullet(T^*X\oplus N^*)\otimes F))$. 
% Let $\mH$ be the 
% vector space of smooth sections over $X$ of the vector bundle 
% $C^\infty(TX\oplus N,\widehat{\pi}^*
% (\Lambda^\bullet(T^*X\oplus N^*)\otimes F))$.

As we saw in Section \ref{s2.1}, the connection form $\omega^\kk$
on $K$-principal bundle $p: G\rightarrow X=G/K$ also 
induces a connection on $C^\infty(TX\oplus N,\widehat{\pi}^*
(\Lambda^\bullet(T^*X\oplus N^*)\otimes F))$
over $X$, which is denoted by
$\nabla^{C^\infty(TX\oplus N,\widehat{\pi}^*
	(\Lambda^\bullet(T^*X\oplus N^*)\otimes F))}$.
In particular, for the canonical section $Y^{TX}$ of
$\widehat{\pi}^*(TX)$ over $\widehat{\mathcal{X}}$, 
the covariant differentiation with respect to the given 
canonical connection in the horizontal direction corresponding to 
$Y^{TX}$ is 
\begin{equation}\label{eq3.5.26}
\nabla_{Y^{TX}}^{C^\infty(TX\oplus N,\widehat{\pi}^*
(\Lambda^\bullet(T^*X\oplus N^*)\otimes F))}.
\end{equation}
	
Since the operator $\mathfrak{D}_b$ is $K$-invariant,
by (\ref{eq:3.54}), it 
descends to an operator $\mathfrak{D}^X_b$ acting on 
$C^\infty(\widehat{\mathcal{X}}, \widehat{\pi}^*
(\Lambda^\bullet(T^*X\oplus N^*)\otimes F))$. %$\mathcal{H}$.	
It is the same for the operator $\widehat{D}^\kg$, 
which descends to an 
operator $\widehat{D}^{\kg,X}$ over $X$.

Recall that $A$ is the self-adjoint $K$-invariant endomorphism
of $E$ in Section \ref{s2.1}.
For $b>0$, let $\mL^X_b, \mL^X_{A,b}$ act on 
%$C^\infty(TX\oplus N,\widehat{\pi}^*
%(\Lambda^\bullet(T^*X\oplus N^*)\otimes F))$
$C^\infty(\widehat{\mathcal{X}}, \widehat{\pi}^*
(\Lambda^\bullet(T^*X\oplus N^*)\otimes F))$ by %the formula
\begin{equation}\label{eq3.5.28}\begin{split}
&\mL_b^X=-\frac{1}{2}\widehat{D}^{\kg,X,2}
+\frac{1}{2}\mathfrak{D}_b^{X,2},\\
&\mL^X_{A,b}= \mL^X_{b} +A.
\end{split}\end{equation}
	
Let $\langle\cdot,\cdot\rangle$ be the 
usual $L_2$ Hermitian product 
on the vector space of smooth compactly supported sections of 
$\widehat{\pi}^*(\Lambda^\bullet(T^*X\oplus N^*)\otimes F)$ over 
$\widehat{\mathcal{X}}$.
Set
\begin{equation}
\label{eq3.5.27}
\begin{split}
\alpha=&\frac{1}{2}\Big(-\Delta^{TX\oplus N}
+\abs{Y}^2-m-n\Big)+N^{\Lambda^\bullet(T^*X\oplus N^*)},\\
\beta=&\nabla^{C^\infty(TX\oplus N,\widehat{\pi}^*
(\Lambda^\bullet(T^*X\oplus N^*)\otimes F))}_{Y^{TX}}
+ \widehat{c}(\mathrm{ad}(Y^{TX}))\\
&-c(\mathrm{ad}(Y^{TX})+i\theta \mathrm{ad}(Y^N))
-i\rho^E(Y^N),\\
\vartheta =& \frac{1}{2}\Big|[Y^N,Y^{TX}]\Big|^2.
\end{split}
\end{equation}

\begin{theo}\label{t3.11}
    \cite[Theorems 2.12.5, 2.13.2]{B11b}
We have
\begin{equation}\label{eq3.5.29}
\mL_b^X=\frac{\alpha}{b^2}+\frac{\beta}{b}+\vartheta.
\end{equation}
The operator $\frac{\partial}{\partial t}+\mL^X_b$ is 
hypoelliptic.
	
Also $\dfrac{1}{b}\nabla^{C^\infty(TX\oplus N,\widehat{\pi}^*
(\Lambda^\bullet(T^*X\oplus N^*)\otimes F))}_{Y^{TX}}$ 
is formally skew-adjoint with respect to 
$\langle\cdot,\cdot\rangle$ and 	
$\mL^X_b-\dfrac{1}{b}\nabla^{C^\infty(TX\oplus N,\widehat{\pi}^*
(\Lambda^\bullet(T^*X\oplus N^*)\otimes F))}_{Y^{TX}}$ 
is formally self-adjoint with respect to $\langle\cdot,\cdot\rangle$.
\end{theo}

\begin{rema}\label{t3.10a}
We will now explain the presence of the term 
$ic([Y^\kk,Y^\kp])$ in the right-hand side of (\ref{eq3.5.15}).
Instead of $\mathfrak{D}_b$, we could consider the operator
\begin{equation}\label{eq3.5.10}
D_b= \widehat{D}^\kg + \frac{1}{b}(\overline{d}^\kp
-i\overline{d}^\kk+\overline{d}^{\kp*}+i \overline{d}^{\kk *}).
\end{equation}
From \eqref{eq3.3.7}, \eqref{eq3.3.9}, \eqref{eq3.5.8} 
and \eqref{eq3.5.10}, we get
\begin{equation}\label{eq3.5.11}
D^2_b=\widehat{D}^{\kg,2}
+ \frac{1}{2b^2}(-\Delta^{\kp\oplus \kk}) 
+\frac{\sqrt{2}}{b} (\underline{Y}^\kp+ i\underline{Y}^\kk) 
+ \mathrm{zero\;order\; terms}.
\end{equation}

If $e\in\kk$, let $\nabla_{e,l}$ be the differentiation operator 
with respect to the left invariant vector field $e$, 
by \eqref{eq3.6.2}, for 
$s\in C^\infty(G, C^\infty(\kg,\Lambda^\bullet(\kg^*)\otimes E))^K$,
\begin{equation}\label{eq3.5.13}
\nabla_{e,l}s= (L^V_{[e,Y]}-\rho^E(e))s.
\end{equation}
Here $[e,Y]$ is a Killing vector field on $\kg=\kp\oplus \kk$ 
and the corresponding Lie derivative $L^V_{[e,Y]}$ 
acts on $C^\infty(\kg,\Lambda^\bullet(\kg^*))$. 
%We compute $L^V_{[e,Y]}$ easily
By \cite[(2.12.4)]{B11b}, we have the formula 
\begin{equation}\label{eq3.5.14}
L^V_{[e,Y]}=\nabla^{V}_{[e,Y]}-(c+\widehat{c})(\mathrm{ad}(e)).
\end{equation}
When we use the identification (\ref{eq:3.54}),
the operator $i\underline{Y}^\kk$ contributes the first order
differential operator $i\nabla^V_{[Y^N,Y^{TX}]}$ along %the 
%$\kp$. %and it is quadratic on $Y$. 
$TX$. This term is very difficult to control analytically. 
%Thus we need to modify the operator $D_b$ to eliminate this term. 

The miraculous fact is that after adding %we add the term
$ic([Y^\kk,Y^\kp])$ to $D_{b}$, 
%to eliminate the bad term $i\ul{Y}^{\kk}$, 
in the  operator $\mL^X_b$, we have eliminated
$i\nabla^V_{[Y^N,Y^{TX}]}$ and we add instead the term  
%we find  an extra fourth degree polynomial 
$\vartheta = \frac{1}{2}\abs{[Y^N,Y^{TX}]}^2$,
which is nonnegative.  This ensures that the operator 
$\frac{\alpha}{b^2}+\vartheta$ is %low bounded, 
bounded below. The operator $\mL^X_b$ is a nice operator.
%and functional analytic or probabilistic techniques can be applied. 
\end{rema}

\begin{prop}\label{t3.11a} \cite[Proposition 2.15.1]{B11b}
    We have the identity
\begin{align}\label{eq:3.5.31a}
\left[\mathfrak{D}^{X}_b, \mathcal{L}^X_{A,b}\right]=0.
\end{align}
\end{prop} 
\begin{proof} The classical Bianchi identity say that 
\begin{align}\label{eq:3.5.32a}
\left[\mathfrak{D}^{X}_b,\mathfrak{D}^{X,2}_{b}\right]=0.
\end{align}
By \eqref{eq3.5.9}, $\widehat{D}^{\kg,X,2}$
is the Casimir operator (up to a constant), so that
\begin{align}\label{eq3.5.31}
\left[\mathfrak{D}^{X}_b,\widehat{D}^{\kg,X,2}\right]=0.
\end{align}

We have the trivial $[\mathfrak{D}^{X}_b, A]=0$.
From \eqref{eq3.5.28}, \eqref{eq:3.5.32a} and \eqref{eq3.5.31},
we get \eqref{eq:3.5.31a}.
\end{proof}
% To establish Theorem \ref{t2.7}, we will view the heat kernel 
% $e^{-t\mathcal{L}^X_A}$ as an element $h$ of a compact group 
% in \eqref{eq1.10}, and $C^\infty(X,F)$ as the cohomology group
% of the fiberwise cohomology with fiber $TX\oplus N$. 
% Then Theorem \ref{t3.11a} allows us to apply 
% the argument in \eqref{eq1.12} to establish certain traces of 
% the operator $e^{-t\mathcal{L}_A^X}e^{-t\mathfrak{D}^{X,2}_b/2}
% =e^{-t\mathcal{L}^X_{A,b}}$ does not depend
% on the parameter $b>0$.

By analogy with \eqref{eq1.13}, we will need to show %ensure
that as $b\to 0$, in a certain sense,
\begin{align}
\label{eq3.5.32}
e^{-t \mathcal{L}^X_{A,b}}\to e^{-t \mathcal{L}^X_A}.
\end{align}
% The final step is to study the asymptotics as $b\to +\infty$
% of a certain trace of $e^{-t \mathcal{L}^X_{A,b}}$ associated with
% the semisimple element $\gamma\in G$, which should be localized 
% on the geodesics conjugate to $\gamma$, i.e., on the submanifold 
% $X(\gamma)$. 

We explain here an algebraic argument %with make 
which gives evidence for \eqref{eq3.5.32}. 
This will be  the analogue of \eqref{3.4}. %in the current situation. 
We denote by $H$ the fiberwise kernel of $\alpha$, 
so that 
\begin{align}
\label{eq3.5.33}
H=e^{-|Y|^2/2}\otimes F.
\end{align}
Let $H^\bot$ be the orthogonal to $H$ in 
$L_2(\widehat{\mathcal{X}},
\widehat{\pi}^*(\Lambda^\bullet(T^*X\oplus N^*)\otimes F))$. 

Note that $\beta$ maps $H$ to $H^\bot$. Let $\alpha^{-1}$
be the inverse of $\alpha$ restricted to $H^\bot$. 
Let $P, P^\bot$ be the orthogonal projections on 
$H$ and $H^\bot$ respectively. We embed $L_2(X,F)$ into 
$L_2(\widehat{\mathcal{X}},
\widehat{\pi}^*(\Lambda^\bullet(T^*X\oplus N^*)\otimes F))$
isometrically   via 
$s\to \widehat{\pi}^*s\,  e^{-|Y|^2/2}/\pi^{(m+n)/4}$. 

\begin{theo}\label{thm}\cite[Theorem 2.16.1]{B11b} 
    The following identify holds:
\begin{align}
\label{eq3.5.34}
P(\vartheta-\beta \alpha^{-1}\beta)P=\mL^X.
\end{align}
\end{theo}
\begin{proof}
 From \eqref{eq3.5.15}, we can write 
 \begin{align}\label{eq3.5.35} 
     \frac{1}{\sqrt{2}}\mathfrak{D}^X_b
=E_{1}+\frac{F_{1}}{b}, \quad \hbox{with }E_{1}
=\frac{1}{\sqrt{2}}\left(\widehat{D}^{\kg,X}
+ic([Y^{N}, Y^{TX}])\right).
 \end{align}
 Then comparing \eqref{eq3.5.28}, \eqref{eq3.5.29} and 
 \eqref{eq3.5.35}, we get 
 \begin{align}\label{eq3.5.36}
\alpha=F_{1}^2,\quad \beta=[E_{1},F_{1}],\quad 
\vartheta =E_{1}^2-\frac{1}{2}\widehat{D}^{\kg,X,2}. 
 \end{align}
 Since $H$ is the kernel of $F_{1}$, we have  $PF_{1}=F_{1}P=0$. 
We obtain thus
 \begin{align}\label{eq3.5.37}
 P(\vartheta-\beta \alpha^{-1}\beta)P
 =P\left(E_{1}^2-E_{1}P^\bot E_{1} 
 -\frac{1}{2}\widehat{D}^{\kg,X,2}\right)P
 =(PE_{1}P)^2-\frac{1}{2}P\widehat{D}^{\kg,X,2}P. 
 \end{align}
 But $H$ is of degree $0$ in $\Lambda^\bullet(\kg^*)$, 
 $\widehat{D}^{\kg}+ic([Y^{\kk}, Y^{\kp}])$ is of odd degree, 
 we know that $PE_{1}P=0$. Thus, \eqref{eq3.5.34} holds.
\end{proof}

%%%%%%%%%%%%%%%%%%%%%%%%%%%%%%%%%%
\subsection{Hypoelliptic orbital integrals}\label{s3.6}
Under the formalism of Section \ref{s2.2}, we replace now 
the finite dimensional vector space $E$ by the infinite dimensional 
vector space
\begin{align}
\label{eq:3.6.1}
\mathcal{E}=\Lambda^\bullet(\kp^*\oplus \kk^*)\otimes
S^\bullet(\kp^*\oplus \kk^*)\otimes E. 
\end{align}
Using \eqref{eq:3.54}, from now on, we 
will work systematically on $C^\infty(\widehat{\mathcal{X}}, \widehat{\pi}^*
(\Lambda^\bullet(T^*X\oplus N^*)\otimes F))$.

Let $dY$ be the volume element of $\kg=\kp\oplus \kk$ 
with respect to the scalar product 
$\langle\cdot ,\cdot \rangle=-B(\cdot,\theta \cdot)$.
It defines a fiberwise volume element  on the fiber 
$TX\oplus N$, which we still denote by $dY$. 
Our kernel $q(g)$ now acts as an endomorphism of 
$\mathcal{E}$ and verifies \eqref{eq2.13} and \eqref{eq2.15}.
In what follows, the operator $q(g)$ is given by continuous kernels
$q(g,Y,Y')$, $Y,Y'\in \kg$. 
Let $q((x,Y),(x',Y'))$, $(x,Y),(x',Y')\in \widehat{\mathcal{X}}$
be the corresponding kernel on $\widehat{\mathcal{X}}$. 

\begin{defi}\label{d3.6.1}\cite[Definition 4.3.3]{B11b}
For  a semisimple element $\gamma\in G$, we define 
$\tr_s^{[\gamma]}[Q]$ as in \eqref{eq2.3.15}, 
\begin{align}
\label{eq3.6.1}
\tr_s^{[\gamma]}[Q]=\int_{\kp^\bot(\gamma)\times \kg} 
\tr_s^{\Lambda^\bullet(\kp^*\oplus \kk^*)\otimes E}
\left[q(e^{-f}\gamma e^f,Y,Y)\right]r(f)dfdY
\end{align}
once it is well-defined. Note here 
$\tr_s^{\Lambda^\bullet(\kp^*\oplus \kk^*)\otimes E}[\cdot]
=\tr^{\Lambda^\bullet(\kp^*\oplus \kk^*)\otimes E}
[(-1)^{N^{\Lambda^\bullet(\kp^*\oplus \kk^*)}}\cdot]$,
i.e., we use the natural $\mathbb{Z}_2$-grading on 
$\Lambda^\bullet(\kp^*\oplus \kk^*)$.  
\end{defi}
% In \cite[\S 4.3]{B11b}, Bismut showed  that 
% the hypoelliptic orbital integral (\ref{eq3.6.1}) 
% is a real trace on certain nice algebra of kernel functions.

\begin{defi}
 \label{d3.6.2}Let ${\bf P}$ be the projection from 
 $\Lambda^\bullet(T^*X\oplus N^*)\otimes F$ on 
 $\Lambda^0(T^*X\oplus N^*)\otimes F$. 
\end{defi}

Recall that $e^{-t\mL^X_{A}}(x,x')$
is the heat kernel of $\mL^X_{A}$ in Section \ref{s2.1}.
For $t>0$, $(x,Y), (x', Y')\in \widehat{\mathcal{X}}$, put 
\begin{align}
\label{eq3.6.5}
q^X_{0,t}\big((x,Y), (x', Y')\big)={\bf P}
e^{-t\mL^X_{A}}(x,x')\pi^{-(m+n)/2}
e^{-\frac{1}{2}\left(|Y|^2+|Y'|^2\right)}{\bf P}.
\end{align}

Let $e^{-t\mathcal{L}_{A,b}^X}$ be the heat operator of 
$\mathcal{L}_{A,b}^X$ and 
$q^X_{b,t}\big((x,Y), (x', Y')\big)$ be the kernel of the 
heat operator $e^{-t\mathcal{L}_{A,b}^X}$ 
associated with the volume form $dx' dY'$.
In \cite[\S 11.5, 11.7]{B11b}, Bismut studied in detail %arguments on
the smoothness of $q^X_{b,t}((x,Y), (x', Y'))$
for $t>0, b>0$, $(x,Y), (x', Y')\in \widehat{\mathcal{X}}$.
In particular,  he showed that 
it is rapidly decreasing in the variables $Y,Y'$.

Now we state an important result 
\cite[Theorem 4.5.2]{B11b} whose proof was given in 
\cite[\S 14]{B11b} where Theorem \ref{thm} 
plays an important role. It ensures that
the hypoelliptic orbital integral is well-defined for 
$e^{-t\mathcal{L}_{A,b}^X}$ and that the analogue of 
Theorem \ref{t1.1} holds for $h= e^{-t\mathcal{L}^X_A}$ 
and $\mathfrak{D}_b^X$.

\begin{theo}
 \label{t3.6.3} Given $0<\epsilon \le M$, there exist $C,C'>0$ 
 such that for $0<b\le M, \epsilon \le t \le M$, 
 $(x,Y), (x', Y')\in \widehat{\mathcal{X}}$,
 \begin{align}
 \label{eq3.6.7}
 \Big|q^X_{b,t}\big((x,Y), (x', Y')\big)\Big|
 \le C\exp\Big(-C' \Big(d^2(x,x')+|Y|^2+|Y'|^2\Big)\Big).
 \end{align}
 Moreover, as $b\to 0$, 
 \begin{align}
 \label{eq3.6.8}
q^X_{b,t}\big((x,Y), (x', Y')\big)\to q^X_{0,t}\big((x,Y), (x', Y')\big).
 \end{align}
\end{theo}

The formal analogue %of the proof 
of Theorem \ref{t1.1} is as follows.
%\eqref{eq3.6.10}, we get:
\begin{theo}
 \label{t3.6.4}\cite[Theorem 4.6.1]{B11b}
 For any $b>0, t>0$, we have 
 \begin{align}\label{eq3.6.11}
 \tr^{[\gamma]}\left[e^{-t\mathcal{L}_{A}^X}\right]
 =\tr_s^{[\gamma]}\left[e^{-t\mathcal{L}_{A,b}^X}\right] .
 \end{align}
\end{theo}
\begin{proof} In \cite[\S 4.3]{B11b}, Bismut showed  that 
the hypoelliptic orbital integral (\ref{eq3.6.1}) 
is a trace on certain algebras of operators given by smooth kernels
which exhibit a Gaussian decay like in (\ref{eq3.6.7}).
By Theorem \ref{t3.6.3}, 
    the kernel function $q^X_{b,t}$ is in this algebra.
 As in (\ref{eq1.12}),  by Proposition \ref{t3.11a},
    \begin{align}\label{eq3.6.12}\begin{split}
\frac{\partial}{\partial b} 
\tr_s^{[\gamma]}\left[e^{-t\mathcal{L}_{A,b}^X}\right] 
&= \tr_s^{[\gamma]}\left[-t \Big(\frac{\partial}{\partial b} 
\mathcal{L}_{A,b}^X\Big) e^{-t\mathcal{L}_{A,b}^X}\right] \\
&= -t \tr_s^{[\gamma]}\left[ \frac{1}{2}\Big[\mathfrak{D}^{X}_b,
\frac{\partial}{\partial b} \mathfrak{D}^{X}_b\Big]
e^{-t\mathcal{L}_{A,b}^X}\right] \\
&=  -\frac{t}{2}\tr_s^{[\gamma]}\left[ \Big[\mathfrak{D}^{X}_b,
(\frac{\partial}{\partial b} \mathfrak{D}^{X}_b)
e^{-t\mathcal{L}_{A,b}^X}\Big]\right] =0.
\end{split}\end{align}
By Theorem \ref{t3.6.3}, we have
\begin{align}\label{eq3.6.10}
\lim_{b\to 0} \tr_s^{[\gamma]}
\left[ e^{-t\mathcal{L}_{A,b}^X}\right]= 
\tr^{[\gamma]}\left[e^{-t\mathcal{L}_A^X}\right]. 
\end{align}
  From  (\ref{eq3.6.12}), (\ref{eq3.6.10}), we get
  (\ref{eq3.6.11}).
\end{proof}

%%%%%%%%%%%%%%%%%%%%%%%%%%%%%%%%%
\subsection{Proof of Theorem \ref{t2.7}}\label{s3.7}
For $b>0$, $s(x,Y)\in C^\infty(\widehat{\mathcal{X}},
\widehat{\pi}^*(\Lambda^\bullet(T^*X\oplus N^*)\otimes F))$, 
set 
\begin{align}
\label{eq3.8.1}
F_bs(x,Y)=s(x,-bY).
\end{align}
Put 
\begin{align}\label{eq3.8.2a}
\ul{\mathcal{L}}^X_{A,b}=F_b\mathcal{L}_{A,b}^XF_b^{-1}.
\end{align}
Let $\ul{q}^X_{b,t}((x,Y),(x',Y'))$ be the kernel associated with 
$e^{-t\ul{\mathcal{L}}_{A,b}^X}$. When $t=1$, 
we will write $\ul{q}^X_{b}$ instead of $\ul{q}^X_{b,1}$. 
Then from \eqref{eq3.8.2a}, we have
\begin{align}
\label{eq3.8.2}
\ul{q}^X_{b,t}\big((x,Y),(x',Y')\big)
=(-b)^{m+n} q^X_{b,t}\big((x,-bY),(x',-bY')\big).
\end{align}
Let $a^{TX}$ be the vector field on $X$ associated with 
$a$ in \eqref{eq2.3.4} induced by the left action of $G$ 
on $X$ (cf. \eqref{eq2.8}). Let $d(\cdot,X(\gamma))$ 
be the distance function to $X(\gamma)$.

\begin{theo}
 \label{t3.8.1} \cite[Theorem 9.1.1, (9.1.6)]{B11b}
 Given $0<\epsilon\le M$, there exist $C,C'>0$, such that 
 for any $b\ge 1$, $\epsilon \le t\le M$, 
 $(x,Y),(x',Y')\in \widehat{\mathcal{X}}$, 
 \begin{align}\label{eq3.8.3}
 \left|\ul{q}^X_{b,t}\big((x,Y),(x',Y')\big)\right|\le Cb^{4m+2n}
 \exp\left(-C\Big(d^2(x,x')+|Y|^2+|Y'|^2\Big)\right).
 \end{align}
 Given $\delta>1, \beta>0$, $0<\epsilon\le M$, there exist
 $C,C'>0$, such that 
 for any $b\ge 1$, $\epsilon \le t\le M$, 
 $(x,Y)\in \widehat{\mathcal{X}}$, if $d(x,X(\gamma))\ge \beta$,
 \begin{align}\label{eq3.8.4}
 \left|\ul{q}^X_{b,t}\big((x,Y),\gamma(x,Y)\big)\right|
 \le Cb^{-\delta}\exp\left(-C'\Big(d^2_\gamma(x)+|Y|^2\Big)\right).
 \end{align}
 Given $\delta> 1, \beta>0, \mu>0$, there exist $C,C'>0$ 
 such that for any $b\ge 1$, $(x,Y)\in \widehat{\mathcal{X}}$,
 if $d(x,X(\gamma))\le \beta$, and $|Y^{TX}-a^{TX}(x)|\ge \mu$, 
 \begin{align}\label{eq3.8.5}
 \left|\ul{q}^X_{b,t}\big((x,Y),\gamma(x,Y)\big)\right|
 \le Cb^{-\delta}e^{-C'|Y|^2}.
 \end{align}
\end{theo}

In view of Theorem \ref{t3.8.1}, the proof of Theorem \ref{t2.7} 
consists in obtaining %getting
the asymptotics of 
$\tr_s^{[\gamma]}[e^{-t\mathcal{L}^X_{A,b}}]$ 
as $b\to +\infty$. By \cite[(2.14.4)]{B11b}, 
the operator in \eqref{eq3.5.28} associated with $B/t$ is 
up to conjugation, $t\mathcal{L}^X_{\sqrt{t}b}$. Observe that 
$J_{\gamma}(Y^{\kk}_0)$ is unchanged when replacing
the bilinear form $B$ by $B/t$, $t>0$. Thus we only need to 
establish the corresponding result for $t=1$. 

When $f\in \kp^\bot(\gamma)$, 
we identify $e^f$ with $e^fp1$. For $f\in \kp^\bot(\gamma)$, 
$Y\in (TX\oplus N)_{e^f}$, set 
\begin{align}\label{eq3.8.6}
\ul{Q}^X_{b}(e^f,Y)=
\tr_s^{\Lambda^\bullet(T^*X\oplus N^*)\otimes F}
\left[\gamma \ul{q}_b^X\big((e^f,Y),\gamma(e^f,Y)\big)\right].
\end{align}
Then 
\begin{align}\label{eq3.8.7}
\tr_s^{[\gamma]}\left[e^{-\mathcal{L}^X_{A,b}}\right]
=\int_{(e^f,Y)\in \widehat{\pi}^{-1}
\kp^\bot(\gamma)}\ul{Q}^X_{b}(e^f,Y)r(f)dfdY.
%\\=\int_{|f|\ge \beta}+\int_{|f|<\beta,|Y^{TX}-
%a^{TX}(e^f)|\ge \mu }+\int_{|f|<\beta,|Y-a^{TX}(e^f)|< \mu }.
\end{align}
Take $\beta\in]0,1]$. By Theorem \ref{t3.8.1}, as $b\to +\infty$,
\begin{align}
\label{eq3.8.8}
\begin{aligned}
&\int_{(e^f,Y)\in \widehat{\pi}^{-1}\kp^\bot(\gamma), |f|
\ge \beta}\ul{Q}^X_{b}(e^f,Y)r(f)dfdY\to 0,\\
&\int_{(e^f,Y)\in \widehat{\pi}^{-1}\kp^\bot(\gamma), 
|f|< \beta,|Y^{TX}-a^{TX}(e^f)|\ge \mu}\ul{Q}^X_{b}(e^f,Y)
r(f)dfdY\to 0.
\end{aligned}
\end{align}
We need to understand the integral on the domain 
$|f|<\beta$, $|Y^{TX}-a^{TX}(e^f)|<\mu$, when $b\to +\infty$. 

Let $\pi: \mathcal{X}\to X$ be the total space of 
the tangent bundle $TX$ to $X$. Let $\varphi_{t}|_{t\in \R}$ be 
the group of diffeomorphisms of $\mathcal{X}$ induced by 
the geodesic flow. By \cite[Proposition 3.5.1]{B11b}, 
$\varphi_{1}(x, Y^{TX})= \gamma \cdot(x, Y^{TX})$ 
is equivalent to   %means exactly
$x\in X(\gamma)$ and $Y^{TX}= a^{TX}(x)$. Equation
\eqref{eq3.8.8} shows that as $b\to +\infty$, 
%we are localizing the contribution of 
the right-hand side of \eqref{eq3.8.7} localizes near
 the minimizing geodesic $x_{t}$ connecting $x$ and $\gamma x$
so that $\dot{x} = a^{TX}$.

Let $N(\gamma)$ be the vector bundle on $X(\gamma)$ which is the
analogue of the vector bundle %as the analogue of 
$N$ on $X$ in \eqref{eq3.5.25}. 
Then $N(\gamma)\subset N|_{X(\gamma)}$. 
Let $N^\bot(\gamma)$ be the orthogonal to $N(\gamma)$ in 
$N|_{X(\gamma)}$. Clearly,
\begin{align}
N^\bot(\gamma)=Z(\gamma)\times_{K(\gamma)}\kk^\bot(\gamma).
\end{align}
Let $p_\gamma: X\to X(\gamma)$ be the projection 
defined by \eqref{eq2.3.11} and \eqref{2.3.14}.
We trivialize the vector bundles $TX,N$ by parallel transport
along the geodesics orthogonal to $X(\gamma)$ 
with respect to the connection $\nabla^{TX}, \nabla^{N}$, 
so that $TX, N$ can be identified with 
$p^*_\gamma TX|_{X(\gamma)}, p^*_\gamma N|_{X(\gamma)}$. 
At $x=p1$, we have
\begin{align}
\label{eq3.8.11}
N(\gamma)=\kk(\gamma), \qquad N^\bot(\gamma)
=\kk^\bot(\gamma).
\end{align}
Therefore at $\rho_\gamma(1,f)$, we may write $Y^N\in N$
in the form 
\begin{align}
\label{eq3.8.12}
Y^N=Y_0^\kk+Y^{N,\bot}, \qquad  \hbox{with }
Y_0^\kk\in \kk(\gamma),Y^{N,\bot} \in \kk^\bot(\gamma).
\end{align}
Let $dY_0^\kk$, $dY^{N,\bot}$ be the volume elements on 
$\kk(\gamma)$, $\kk^\bot(\gamma)$, so that 
\begin{align}
\label{eq3.8.13}
dY^N=dY_0^\kk\, dY^{N,\bot}.
\end{align}
To evaluate the limit of \eqref{eq3.8.7} as $b\to +\infty$
for $\beta>0$, we may by \eqref{eq3.8.8}, as well 
 consider the integral
\begin{multline}
\label{eq3.8.14}
\int_{|f|< \beta,|Y^{TX}-a^{TX}(e^f)|< \mu} 
\ul{Q}^X_{b}\left(e^f,Y\right)r(f)dfdY^{TX}dY_0^\kk dY^{N,\bot}\\
=b^{-4m-2n+2r}\int_{|f|< b^2\beta,|Y^{TX}|< b^2\mu}
 \ul{Q}^X_{b}\left(e^{f/b^2},\frac{Y^{TX}}{b^2}+a^{TX}(e^{f/b^2}),
 Y_0^\kk+\frac{Y^{N,\bot}}{b^2}\right)\\
 r(f/b^2)dfdY^{TX}dY_0^\kk dY^{N,\bot}.
\end{multline}
	
Let $\ul{\kz}(\gamma)$ be the another copy of $\kz(\gamma)$, 
and let $\ul{\kz}(\gamma)^*$ be the corresponding copy 
of the dual of $\ul{\kz}(\gamma)$. Also, for $u\in \kz(\gamma)^*$, 
we denote by $\ul{u}$ the corresponding element in 
$\ul{\kz}(\gamma)^*$. Let $e_1,\cdots,e_r$ be a basis 
of $\kz(\gamma)$, let $e^1,\cdots,e^r$ be 
the corresponding dual basis of $\kz(\gamma)^*$. 

Put $\mG=\End(\Lambda^\bullet(\kg^*))\widehat{\otimes} 
\Lambda^\bullet(\ul{\kz}(\gamma)^*)$. 
Let $e_{r+1},\cdots,e_{m+n}$ be a basis of 
$\kz^{\bot}(\gamma)$, and let  $e^*_{r+1},\cdots,e^*_{m+n}$
be the dual basis to 
$e_{r+1},\cdots,e_{m+n}$  with respect to 
$B|_{\kz^{\bot}(\gamma)}$. 
Then $\mG$ is generated by all the monomials in
$c(e_i),\widehat{c}(e_i), 1\le i\le m+n$, $\ul{e}^j$, $1\le j\le r$. 
Let $\widehat{\tr}_s$ be the linear map from $\mG$ into $\R$
that, up to permutation, vanishes on all monomials except 
those of the following form:
\begin{align}
\label{eq3.8.16}
\widehat{\tr}_s\left[c(e_1)\ul{e}^1\cdots c(e_r)\ul{e}^r 
c(e^*_{r+1})\widehat{c}(e_{r+1})\cdots c(e^*_{m+n})
\widehat{c}(e_{m+n})\right]=(-1)^r(-2)^{m+n-r}.
\end{align}
For $u\in \mG$, $v\in \End(E)$, we define 
\begin{align}
\label{eq3.8.17}
\widehat{\tr}_s[uv]=\widehat{\tr}_s[u]\tr^E[v].
\end{align}
Set 
\begin{align}\label{eq3.8.15}
{\bf \alpha}=\sum_{i=1}^r c(e_i)\ul{e}^i\in c(\kz(\gamma))
\widehat{\otimes} \Lambda^\bullet(\ul{\kz}(\gamma)^*).
\end{align}
\begin{defi}\label{d3.8.3}  
 Let $\kq_{b}^X((x,Y),(x',Y'))$ denote the smooth kernel 
 associated with $e^{-\ul{\mathcal{L}}_{A,b}^X-{\bf \alpha}}$, and 
 \begin{align}
 \label{eq3.8.18}
 Q_{b}^X(x,Y)=\gamma \kq_{b}^X\big((x,Y),\gamma(x,Y)\big).
 \end{align}
\end{defi}
Since $\ul{\mL}_{A,b}^X+{\bf \alpha}$ can be obtained from
$\ul{\mL}_{A,b}^X$ by a conjugation, by 
a simple argument on Clifford algebras, we get:
\begin{prop}
 \label{t3.8.4} \cite[Proposition 9.5.4]{B11b}. For $b>0$, 
 the following identity holds:
 \begin{align}\label{eq3.8.19}
\ul{Q}^X_{b}(x,Y)
 =b^{-2r}\widehat{\tr}_s\left[Q_{b}^X(x,Y)\right]. 
 \end{align}
\end{prop}

Now we define a limit operator acting on 
$C^\infty(\kp\times \kg,\Lambda^\bullet(\kg^*)
\widehat{\otimes} \Lambda^\bullet(\ul{\kz}(\gamma)^*)
\otimes E)$. We denote by $y$ the tautological 
section of the first component of $\kp$ in $\kp\times \kg$, 
and by $Y=Y^\kp+Y^\kk$ the tautological section 
of $\kg=\kp\oplus \kk$. Let $dy$ the volume form on $\kp$ 
and let $dY$ the volume form on $\kg=\kp\oplus \kk$.
	
\begin{defi}
 Given $Y_0^{\kk}\in \kk(\gamma)$, set 
 \begin{multline}
 \label{eq3.8.20}
 \mP_{a,Y^\kk_0}=\frac{1}{2}\Big|[Y^\kk,a]
 +[Y_0^\kk,Y^\kp]\Big|^2-\frac{1}{2}\Delta^{\kp\oplus\kk}
 + \sum_{i=1}^r c(e_i)\ul{e}^{i}      %{\bf \alpha}
 -\nabla^H_{Y^\kp}-\nabla^V_{[a+Y_0^\kk,[a,y]]}\\
 -\widehat{c}(\ad(a))+c(\ad(a)+i\theta \ad(Y_0^\kk))
 \end{multline}
 acting on $C^\infty(\kp\times \kg,\Lambda^\bullet(\kg^*)
 \widehat{\otimes}\Lambda^\bullet(\ul{\kz}(\gamma)^*)
 \otimes E)$.
\end{defi}

Let $R_{Y_0^\kk}((y,Y),(y',Y^{\prime}))$ 
be the smooth kernel of $e^{-\mP_{a,Y^\kk_0}}$ 
with respect to the volume form $dydY$ on $\kp\times \kg$. 
Then 
\begin{align}
\label{eq3.8.21}
R_{Y_0^\kk}((y,Y),(y',Y^{\prime}))\in 
\End(\Lambda^\bullet(\kz^{\bot}(\gamma)^*))
\widehat{\otimes}c(\kz(\gamma))\widehat{\otimes} 
\Lambda^\bullet(\ul{\kz}(\gamma)^*).
\end{align}

The following result gives an estimate and pointwise asymptotics 
of $Q_{b}^X$.
\begin{theo}
 \label{t3.8.7}\cite[Theorems 9.5.6, 9.6.1]{B11b} 
 Given $\beta>0$, there exist $C,C_\gamma'>0$ 
 such that for $b\ge1$, $f\in \kp^\bot(\gamma)$, 
 $|f|\le \beta b^2$, and $|Y^{TX}|\le \beta b^2$,
\begin{multline}\label{eq3.8.21a}
b^{-4m-2n}\left|Q^X_{b}\left(e^{f/b^2},a^{TX}(e^{f/b^2})
+Y^{TX}/b^2,Y_0^\kk+Y^{N,\bot}/b^2\right)\right|\\
\le C\exp\left(-C'\left|Y_0^\kk\right|^2
-C_\gamma'\left(|f|^2+\left|Y^{TX}\right|^2
+\left|(\Ad(k^{-1})-1)Y^{N,\bot}\right|
+\left|[a,Y^{N,\bot}]\right|\right)\right).
\end{multline}
As $b\to +\infty$,
\begin{multline}
\label{eq3.8.22}
b^{-4m-2n}Q^X_{b}\left(e^{f/b^2},a^{TX}(e^{f/b^2})
+Y^{TX}/b^2,Y_0^\kk+Y^{N,\bot}/b^2\right)\\
\to e^{-(|a|^2+|Y_0^\kk|^2)/{2}}
\Ad\left(k^{-1}\right)R_{Y_0^\kk}\left((f,Y),
\Ad(k^{-1})(f,Y)\right)\rho^E\left(k^{-1}\right)
e^{-i\rho^E(Y_0^\kk)-A}.
\end{multline}
\end{theo}
A crucial computation in \cite[Theorem 5.5.1, (5.1.11)]{B11b}
gives the following key result.
%detailed computation %calculation:
\begin{theo}
\label{t3.8.8}
For $Y_0^\kk\in \kk(\gamma)$, we have the identity
\begin{align}\label{eq3.8.23}
(2\pi)^{r/2}\int_{\kp^{\bot}(\gamma)\times 
(\kp\oplus \kk^\bot(\gamma))}\widehat{\tr}_s
\left[\Ad(k^{-1}) R_{Y_0^\kk}\left((y,Y),\Ad(k^{-1})(y,Y)\right)
\right]dydY= J_\gamma\left(Y_0^\kk\right).
\end{align}
\end{theo}
From Theorems \ref{t3.8.7}, \ref{t3.8.8}, 
\eqref{eq3.8.6}-\eqref{eq3.8.8}, 
and \eqref{eq3.8.19}, 
we obtain Theorem \ref{t2.7}.

\begin{exem}\cite[\S 10.6]{B11b}, \cite[\S 5.1]{B12}.
 In  Example \ref{e2.9}, we have %$X=\R$, 
$N=0$, 
$\widehat{\mathcal{X}}=TX\oplus N=T\R=\R\oplus \R.$ 
Using the coordinates $(x,y)\in \R\oplus \R$, we get
\begin{align}
\label{eq3.9.1}
&\mL^X_b=M_b+\frac{N^{\Lambda^\bullet(\R)}}{b^2},
& \hbox{with } M_b=\frac{1}{2b^2}
\left(-\frac{\partial^2}{\partial y^2}+y^2-1\right)
+\frac{y}{b}\frac{\partial}{\partial x}.
\end{align}
% The heat kernel $p_t(x,x')$ associated with 
% $e^{t\Delta^{\R}/2}$, where 
% $\Delta^{\R}=\frac{\partial^2}{\partial x^2}$, is given by 
% \begin{align}
% \label{eq3.9.2}
% p_t(x,x')=\frac{1}{\sqrt{2\pi t}}e^{-\frac{(x-x')^2}{2t}}.
% \end{align}
The heat kernel $p_{b,t}((x,y),(x',y'))$ associated with 
$e^{-tM_b}$ depends only on $x'-x,y,y'$. %Note that 
The heat kernel of the  operator 
$-\frac{\partial^2}{\partial y^2}+y\frac{\partial}{\partial x}$
was first calculated     %evaluated
by Kolmogorov \cite{K34},
and $p_{b,t}((x,y),(x',y'))$ has been computed explicitly in 
\cite[Proposition 10.5.1]{B11b}.

Let $a\in \kp=\R$. Then $a$ acts as %by
translation by $a$ on the first component % copy 
of $\R\oplus \R$.  
From \eqref{eq3.6.1} and \eqref{eq3.9.1}, we deduce that
\begin{align}
\label{eq3.9.3}
\begin{aligned}
&\tr_s^{[a]}\left[e^{-t\mL_{b}^X}\right]
=\left(1-e^{-t/b^2}\right)\int_{\R}p_{b,t}\big((0,Y),(a,Y)\big)dY.
\end{aligned}
\end{align}
Theorem \ref{t3.6.4} can be stated in the special case
of this example as follows.
\begin{theo}
 \label{t3.9.1} For any $t>0, b>0$, we have
 \begin{align}
 \label{eq3.9.4}
\tr^{[a]}\left[e^{t\Delta^{\R}/2}\right]
=\tr_s^{[a]}\left[e^{-t\mL_{b}^X}\right].
 \end{align}
\end{theo}
\begin{proof}
We give a simple direct    %formal
proof which can be ultimately easily justified. 
Note that
\begin{align}
\label{eq3.9.5}
M_b=\frac{1}{2b^2}\left(-\frac{\partial^2}{\partial y^2}
+\left(y+b\frac{\partial}{\partial x}\right)^2-1\right)
-\frac{1}{2}\frac{\partial^2}{\partial x^2}.
\end{align}
By \eqref{eq3.9.5}, we get 
\begin{align}
\label{eq3.9.6}
e^{-b\frac{\partial^2}{\partial x\partial y}}M_b
e^{b\frac{\partial^2}{\partial x\partial y}}
=\frac{1}{2b^2}\left(-\frac{\partial^2}{\partial y^2}
+y^2-1\right)-\frac{1}{2}\frac{\partial^2}{\partial x^2}.
\end{align}
Using the fact that $p_{b,t}((x,y),(x',y'))$ only depends on 
$x'-x$, $y$, $y'$, we deduce from \eqref{eq3.9.6} that 
\begin{align}
\label{eq3.9.7}\begin{split}
\int_{\R}p_{b,t}\big((0,y),(a,y)\big)dy
&=\tr\left[e^{-\frac{t}{2b}
\left(-\frac{\partial^2}{\partial y^2}+y^2-1\right)}\right]
\tr^{[a]}\left[e^{t\Delta^{\R}/2}\right]\\
&= \frac{1}{1-e^{-t/b^2}} \tr^{[a]}\left[e^{t\Delta^{\R}/2}\right],
\end{split}\end{align}
since the spectrum of the harmonic oscillator 
$\frac{1}{2}(-\frac{\partial^2}{\partial y^2}+y^2-1)$ is $\N$.
By (\ref{eq3.9.3}) and  \eqref{eq3.9.7}, we get \eqref{eq3.9.4}. 
\end{proof}

By (\ref{eq3.9.4}), we can compute the limit 
as $b\to +\infty$  of the right-hand side of (\ref{eq3.9.3})
from the explicit formula of $p_{b,t}((x,y),(x',y'))$,
and in this way we get (\ref{eq2.4.11}).
In other words, we interpret %really  
(\ref{eq2.4.11}) as a consequence of a local index theorem.
% in view the heat kernel as certain index, 
% then we evaluate it via the hypoelliptic version 
% of the local index technique.
\end{exem}

%%%%%%%%%%%%%%%%%%%%%%%%%%%%%%%
\subsection{A brief idea on the proof of Theorems
\ref{t3.6.3}, \ref{t3.8.1}, \ref{t3.8.7}}\label{s3.8}

The wave operator for %of
the elliptic Laplacian has the property
of finite propagation speed, %and this plays an important role on 
which explain the Gaussian decay %property
of the elliptic heat kernel.
%But there is no wave operator associated with 
The hypoelliptic Laplacian does not have a wave equation.
%as it is only a first order differential operator along 
%the horizontal direction.

One difficult point in Theorems \ref{t3.6.3}, \ref{t3.8.1}
and \ref{t3.8.7} is to get the uniform Gaussian-like  estimate.  
%  the distance function $d(x,x')$
% uniformly on the parameter $b$ for 
% hypoelliptic heat kernels, in particular comparing 
% with Bismut-Lebeau's book \cite{BL08}, where the base manifold
% is compact.
% As the hypoelliptic Laplacian is second order along the fiber and 
% first order along the horizontal direction,
% roughly we can understand that it has some natural as the 
% heat kernel along the fiber and as he wave operator 
% along the horizontal direction, this explains at least formally
% why we can get the Gaussian-like decay estimate properties 
% on $d(x,x')$. 
% Here is a more intuitive version: 
Let us give an argument back to \cite[\S 12.3]{B11b}
which explains some heuristic relations
of the hypoelliptic heat equation to the wave equation on $X$.
%assume temporarily  that 
Here $q^X_{b,t}$ will denote scalar hypoelliptic heat kernel
on the total space $\mathcal{X}$ of the tangent bundle $TX$. Put
\begin{align}\label{eq:3.9.10}\begin{split}
&\sigma_{b,t}((x,Y), x') = \int_{Y'\in  T_{x'}X}
q^X_{b,t}\big((x,Y),(x',Y')\big)  dY',\\
&M_{b,t}((x,Y), x') = 
\frac{1}{\sigma_{b,t}((x,Y), x') }
\int_{Y'\in T_{x'}X}
q^X_{b,t}\big((x,Y),(x',Y')\big) ( Y'\otimes Y') dY'.
\end{split}\end{align}
Then $M_{b,t}((x,Y), x')$ takes its values in symmetric positive 
endomorphisms of $T_{x'}X$. We can associate to $M_{b,t}$
the second order elliptic operator acting on $C^\infty (X,\R)$,
\begin{align}\label{eq:3.9.11}
{\bf M}_{b,t}(x,Y) g(x')= \left\langle  
\nabla^{TX}_{\cdot}\nabla_{\cdot},
M_{b,t}((x,Y), x')  g(x')\right\rangle,
\end{align}
where the operator $\nabla^{TX}_{\cdot}\nabla_{\cdot}$
acts on the variable $x'$. Then we have \cite[(12.3.12)]{B11b}
% \begin{align}\label{eq:3.9.14}
% P^X_{b,t}(x,x')=\int_{Y\in \widehat{\mathcal{X}}_{x}}
% \int_{Y'\in \widehat{\mathcal{X}}_{x'}} 
% q^X_{b,t}\big((x,Y),(x',Y')\big) dY \, dY',
% \end{align}
%then we can verify that
\begin{align}\label{eq:3.9.15}
\left(b^2  \frac{\partial^2}{\partial t^2}
+ \frac{\partial}{\partial t} -{\bf M}_{b,t}(x,Y) \right)
\sigma_{b,t}((x,Y), \cdot)=0.
\end{align}
This is a hyperbolic equation. 
% But when only the variable $x$ is considered, there is a 
% wave equation quality to the heat kernel operator 
% for the hypoelliptic Laplacian.
As $b\to 0$, it converges in the proper sense to the
standard parabolic
heat operator
\begin{align}\label{eq:3.9.16}
\Big(\frac{\partial}{\partial t} - \frac{1}{2}\Delta\Big) 
p_{t}(x,\cdot)=0.
\end{align}

The above consideration plays an important role
in the proof given in \cite{B11b} of the estimates 
(\ref{eq3.6.7}), (\ref{eq3.8.3}), (\ref{eq3.8.4}) 
and (\ref{eq3.8.21a}).
% From (\ref{eq:3.9.15}), it is  if $b\to +\infty$, 
% $P^X_{b,t}$ is a wave operator and if $b=0$,
% it becomes a heat operator.
%\newpage
%%%%%%%%%%%%%%%%%%%%%%%%%%%%%%%
\section{Analytic torsion and dynamical zeta function}\label{s4}
%%%%%%%%%%%%%%%%%%%%%%%%%%%%

Recall that a flat vector bundle $(F,\nabla)$ with flat connection 
$\nabla$ over a smooth manifold $M$
%is equivalent that there exists 
comes from a representation 
$\rho: \pi_{1}(M)\to {\rm GL}({\bf q},\C)$ so that 
if  $\widetilde{M}$ is the universal cover of $M$,  then 
$F= \widetilde{M}\times_{\rho} \C^{\bf q}$.
%the induced vector bundle by $\rho$, here
The analytic torsion associated with 
a flat vector bundle on a smooth compact 
Riemannian manifold $M$ is a classical spectral invariant 
defined by Ray and Singer \cite{RS71} in 1971. 
It is a regularized determinant of the Hodge Laplacian for 
the de Rham complex associated with this flat vector bundle.

For $\Gamma\subset G$ a discrete cocompact torsion free 
subgroup of a connected reductive Lie group $G$, 
if $Z= \Gamma\backslash G/K$ is the locally symmetric space 
as in (\ref{eq2.9}), then $\Gamma = \pi_{1}(Z)$. 
By the superrigidity theorem of Margulis 
 \cite[Chap. VII, §5]{Margulis91}, 
if the real rank of $G$ is $\geq 2$, a general 
representation of $\Gamma$ is not too far from 
a unitary representation of $\Gamma$ or the restriction to
$\Gamma$ of a representation of $G_{\C}$,
the complexification of $G$. 
See \cite[Chap. XIII, 4.6]{BorelWallach00} for more details.

Assume that the difference of the complex ranks of 
$G$ and $K$ %the rank of its maximal compact group 
is different from $1$. For a flat vector bundle
induced by a $G_{\C}$-representation, 
as an application of Theorem \ref{t0.1}, 
we obtain a vanishing result of individual 
orbital integrals that appear %appeared
in the supertrace of the heat kernel %to define 
from which the analytic torsion can be obtained.
In particular, this implies that the associated analytic torsion
is equal to $1$ (cf. Theorem \ref{t4.2.4}). 
% For a unitary representation, by using the above vanishing result 
% of orbital integrals for the trivial representation, 
% this implies again the analytic torsion is $1$
% (cf. Theorem \ref{t4.3.1}).

We explain finally Shen's recent solution on Fried's conjecture 
for locally symmetric spaces: 
for any unitary representation of $\Gamma$
such that the cohomology of the associated flat vector bundle
on $Z=\Gamma\backslash G/K$ vanishes, % it identifies
 the value at zero of a Ruelle dynamical zeta function
%associated with the unitary representation of $\Gamma$,
identifies to the associated analytic torsion.

This section is organized as follows. In Section \ref{s3.1}, 
we introduce %a classical spectral invariant:
the Ray-Singer analytic torsion. 
%for a flat vector bundle on a smooth compact manifold. 
In Section \ref{s4.2}, 
we study the analytic torsion on locally symmetric spaces
for flat vector bundles  induced by 
a representation of  $G_{\C}$.
%the complexification of the total group,
%as an application of Theorem \ref{t0.1}.
% we deduce that 
% the analytic torsion is $1$ if the difference of the complex rank of 
% the total group and its maximal compact group is different to $1$. 
Finally in Section \ref{s4.3},
%we explain first  the corresponding 
%result in Section \ref{s4.2} for a unitary flat vector bundle.
 we describe Shen's solution of 
Fried's conjecture in the case of locally symmetric spaces. 
In Section \ref{s4.4}, we make some remarks on
 related research directions.
 
 %%%%%%%%%%%%%%%%%%%%%%%%%%%%%%%%
\subsection{Analytic torsion}
\label{s3.1}

Let $M$ be a compact manifold of dimension $m$.
Let $(F,\nabla)$ be a flat complex vector
bundle on $M$  with flat connection $\nabla$ 
(i.e., its curvature $(\nabla)^2=0$). 
%We denote by $\mathrm{rk} (F)$ the rank of $F$.
The flat connection $\nabla$ induces 
an exterior differential operator $d$ on $\Omega^\bullet(M,F)$,
the vector space of differential forms on $M$ with values in $F$,
and $d^2=0$. Let $H^\bullet (M,F)$ be the cohomology group of 
the complex $(\Omega^\bullet(M,F),d)$ as in \eqref{eq3.1.1}.

Let $h^F$ be a Hermitian metric on $F$. Then as explained in 
Section \ref{s1.2}, 
$g^{TM}$ and $h^F$ induce naturally a Hermitian product on 
$\Omega^\bullet(M,F)$.
%and we still have the corresponding formal adjoint $d^*$. Put
Let $D$ be as in \eqref{eq3.1.3}. 

We introduce here a refined spectral invariant of $D^2$ which is 
particularly interesting. 

Let $P$ be the orthogonal projection from $\Omega^\bullet(M,F)$ 
onto $\Ker (D)$ and let $P^\perp=1-P$. 
Let $N$ be the number operator acting on $\Omega^\bullet(M,F)$,
i.e., multiplication by $j$ on $\Omega^j(M,F)$. For $s\in \C$ 
and $\Re(s)>\dfrac{m}{2}$, set
\begin{equation}
\label{eq3.1.8}
\begin{split}
\theta(s)&= - \sum^m_{j=0}(-1)^jj\tr|_{\Omega^j(M,F)}
[(D^2)^{-s}P^\perp]\\
&=-\frac{1}{\Gamma(s)}\int_0^{+\infty }
\tr_s[Ne^{-tD^2}P^\perp]t^s\frac{dt}{t},
\end{split}
\end{equation}
where $\Gamma(\cdot)$ is the Gamma function. 

From the small time heat kernel expansion
(cf. \cite[Theorem 2.30]{BeGeVe04}), we know that 
$\theta(s)$ is well-defined for $\Re(s)>\dfrac{m}{2}$ and 
extends holomorphically near $s=0$.
\begin{defi}\label{d3.1} \cite{RS71}
The (Ray-Singer) analytic torsion is defined as 
\begin{align}\label{eq:3.1.9a}
T(g^{TM},h^{F})=\exp\Big(\dfrac{1}{2}
\dfrac{\partial \theta}{\partial s}(0)\Big).
\end{align}
\end{defi} 

We have the formal identity, 
\begin{equation}\label{eq3.1.9}
T(g^{TM},h^{F})
=\prod^m_{j=0} \det (D^2|_{\Omega^j(M,F)})^{(-1)^jj/2}.
\end{equation}
\begin{rema}
 \label{r3.1a}
 \begin{enumerate}
\item[a)] If $h^F$ is parallel with respect to $\nabla$, then 
$F$ is induced by a unitary representation of $\pi_1(M)$,
and we say that $(F,\nabla,h^F)$ is a unitary flat vector bundle.
 In this case, if $m$ is even and $M$ is orientable,
 by a Poincar\'{e} duality argument, we have $T(g^{TM},h^F)=1$.
\item[b)] If $m$ is odd, and $H^\bullet(M,F)=0$,
then $T(g^{TM},h^{F})$  does not depend on the choice 
of $g^{TM},\; h^F$, thus it is a 
topological invariant (cf. \cite[Theorem 0.1]{BZ92}).
\end{enumerate}
\end{rema}
%%%%%%%%%%%%%%%%%%%%%%%%%%%%%%%%%%%
\subsection{Analytic torsion for locally symmetric spaces}\label{s4.2}
We use the same notation and assumptions as in Section \ref{s2}. 
Recall that $\rho^E: K\to \mathrm{U}(E)$ is a finite dimensional 
unitary representation of $K$, and $F=G\times_K E$ is the induced
Hermitian vector bundle on the symmetric space $X=G/K$. 
Assume form now on
that the complexification $G_\C$ of $G$ exists,  
and the representation $\rho^E$ is induced by a holomorphic 
representation of $G_\C\to  \mathrm{Aut}(E)$, that is still 
denoted by $\rho^E$.  We have the canonical identification of 
$G\times_K E$ as a trivial bundle $E$ on $X$:
\begin{align}
\label{eq4.2.1}
F= G\times_K E\to X\times E, \quad (g,v)\to \rho^E(g)v. 
\end{align}
This induces a canonical flat connection $\nabla$ on $F$ such that
\begin{align}
\label{eq4.2.2}
\nabla=\nabla^F+\rho^E\omega^{\kp}.
\end{align}

\begin{rema}\label{t4.2.1}
Let $U$ be a maximal compact subgroup of $G_{\C}$. Then $U$ is
the compact form of $G$ and $\mathfrak{u}=i\kp\oplus \kk$
is its Lie algebra. 
By Weyl's unitary trick \cite[Proposition 5.7]{Knapp01}, if 
$U$ is simply connected, it is equivalent to consider 
representations of $G$, of $U$ on $E$, or holomorphic 
representations of the complexification $G_{\C}$ of $G$ on $E$, 
or representations of $\kg$, or $\mathfrak{u}$ on $E$.
\end{rema}

We fix a $U$-invariant Hermitian metric on $E$.
This implies in particular it is  $K$-invariant, 
and $\rho^E(v)\in \End(E)$ is symmetric for $v\in \kp$. 
This induces a Hermitian metric $h^F$ on $F$.  
As in Section \ref{s3.1}, we consider now the operator 
$D$ acting on $\Omega^\bullet(X,F)$ induced by $g^{TX},h^F$.

Let $C^{\kg,X}$ be the Casimir operator of $G$ acting on 
$C^\infty(X,\Lambda^\bullet(T^*X)\otimes F)$ as in \eqref{eq2.5}. 
Then by \cite[(2.6.11)]{B11b} and \cite[Proposition 8.4]{BMZ16}, 
we have 
\begin{align}
\label{eq4.2.3}
D^2=C^{\kg,X}-C^{\kg,E}.
\end{align}

Let $T$ be  a maximal torus in $K$ and let $\kt\subset \kk$ 
be its Lie algebra. Set 
\begin{align}
\label{eq4.2.4}
\kb=\{v\in \kp: [v,\kt]=0\}.
\end{align}
Put 
\begin{align}
\label{eq4.2.5}
\mathfrak{h}=\kb\oplus \kt.
\end{align}
By \cite[p. 129]{Knapp01}, we know that $\mathfrak{h}$ is 
a Cartan subalgebra of $\kg$ and that $\dim \kt$ is 
the complex rank of $K$ and $\dim \mathfrak{h}$ is 
the complex rank of $G$. Also, $m$ and $\dim \kb$ 
have the same parity.

For $\gamma= e ^{a} k^{-1} \in G$ a semisimple element
as in (\ref{eq2.3.4}), let $K^{0}(\gamma)\subset K(\gamma)$ be 
the connected component of the identity. 
Let $T(\gamma)\subset K^{0}(\gamma)$ be a maximal torus in 
$K^{0}(\gamma)$, and let $\kt(\gamma)\subset \kk(\gamma)$
be its Lie algebra. By (\ref{eq2.3.5}) and (\ref{eq2.3.9}), 
$k$ commutes with $T(\gamma)$,
thus by \cite[Theorem 4.21]{Knapp01}, there exists $k_{1}\in K$
    such that $k_{1}T(\gamma)k_{1}^{-1}\subset T$, 
    $k_{1}k k_{1}^{-1}\subset T$. By working on 
$k_{1}\gamma k_{1}^{-1}=e ^{\Ad(k_{1}) a} (\Ad(k_{1}) k)^{-1}$
instead of $\gamma$,
we may and we will assume that $T(\gamma)\subset T$, $k\in T$.
In particular $\kt(\gamma)\subset \kt$. Set
\begin{align}\label{eq4.2.4a}
\kb(\gamma)=\{v\in \kp: [v,\kt(\gamma)]=0, \, \Ad(k) v=v\}.
\end{align}
Then 
\begin{align}\label{eq4.2.5a}
\kb\subset \kb(\gamma)\quad \text{ and  } \,  \,
\kb(1)= \kb.
\end{align}
Recall that $N^{\Lambda^\bullet(T^*X)}$ is the number operator
on $\Lambda^\bullet(T^*X)$.
\begin{theo}\label{t4.2.2}
\cite[Theorem 7.9.1]{B11b}, \cite{BMZ11},
\cite[Theorem 8.6, Remark 8.7]{BMZ16}. 
For any semisimple element $\gamma\in G$, 
if $m$ is even, or if $m$ is odd and $\dim \kb(\gamma)\ge 2$, 
then for any $t>0$, we have 
\begin{align}
\label{eq4.2.6}
\tr_s^{[\gamma]}\left[\left(N^{\Lambda^\bullet(T^*X)}
-\frac{m}{2}\right)e^{-\frac{t}{2}D^2}\right]=0.
\end{align}
\end{theo}
\begin{proof}
By Theorem \ref{t2.7}, \eqref{eq2.5} and \eqref{eq4.2.3}, 
for any $t>0$ and any semisimple element $\gamma\in G$,
\begin{multline}\label{eq4.2.7}
\tr_s^{[\gamma]}\left[\left(N^{\Lambda^\bullet(T^*X)}
-\frac{m}{2}\right)e^{-\frac{t}{2}D^2}\right]
=\frac{e^{-|a|^2/2t}}{(2\pi t)^{p/2}}
e^{\frac{t}{16}\tr^\kp[C^{\kk,\kp}]
+\frac{t}{48}\tr^{\kk}[C^{\kk,\kk}]}\\
\int_{\kk(\gamma)} J_\gamma(Y^\kk_0) 
\tr_s^{\Lambda^\bullet(\kp^*)\otimes E}
\left[\left(N^{\Lambda^\bullet(\kp^*)}-\frac{m}{2}\right)
\rho^{\Lambda^\bullet(\kp^*)\otimes E}(k^{-1})
e^{-i \rho^{\Lambda^\bullet(\kp^*)\otimes E}(Y^\kk_0)
+\frac{t}{2}C^{\kg,E}}\right]\\
e^{-|Y^\kk_0|^2/2t} \frac{dY^\kk_0}{(2\pi t)^{q/2}}.
	\end{multline}
But 
\begin{multline}
\label{eq4.2.8}
\tr_s^{\Lambda^\bullet(\kp^*)\otimes E}
\left[\left(N^{\Lambda^\bullet(\kp^*)}-\frac{m}{2}\right)
\rho^{\Lambda^\bullet(\kp^*)\otimes E}(k^{-1})
e^{-i \rho^{\Lambda^\bullet(\kp^*)\otimes E}(Y^\kk_0)
+\frac{t}{2}C^{\kg,E}}\right]\\
=\tr_s^{\Lambda^\bullet(\kp^*)}
\left[\left(N^{\Lambda^\bullet(\kp^*)}-\frac{m}{2}\right)
\rho^{\Lambda^\bullet(\kp^*)}(k^{-1})
e^{-i \rho^{\Lambda^\bullet(\kp^*)}(Y^\kk_0)}\right]
\tr^{E}\left[\rho^E(k^{-1})e^{-i \rho^E(Y^\kk_0)
+\frac{t}{2}C^{\kg,E}}\right].
\end{multline}
	
If $u$ is an isometry of $\kp$, we have 
\begin{align}
\label{eq4.2.9}
\begin{aligned}
&\tr_s^{\Lambda^\bullet(\kp^*)}[u]=\det\left(1-u^{-1}\right),\\
&\tr_s^{\Lambda^\bullet(\kp^*)}
\left[N^{\Lambda^\bullet(\kp^*)}u\right]
=\frac{\partial }{\partial s} \det\left(1-u^{-1}e^s\right)(0).
\end{aligned}
\end{align}
If the eigenspace associated with the eigenvalue $1$
is of dimension $\ge1$, the fist quantity in \eqref{eq4.2.9} 
vanishes. If it is of dimension $\ge 2$, the second expression 
in \eqref{eq4.2.9} also vanishes. Also, if $m$ is even and 
$u$ preserves the orientation, then 
\begin{align}\label{eq4.2.10}
\tr_s^{\Lambda^\bullet(\kp^*)}
\left[\left(N^{\Lambda^\bullet(\kp^*)}-\frac{m}{2}\right)u\right]=0.
\end{align}

From \eqref{eq4.2.7}, \eqref{eq4.2.9} and \eqref{eq4.2.10}, 
we get \eqref{eq4.2.6}.
\end{proof}

Now let $\Gamma$ be a discrete torsion free
cocompact subgroup of $G$. 
Set $Z=\Gamma\backslash X$. Then
$\pi_{1}(Z)=\Gamma$ and the flat vector bundle
$F$ descents as a flat vector bundle $F$ over $Z$.

\begin{theo}\label{t4.2.4}\cite[Remark 8.7]{BMZ16}.
For a flat vector bundle $F$ on $Z=\Gamma\backslash X$ 
induced by a holomorphic representation of $G_{\C}$,
    if $m$ is even, or if $m$ is odd and $\dim \mathfrak{b}\ge3$, 
    then $T(g^{TZ},h^F)=1$.
\end{theo}
\begin{proof} Under the condition of Theorem \ref{t4.2.4},
 from Theorem \ref{t4.2.2},  \eqref{eq2.4.10}
 and \eqref{eq4.2.5a}, we get 
\begin{align}\label{eq4.2.11}
\tr_s\left[\left(N^{\Lambda^\bullet(T^*Z)}-\frac{m}{2}\right)
e^{-\frac{t}{2}D^{Z,2}}\right]=0.
\end{align}
Now Theorem \ref{t4.2.4} is a direct consequence of 
\eqref{eq1.11} for $h=1$, \eqref{eq3.1.8}
and \eqref{eq4.2.11}.
\end{proof}

\begin{rema}
 \label{r4.2.5}
 \begin{itemize}
\item[a)] If $F$ is trivial, i.e., it is induced by the trivial 
representation of $G$, then Theorem \ref{t4.2.2}
under the condition of Theorem  \ref{t4.2.4} was first 
obtained by Moscovici-Stanton \cite[Theorem 2.1]{MS91}.

\item[b)] Assume $G$ is semisimple, then the induced metric
$h^F$ on $F$ is unimodular, i.e., the metric $h^{\det F}$ 
on $\det F:= \Lambda^\mathrm{max}F$  induced by $h^F$
is parallel with respect to the flat connection on 
$\det F$. %on $X$.
In this case, 
Theorem \ref{t4.2.2} for $\gamma=1$ was first obtained
by Bergeron-Venkatesh \cite[Theorem 5.2]{BergeronV13}, 
and M\"{u}ller and Pfaff  %\cite{MullerPfaff13}, 
\cite{MullerPfaff13a}
gave a new proof of Theorem \ref{t4.2.4}. 
%\item[c)] 
%There are two classes of finite-dimensional 
%representations of $\Gamma$ for which we can work with
%self-adjoint operators and apply the usual Selberg trace formula. 
%These are unitary representations and restrictions of 
%rational representations of $G$. In general, not every 
%representation of $\Gamma$ belongs to one of these classes.
% If $\rank(G)\geq  2$, the superrigidity theorem of 
% Margulis  \cite[Chap. VII, §5]{Margulis91} implies that a general 
% representation of $\Gamma$ is not too far from a representation 
% which is either unitary or the restriction of 
% a $G_{\C}$ representation.  
% See \cite[Chap. XIII, 4.6]{BorelWallach00} for more details.

\item[c)] We can drop the condition on torsion freeness of 
$\Gamma$ in \eqref{eq4.2.11}. 
\end{itemize}
\end{rema}

\begin{rema}\label{r4.2.6}
For $p,q\in \N$, let $\mathrm{SO}^0(p,q)$ 
be the connected component of the identity in the real group 
$\mathrm{SO}(p,q)$. By \cite[Table V p. 5.18]{Helgason01} and
\cite[Table C1 p. 713, Table C2 p. 714]{Knapp01}, 
among the noncompact simple connected complex groups 
such that $m$ is odd and $\dim \kb=1$, there is only 
$\mathrm{SL}_2(\C)$, and among the noncompact simple real 
connected groups, there are only $\mathrm{SL}_3(\R)$, 
$\mathrm{SL}_4(\R)$, $\mathrm{SL}_2(\mathbb{H})$, 
and $\mathrm{SO}^0(p,q)$ with $pq$ odd $>1$. 
Also, by \cite[p. 519, 520]{Helgason01}, 
$\mathfrak{sl}_2(\C)=\mathfrak{so}(3,1)$, 
$\mathfrak{sl}_4(\R)=\mathfrak{so}(3,3)$, 
$\mathfrak{sl}_2(\mathbb{H})=\mathfrak{so}(5,1)$. 
Therefore the above  list can be reduced to 
$\mathrm{SL}_3(\R)$ and $\mathrm{SO}^0(p,q)$
with $pq$ odd $>1$.
\end{rema}

Assume from now on that
$\rho:\Gamma\to U({\bf q})$ be a unitary representation. 
Then $F=X\times_{\Gamma} \C^{\bf q}$ is a flat vector bundle on
$Z=\Gamma\backslash X$ with metric $h^F$ induced by 
the canonical metric on $\C^{\bf q}$,
%It is parallel with respect to the flat connection $\nabla^{F,f}$. 
i.e.,  $F$ is a unitary flat vector bundle on $Z$ 
with holonomy $\rho$. 
By Remark \ref{r3.1a}, if $m$ is even, then 
$T(g^{TZ},h^{F})=1$.
Thus we can simply assume that $m$ is odd.

Observe that the pull back
of $(F,h^F)$ over $X$ is $\C^{\bf q}$ with canonical metric, thus 
the heat kernel on $X$ is given by
\begin{align}
\label{eq4.3.1}
e^{-tD^2}(x,x')=e^{-tD_0^2}(x,x')\otimes \mathrm{Id}_{\C^{\bf q}}
\end{align}
where $e^{-tD_0^2}(x,x')\in \Lambda^\bullet(T_x^*X)
\otimes \Lambda^\bullet(T_{x'}^*X)^*$ is the heat kernel on $X$ 
for the trivial representation 
$G\to \mathrm{Aut}(\C)$. %as in Section \ref{s3.1}. 
Thus,
\begin{multline}
\label{eq4.3.2}
\tr_s\left[\left(N^{\Lambda^\bullet(T^*X)}-\frac{m}{2}\right)
\gamma e^{-tD^{2}}(\gamma^{-1}\widetilde{z},
\widetilde{z})\right]\\
=\tr[\rho(\gamma)]\tr_s\left[\left(N^{\Lambda^\bullet(T^*X)}
-\frac{m}{2}\right)\gamma e^{-tD^{2}_0}
(\gamma^{-1}\widetilde{z},\widetilde{z})\right].
\end{multline}
Note that for $\gamma\in \Gamma$, $\tr[\rho(\gamma)]$
depends only on the conjugacy class of $\gamma$, 
thus form \eqref{eq2.10}, \eqref{eq2.21} and \eqref{eq4.3.2}, 
we get the analogue of Theorem \ref{t2.5a},
\begin{multline}
\label{eq4.3.3}
\tr_s\left[\left(N^{\Lambda^\bullet(T^*Z)}-\frac{m}{2}\right)
e^{-tD^{Z,2}}\right]\\
=\sum_{[\gamma]\in [\Gamma]}\vol \big(\Gamma\cap Z(\gamma)
\backslash X(\gamma)\big)
\tr[\rho(\gamma)]\tr_s^{[\gamma]}
\left[\left(N^{\Lambda^\bullet(T^*X)}-\frac{m}{2}\right)
e^{-tD^{2}_0}\right].
\end{multline}

Since the metric $h^F$ is given by the unitary representation $\rho$ 
and $g^{TZ}$ is induced by the bilinear form $B$ on $\kg$, 
we denote the analytic torsion in Section \ref{s3.1} by $T(F)$. 
%By Remark \ref{r3.3}, if $m$ is even, then $T(F)=1$;
%if $F$ is acyclic, i.e., $H^\bullet(Z,F)=0$, then 
%$T(F)$ is a topological invariant. 

By Theorem \ref{t4.2.2} for the trivial representation 
$G\to \mathrm{Aut}(\C)$ and \eqref{eq4.3.3}, we get %as in 
a result similar to Theorem \ref{t4.2.4}.
\begin{theo}\label{t4.3.1}  \cite[Corollary 2.2]{MS91}.
 %If $m$ is even, or I
For a unitary flat vector bundle $F$ on 
$Z=\Gamma \backslash X$,
%induced by a unitary representation $\rho: \Gamma \to U({\bf q})$,
if $m$ is odd and $\dim \kb\ge 3$, 
 then for $t>0$, we have 
\begin{align}
\label{eq4.3.4}
\tr_s\left[\left(N^{\Lambda^\bullet(T^*Z)}-\frac{m}{2}\right)
e^{-t D^{Z,2}/2}\right]=0.
\end{align}
In particular,  %the analytique torsion 
$T(F)=1$.
\end{theo}

%%%%%%%%%%%%%%%%%%%%%%%%%%%%%%%%
\subsection{Fried's conjecture for locally symmetric 
spaces}\label{s4.3}

The possible relation of the topological torsion to 
the dynamical systems was first observed by Milnor
 \cite{Milnor68} in 1968. 
A quantitative description of their relation was formulated by 
Fried \cite{Fried86} when $Z$ is a closed oriented hyperbolic 
manifold. Namely, he showed that for an acyclic unitary flat vector
bundle $F$, the value at zero of the Ruelle dynamical zeta function, 
constructed using the closed geodesics in $Z$ 
and the holonomy of $F$, is equal to the associated analytic torsion.  
In \cite[p.66]{Fried87}, Fried's conjectured that a similar result 
still holds for general closed locally homogenous manifolds.
In 1991, Moscovici-Stanton \cite{MS91} made an important 
contribution to the proof of
Fried's conjecture for locally symmetric spaces. 

Let $\Gamma\subset G$ be a discrete cocompact torsion free 
subgroup of a connected reductive Lie group $G$.
Then we get the symmetric space $X= G/K$ 
and the locally symmetric space  $Z= \Gamma\backslash X$.
By Remark \ref{r3.1a}, we may assume 
that $\dim Z=m$ is odd.

Recall that $[\Gamma]$ is the set of conjugacy classes of $\Gamma$.
For $[\gamma]\in [\Gamma]\backslash \{1\}$, denote by 
$B_{[\gamma]}$ 
the space of closed geodesics in $[\gamma]$.  
As a subset of the loop space $LZ$, we equipped 
$B_{[\gamma]}$ the induced topology and smooth structure. 
By Proposition \ref{p2.9},   %\cite{DKV79}, 
$B_{[\gamma]}\simeq \Gamma\cap Z(\gamma)\backslash 
X(\gamma)$ is a compact locally symmetric space, 
and the elements of $B_{[\gamma]}$ have the same length 
$l_{[\gamma]}>0$.  The group $\mathbb{S}^1$ acts on 
$B_{[\gamma]}$ by rotation. This action is locally free. 
Denote by $\chi_{\mathrm{orb}}(\mathbb{S}^1\backslash 
B_{[\gamma]})\in \mathbb{Q}$ the orbifold Euler 
characteristic number
for the quotient orbifold $\mathbb{S}^1\backslash B_{[\gamma]}$. 
Recall that if $e(\mathbb{S}^1\backslash B_{[\gamma]},
\nabla^{T(\mathbb{S}^1\backslash B_{[\gamma]})})
\in \Omega^\bullet(\mathbb{S}^1\backslash B_{[\gamma]},
o(T(\mathbb{S}^1\backslash B_{[\gamma]})))$ 
is the Euler form %class
defined using Chern-Weil theory for the Levi-Civita connection 
$\nabla^{T(\mathbb{S}^1\backslash B_{[\gamma]})}$, then 
\begin{align}
\label{eq4.3.20}
\chi_{\mathrm{orb}}(\mathbb{S}^1\backslash 
B_{[\gamma]})=\int_{\mathbb{S}^1\backslash 
B_{[\gamma]}}e(\mathbb{S}^1\backslash B_{[\gamma]},
\nabla^{T(\mathbb{S}^1\backslash B_{[\gamma]})}).
\end{align}
Let %Denote by 
\begin{align}\label{eq4.3.5}
n_{[\gamma]}=\left|\Ker\big(\mathbb{S}^1\to 
\mathrm{Diff}(B_{[\gamma]})\big)\right|
\end{align}
be the generic multiplicity of $B_{[\gamma]}$.

\begin{defi}\label{d4.3.2}
    Given a representation $\rho:\Gamma\to U({\bf q})$,
 we say that the dynamical zeta function $R_\rho(\sigma)$
 is well-defined 
 if the following properties hold:
 \begin{enumerate}
 \item For $\sigma\in \C$, $\mathrm{Re}(\sigma)\gg1$, the sum 
 \begin{align}
 \label{eq4.3.6}
 \Xi_\rho(\sigma)=\sum_{[\gamma]\in [\Gamma]\backslash \{1\}}
 \tr[\rho(\gamma)]\frac{\chi_{\mathrm{orb}}
 (\mathbb{S}^1\backslash B_{[\gamma]})}{n_{[\gamma]}}
 e^{-\sigma l_{[\gamma]}}.
 \end{align}
 %converges 
 defines to a holomorphic function. 
 \item The function $R_\rho(\sigma)=\exp(\Xi_\rho(\sigma))$ 
 has a meromorphic extension to $\sigma\in \C$.
 \end{enumerate}
\end{defi}

Note that $\gamma\in \Gamma$ is primitive means 
if $\gamma=\beta^k$, 
$\beta\in \Gamma$, $k\in \N^*$, then $\gamma=\beta$. 
\begin{rema}\label{r4.3.3}
If $Z$ is a compact oriented hyperbolic manifold, then 
$\mathbb{S}^1\backslash B_{[\gamma]}$ is a point. 
Moreover, if $\rho$ is the trivial representation, then 
\begin{align}\label{eq4.3.7}
R_\rho(\sigma)=\exp\left(\sum_{[\gamma]\in 
[\Gamma]\backslash \{1\}}
\frac{1}{n_{[\gamma]}}e^{-\sigma l_{[\gamma]}}\right)
=\prod_{[\gamma]\ \mathrm{primitive}, \gamma\neq1}
\left(1-e^{-\sigma l_{[\gamma]}}\right)^{-1}.
\end{align}
\end{rema}

\begin{theo}\label{t4.3.3}\cite{Shen16}
  %Assume $\dim \kb=1$. 
For any unitary flat vector bundle $F$ on $Z$ with holonomy $\rho$,
 the dynamical zeta function $R_\rho(\sigma)$ is
  a well-defined meromorphic function on $\C$
  which is holomorphic for $\mathrm{Re}(\sigma)\gg1$.
  Moreover, there exist explicit constants $C_\rho\in \R$
  and $r_\rho\in \Z $  such that, when  $\sigma\to 0$,
  we have
\begin{align}\label{eq4.3.10}
  R_\rho(\sigma)=C_\rho T(F)^2\sigma^{r_\rho}
  +\mathcal{O}(\sigma^{r_\rho+1}).
\end{align}
 If $H^\bullet (Z,F)=0$,  then
\begin{align}\label{eq4.3.11}
  C_\rho=1,\qquad r_\rho=0,
\end{align}
so that 
\begin{align}\label{eq4.3.9}
R_\rho(0)=T(F)^2.
\end{align}
\end{theo}

\begin{proof}[Proof of Theorem \ref{t4.3.3}]
The most difficult part of the proof is to 
express the $R_{\rho}(\sigma)$ as a product of determinant
of shifted Casimir operators,
\footnote{The gap in Moscovici-Stanton's paper comes from using
an operator $\Delta^{j,l}_\phi$ \cite[p. 206]{MS91}
which could not be defined on $Z$.}
in fact being the analytic torsion.
Shen's idea is to interpret the right-hand side of \eqref{eq4.3.6} 
as the Selberg trace formula (by eliminating the term $\tr^{[1]}$) 
of the heat kernel for some representations of $K$ by using 
Theorems \ref{t2.5a} and \ref{t2.7}.  

By Theorem \ref{t4.3.1}, we can concentrate on the proof in 
the case $\dim \kb=1$.  From now on, we assume $\dim \kb=1$.

 Up to conjugation, there exists 
a unique standard parabolic subgroup $Q\subset G$ with 
Langlands decomposition $Q=M_QA_QN_Q$ such that 
$\dim  A_Q=1$. Let $\mathfrak{m}, \kb, \kn$ be 
the respective Lie algebras of $M_Q, A_Q, N_Q$.
Let $M$ be the connected component of identity of $M_Q$.
Then $M$ is a connected reductive group with  
maximal compact subgroup $K_M=M\cap K$ and 
with Cartan decomposition 
$\mathfrak{m}=\kp_\mathfrak{m}\oplus \kk_\mathfrak{m}$, 
and $K_M$ acts on $\kp_\mathfrak{m}$, $M$ acts on $\kn$ 
and $M$ acts trivially on $\kb$.  
We have an identity of real $K_M$-representations
\begin{align}\label{eq:introdp}
  \kp\simeq\kp_\mathfrak{m}\oplus\kb \oplus \kn,
\end{align}
and $\dim \kn$ is even. Moreover there exists $\nu\in \kb^*$
such that (cf.  \cite[Proposition 6.2]{Shen16}) 
\begin{align}\label{eq4.3.13a}
[a,f]=\langle \nu,a\rangle f
\quad \text{ for any  }  f\in \kn,  a\in \kb.
\end{align}
For $\gamma\in G$ semisimple, Shen  \cite[Proposition 4.11]{Shen16} 
observes that
\begin{align}\label{eq4.3.14a}
\gamma \text{  can be conjugated
into }  H:=\exp(\kb)T \, \, (cf. \eqref{eq4.2.5}) 
\text{  iff }  \dim \kb(\gamma)=1.
\end{align}

Let $R(K,\R), R(K_M,\R)$ be the real representation rings  of 
$K$ and $K_M$. We can prove that the restriction 
$R(K,\R)\to R(K_M,\R)$ is injective.
The key result \cite[Theorem 6.11]{Shen16} is that 
the $K_M$-representation on $\kn$ has a unique lift in $R(K,\R)$. 
As $\kp$ is a $K$-representation, %for any $i,j$, 
%$\Lambda^i(\kp_\mathfrak{m}^*)$, $\Lambda^j(\kn^*)$
$K_{M}$-action on $\kp_\mathfrak{m}$
also lifts to $K$  %lift in $R(K,\R)$ 
by lifting $\kb$ as a trivial $K$-representation in \eqref{eq:introdp}. 

For a real finite dimensional representation $\varsigma$ of 
$M$ on the vector space $E_\varsigma$, such that 
$\varsigma|_{K_M}$ %on $K_M$ 
can be lifted into $R(K,\R)$, this implies that there exists 
a real finite dimensional $\Z_2$-representation 
$\widehat{\varsigma}=\widehat{\varsigma}^+
- \widehat{\varsigma}^-$ of $K$ on 
$E_{\widehat{\varsigma}}
=E_{\widehat{\varsigma}}^+ -  E_{\widehat{\varsigma}}^-$
such that we have the equality in $R(K_M,\R)$,
\begin{align}\label{eq4.3.13}
%\left(E_{\widehat{\varsigma}}^+- E_{\widehat{\varsigma}}^-\right)
E_{\widehat{\varsigma}}|_{K_M}
=\sum_{i=0}^{\dim \kp_\mathfrak{m}}
(-1)^i\Lambda^i(\kp^*_\mathfrak{m})\otimes (E_\varsigma|_{K_M}).
\end{align}
Let $\mathcal{E}_{\widehat{\varsigma}}
=G\times_K E_{\widehat{\varsigma}}$ 
be the induced $\Z_2$-graded vector bundle on $X$, and 
$\mathcal{F}_{\widehat{\varsigma}}
=\Gamma\backslash \mathcal{E}_{\widehat{\varsigma}}$. 
Let $C^{\kg,Z,\widehat{\varsigma},\rho}$ be 
the Casimir element of $G$ acting on 
$C^\infty(Z,\mathcal{F}_{\widehat{\varsigma}}\otimes F)$. 
Modulo some technical conditions, 
%by highly nontrivial computations using   %from
%Theorem \ref{t2.5a}, 
with the help of Theorem \ref{t4.2.2} and (\ref{eq4.3.14a}),
Shen  \cite[Theorems 5.3, 7.3]{Shen16} obtains the identity  %got 
\begin{multline}
\label{eq4.3.14}
\tr_s\left[e^{-tC^{\kg,Z,\widehat{\varsigma},\rho}/2}\right]
={\bf q}\vol(Z)\tr^{[1]}\left[e^{-tC^{\kg,X,\widehat{\varsigma}}/2}
\right]\\
+\frac{1}{\sqrt{2\pi t}}e^{-c_{\varsigma}t}
\sum_{[\gamma]\in [\Gamma]\backslash \{1\}}
\tr[\rho(\gamma)]\frac{\chi_{\rm{orb}}
(\mathbb{S}^1\backslash B_{[\gamma]})}{n_{[\gamma]}}
\frac{\tr^{E_{\varsigma}}[\varsigma(k^{-1})]}
{\left|\det(1-\Ad(\gamma))|_{\kz_0^\bot}\right|^{1/2}}
|a|e^{-|a|^2/2t},
\end{multline}
where $c_{\varsigma}$ is some explicit constant,
and $k^{-1}$ is defined in \eqref{eq2.3.4}. 
We do not write here the exact formula for 
$\tr^{[1]}[e^{-tC^{\kg,X,\widehat{\varsigma}}/2}]$.
%as it is less relevant. 

By \eqref{eq4.3.14}, if we set
\begin{align}\label{eq4.3.15}
\Xi_{\varsigma,\rho}(\sigma)
=\sum_{[\gamma]\in [\Gamma]\backslash \{1\}}
\tr[\rho(\gamma)]\frac{\chi_{\rm{orb}}(\mathbb{S}^1
\backslash B_{[\gamma]})}{n_{[\gamma]}}
\frac{\tr^{E_{\varsigma}}[\varsigma(k^{-1})]}
{\left|\det(1-\Ad(\gamma))|_{\kz_0^\bot}\right|^{1/2}}
e^{-|a|\sigma},
\end{align}
we need to eliminate the denominator 
$\left|\det(1-\Ad(\gamma))|_{\kz_0^\bot}\right|^{1/2}$ 
to relate $\Xi_{\varsigma,\rho}$ to $\Xi_{\rho}$,
the logarithm of  the dynamical zeta function $R_{\rho}(\sigma)$. 

Set $2l =\dim \kn$. The following observation is crucial.

\begin{prop}\label{t4.3.6} \cite[Proposition 6.5]{Shen16} For 
$\gamma= e^ak^{-1}\in H:=\exp(\kb)T$, %(cf. \eqref{eq4.2.5}), 
$a\neq0$, with $\nu\in \kb^*$ in (\ref{eq4.3.13a}),  we have
\begin{align}
\label{eq4.3.16}
\left|\det(1-\Ad(\gamma))|_{\kz_0^\bot}\right|^{1/2}
=\sum_{j=0}^{2l}(-1)^j\tr^{\Lambda^j(\kn^*)}
\left[\Ad\left(k^{-1}\right)\right]e^{(l-j)|\nu||a|},
\end{align}
% here $\nu\in \kb^*$ such that $[a,f]=\langle \nu,a\rangle f$
% for any $f\in \kn$ and $l=\frac{1}{2}\dim_{\R}\kn$.
\end{prop}
 
Let $\varsigma_j$ be the representation of $M$ on 
$\Lambda^j(\kn^*)$. By \eqref{eq4.3.15} and \eqref{eq4.3.16},
we have
\begin{align}
\label{eq4.3.17}
\Xi_{\rho}(\sigma)=\sum_{j=0}^{2l}(-1)^j\Xi_{\varsigma_j,\rho}
(\sigma+(j-l)|\nu|).
\end{align}

On the other hand, since $\dim \kb=1$, from \eqref{eq:introdp}, 
%when it is restricted to $K_M$, 
we have the following identity in $R(K_M,\R)$,
\begin{align}\label{eq:intoMS}
\sum_{i=0}^{\dim \kp}(-1)^{i-1}i\Lambda^i(\kp^*)
=\sum_{j=0}^{\dim \kn}  (-1)^j
\Big(\sum^{\dim \kp_\mathfrak{m}}_{i=0} 
(-1)^i\Lambda^i(\kp_\mathfrak{m}^*)\Big)\otimes
\Lambda^j(\kn^*).
\end{align}

Note that $\sum_{i=0}^{\dim \kp}(-1)^{i-1}i\Lambda^i(\kp^*)
\in R(K,\R)$ is used to define the analytic torsion.
From \eqref{eq4.3.17} and \eqref{eq:intoMS}, 
Shen obtains a very interesting expression for $R_{\rho}(\sigma)$
in term of determinants of shifted Casimir operators,
from which he could obtain equation (\ref{eq4.3.10}).

For a representation $\varsigma$ of $M$ in \eqref{eq4.3.13}, set 
\begin{align}
\label{eq4.3.20a}
r_{\varsigma,\rho}
=\dim_{\C} \Ker (C^{\kg,Z,\widehat{\varsigma}^+,\rho})
-\dim_{\C} \Ker (C^{\kg,Z,\widehat{\varsigma}^-,\rho}).
\end{align}
%From the above proof of Equation (\ref{eq4.3.10}) 
%Theorem \ref{t4.3.3},  for the representation
%$\varsigma_j$ of $M$ on $\Lambda^j(\kn^*)$, we get
Then Shen \cite[(5.12), (7.75)]{Shen16} obtains the formula
\begin{align}
\label{eq4.3.21}
&C_\rho=\prod_{j=0}^{l-1}\left(-4|l-j|^2|\nu|^2\right)^{(-1)^{j-1}
r_{\varsigma_j,\rho}},
&r_{\rho}=2\sum_{j=0}^{l}(-1)^{j-1}r_{\varsigma_j,\rho}. 
\end{align}

Shen shows  that if $H^\bullet(Z,F)=0$, %established 
then $r_{\varsigma_j,\rho}=0$ 
for any $0\le j\le l$ by using 
the spectral aspect of the Selberg trace formula
(Theorem \ref{t0.1} and (\ref{0.5})), 
and some deep results on the representation theory 
of reductive groups.

For a $G$-representation $\pi: G\to {\rm Aut}(V)$ and $v\in V$,
recall that $v$ is said to be differentiable if $c_{v}: G\to V$,
$c_{v}(g)= \pi(g)v$ is $C^\infty$,
that $v$ is said to be $K$-finite if it is contained in a finite 
dimensional subspace stable under $K$. We denote
by $V_{(K)}$ the subspace of differentiable and $K$-finite vectors
in $V$.
Let $H^\bullet (\kg,K;V_{(K)})$ be the $(\kg,K)$-cohomology 
of $V_{(K)}$. 

We denote by $\widehat{G}_u$ the unitary dual of $G$, 
that is the set of equivalence classes of complex irreducible 
unitary representation $\pi$ of $G$ on the Hilbert space $V_\pi$. 
Let $\chi_\pi$ be the corresponding infinitesimal character. 

Let $\widehat{p}:\Gamma\backslash G \to  Z$ be 
the natural projection.
The group $G$ acts unitarily on the right on 
$L_{2}(\Gamma\backslash G, \widehat{p}^*F)$,
then $L_{2}(\Gamma\backslash G, \widehat{p}^*F)$
decomposes into a discrete Hilbert direct sum with 
finite multiplicities of irreducible unitary representations of $G$,
\begin{align}\label{eq:4.3 25}
L_{2}(\Gamma\backslash G, \widehat{p}^*F)
= \widetilde{\bigoplus_{\pi\in \widehat{G}_u}}
n_{\rho}(\pi) V_{\pi} \quad \text{ with } \,
n_{\rho}(\pi) <+\infty,
\end{align}
here $\widetilde{\quad}$ means the Hilbert completion.

Set 
\begin{align}\label{eq:4.3.26}
W= \widetilde{\bigoplus_{\pi\in \widehat{G}_u, \chi_\pi \, \text{Êis trivial} }}
n_{\rho}(\pi) V_{\pi},
\end{align}
then $W$ is  the closure  in 
$L_{2}(\Gamma\backslash G, \widehat{p}^*F)$ of $W^\infty$, 
the subspace of $C^\infty(\Gamma\backslash G,\widehat{p}^*F)$ 
on which the center of $U(\mathfrak{g})$ acts by the same scalar 
as in the trivial representation of $\mathfrak{g}$. 
% infinitesimal character vanishes.
By standard arguments 
\cite[Chap. VII, Theorem 3.2, Corollary 3.4]{BorelWallach00}, 
the cohomology $H^\bullet (Z,F)$
is canonically isomorphic to  the $(\kg,K)$-cohomology
$H^\bullet (\kg,K;W_{(K)})$ of $W_{(K)}$,
the vector space of differentiable and $K$-finite vectors of $W$, i.e., 
\begin{align}\label{eq:4.3.30}
H^\bullet (Z,F) = \bigoplus_{\pi\in \widehat{G}_u, \chi_\pi \, \text{Êis trivial} }
n_{\rho}(\pi)  \, H^\bullet (\kg,K;V_{\pi, (K)}).
\end{align}

 Vogan-Zuckerman \cite[Theorem 1.4]{VoganZuckerman} and 
 Vogan \cite[Theorem 1.3]{Vogan2} classified all irreducible unitary 
 representations of $G$ with nonzero  $(\kg,K)$-cohomology. 
 On the other hand, in \cite[Theorem 1.8]{Salamanca}, 
 Salamanca-Riba showed 
 that any irreducible  unitary representation of $G$ with trivial %vanishing 
 infinitesimal character %$\mathcal{Z}(\kg)$-action 
 is in the class specified by  
 Vogan and Zuckerman, which means that it possesses nonzero 
 $(\kg,K)$-cohomology. In summary,   %summarize,
 if $(\pi,V_\pi)\in \widehat{G}_{u}$, then 
 \begin{align}\label{eq:4.3.31}
 \chi_{\pi} \, \hbox{is non-trivial if and only if } 
 H^\bullet (\kg,K; V_{\pi,(K)})=0.
 \end{align}
 
By the above considerations, %the acyclicity of $F$
$H^\bullet(Z,F)=0$ is equivalent 
to $W=0$.  This is the main algebraic ingredient in 
the proof of (\ref{eq4.3.11}). %Theorem \ref{t4.3.4}. 

Shen's contribution \cite[Corollary 8.15]{Shen16} is to give a formula 
for $r_{\varsigma_j,\rho}$ using  %with the help of 
Hecht-Schmid's work  \cite[Theorem 3.6, Corollary 3.32]{HechtSchmid} 
on the $\kn$-homology 
of $W$, Theorem \ref{t2.7} and (\ref{eq2.4.10}). 
From this formula, 
we see immediately that $W=0$ implies $r_{\varsigma_j,\rho}=0$
for $0\le j\le l$.
\end{proof}

%%%%%%%%%%%%%%%%%%%%%%%%%%%%%%%%%
\subsection{Final remarks}\label{s4.4}
1. Theorem \ref{t2.7} gives an explicit formula for 
the orbital integrals for the heat kernel of the Casimir operator
and it holds for any semisimple element $\gamma\in G$. 
A natural question is how to evaluate or define 
the weighted orbital integrals that appear in
%for a nonsemisimple element, or more geometrically, to find
%the precise %formula of the 
Selberg trace formula for 
a discrete subgroup $\Gamma\subset G$ such that 
$\Gamma\backslash G$ has a finite volume. 

2. Bismut-Goette \cite{BG04} introduced a local topological 
invariant for compact manifolds with a compact Lie group action: 
the $V$-invariant. It appears as an exotic term in the difference 
between two natural versions of equivariant analytic torsion. 
The $V$-invariant shares formally many similarities %property for 
with the analytic torsion, and if we apply formally the construction 
of the $V$-invariant to the associated loop space
equipped with its natural 
$\mathbb{S}^1$ action, then we get the analytic torsion.

In Shen's proof of Fried's conjecture for locally symmetric spaces, 
Shen observed that the $V$-invariant for the $\mathbb{S}^1$-action 
on $B_{[\gamma]}$ is exactly 
\begin{align}\label{eq:4.3.35}
-\frac{\chi_{\rm{orb}}(\mathbb{S}^1\backslash B_{[\gamma]})}
{2n_{[\gamma]}}.
\end{align}
 This suggests a general definition of the Ruelle dynamical 
 zeta function for any compact manifold by replacing 
$-\frac{\chi_{\rm{orb}}(\mathbb{S}^1\backslash B_{[\gamma]})}
{2n_{[\gamma]}}$ in Definition \ref{d4.3.2}
by the associated $V$-invariant. 
One could then compare it with the analytic torsion, 
and obtain a generalized version of Fried's conjecture for 
any manifold with non positive curvature . 
Note that for a strictly negative curvature manifold, $B_{[\gamma]}$ 
is a circle and the $V$-invariant is $-\frac{1}{2n_{[\gamma]}}$.
Recently, Giulietti-Liverani-Pollicott \cite{GiuLP13},
 Dyatlov-Zworski \cite{DyaZ16}, Faure-Tsujii \cite{FaureTsuji17}
 established that for the trivial representation of $\pi_{1}(Z)$,
and when $Z$ has strictly negative curvature,
the Ruelle dynamical zeta function is 
a well-defined meromorphic function on $\C$.

%\newpage
%\bibliographystyle{amsplain}\bibliography{Ma,Bismut2,Others2}
%\comment{
\def\cprime{$'$} \def\cprime{$'$}
\providecommand{\bysame}{\leavevmode\hbox to3em{\hrulefill}\thinspace}
\providecommand{\MR}{\relax\ifhmode\unskip\space\fi MR }
% \MRhref is called by the amsart/book/proc definition of \MR.
\providecommand{\MRhref}[2]{%
  \href{http://www.ams.org/mathscinet-getitem?mr=#1}{#2}
}
\providecommand{\href}[2]{#2}

\end{document}